\documentclass[10pt]{article}
\usepackage{amsthm,amsfonts,amsmath}
\usepackage{amscd,amssymb,latexsym,epsfig}
\title{Floer homology and the heat flow}

\author{Dietmar~A.~Salamon\;\;\;\;\;\; Joa Weber\\\\
        ETH-Z\"urich}

\date{28 September 2004}

\newtheorem{theorem}{Theorem}[section]
\newtheorem{corollary}[theorem]{Corollary}

\newtheorem{lemma}[theorem]{Lemma}
\newtheorem{proposition}[theorem]{Proposition}

\newtheorem{definition}[theorem]{Definition}
\newtheorem{remark}[theorem]{Remark}

\newcommand{\1}{{{\mathchoice {\rm 1\mskip-4mu l} {\rm 1\mskip-4mu l}
{\rm 1\mskip-4.5mu l} {\rm 1\mskip-5mu l}}}}

\newcommand{\C}{{\mathbb{C}}}

\newcommand{\N}{{\mathbb{N}}}

\newcommand{\R}{{\mathbb{R}}}

\newcommand{\Z}{{\mathbb{Z}}}
\newcommand{\Aa}{{\mathcal{A}}}   

\newcommand{\Dd}{{\mathcal{D}}}

\newcommand{\Ff}{{\mathcal{F}}}
\newcommand{\Hh}{{\mathcal{H}}}

\newcommand{\Ll}{{\mathcal{L}}}   
\newcommand{\Mm}{{\mathcal{M}}}   

\newcommand{\Pp}{{\mathcal{P}}}

\newcommand{\Ss}{{\mathcal{S}}}
\newcommand{\Tt}{{\mathcal{T}}}

\newcommand{\Vv}{{\mathcal{V}}}
\newcommand{\Ww}{{\mathcal{W}}}

\newcommand{\im}{{\rm im }}        

\newcommand{\INDEX}{{\rm index}}   
\newcommand{\IND}{{\rm ind}}       
\newcommand{\grad}{{\rm grad\,}}    
\newcommand{\Or}{{\rm Or}}            
\newcommand{\HF}{{\rm HF}}            
\newcommand{\HM}{{\rm HM}}            

\newcommand{\eps}{{\varepsilon}}

\renewcommand{\i}{{\iota}}

\newcommand{\om}{{\omega}}

\newcommand{\Om}{{\Omega}}
\newcommand{\Cinf}{C^{\infty}}

\newcommand{\tu}{\tilde{u}}
\newcommand{\tv}{\tilde{v}}
\newcommand{\tw}{\tilde{w}}
\newcommand{\inner}[2]{\langle #1, #2\rangle}

\def\NABLA#1{{\mathop{\nabla\kern-.5ex\lower1ex\hbox{$#1$}}}}
\def\Nabla#1{\nabla\kern-.5ex{}_{#1}}
\def\Tabla#1{\Tilde\nabla\kern-.5ex{}_{#1}}
\def\abs#1{\mathopen|#1\mathclose|}   
\def\Abs#1{\left|#1\right|}            

\def\Norm#1{\left\|#1\right\|}
\def\NORM#1{{{|\mskip-2.5mu|\mskip-2.5mu|#1|\mskip-2.5mu|\mskip-2.5mu|}}}
\renewcommand{\Tilde}{\widetilde}
\renewcommand{\Hat}{\widehat}

\newcommand{\p}{{\partial}}

\newcommand{\IMP}{\Longrightarrow}

\begin{document}

\maketitle


\begin{abstract} 

We study the heat flow in the loop space of a closed
Riemannian manifold $M$ as an adiabatic limit of the Floer 
equations in the cotangent bundle.  Our main 
application is a proof that the Floer homology of the 
cotangent bundle, for the Hamiltonian function kinetic 
plus potential energy, is naturally isomorphic to the 
homology of the loop space.

\end{abstract}



\section{Introduction}\label{sec:intro}

Let $M$ be a closed Riemannian manifold
and denote by $\Ll M$ the free loop space.  
Consider the classical action functional
$$
\Ss_V(x) = \int_0^1 
\left(\frac12|\dot x(t)|^2 - V(t,x(t))\right)\,dt
$$
for $x:S^1\to M$. Here and throughout we identify $S^1=\R/\Z$
and think of $x\in\Cinf(S^1,M)$ as a smooth function
$x:\R\to M$ which satisfies $x(t+1)=x(t)$. 
The potential is a smooth function $V:S^1\times M\to\R$
and we write $V_t(x):=V(t,x)$. The critical points 
of $\Ss_V$ are the $1$-periodic solutions of the ODE
\begin{equation}\label{eq:geo}
\Nabla{t}\dot x = -\nabla V_t(x),
\end{equation}
where $\nabla V_t$ denotes the gradient and $\Nabla{t}\dot x$ 
denotes the Levi-Civita connection.  Let $\Pp = \Pp(V)$ 
denote the set of $1$-periodic solutions $x:S^1\to M$ 
of~(\ref{eq:geo}). In the case $V=0$ these are the closed 
geodesics. Via the Legendre transformation the solutions 
of~(\ref{eq:geo}) can also be interpreted as the 
critical points of the symplectic action 
$\Aa_V:\Ll T^*M\to\R$ given by 
$$
     \Aa_V(z) 
     = \int_0^1 \biggl(
        \inner{y(t)}{\dot x(t)} - H(t,x(t),y(t))\biggr)\,dt
$$
where $z=(x,y):S^1\to T^*M$ and the 
Hamiltonian $H=H_V:S^1\times T^*M\to\R$ is given by 
\begin{equation}\label{eq:H}
     H(t,x,y) = \frac12|y|^2+V(t,x)
\end{equation}
for $y\in T_x^*M$. A loop $z(t)=(x(t),y(t))$ in $T^*M$ 
is a critical point of $\Aa_V$ iff $x$ is a solution
of~(\ref{eq:geo}) and $y(t)\in T_{x(t)}^*M$ is related to
$\dot x(t)\in T_{x(t)}M$ via the isomorphism
$TM\to T^*M$ induced by the Riemannian metric. 
For such loops $z$ the symplectic action $\Aa_V(z)$ 
agrees with the classical action $\Ss_V(x)$. 

The negative $L^2$ gradient flow of the classical action
gives rise to a Morse-Witten complex which computes 
the homology of the loop space.  For a regular value $a$ 
of $\Ss_V$ we shall denote by $\HM_*^a(\Ll M,\Ss_V)$ the 
homology of the Morse-Witten complex of the functional $\Ss_V$ 
corresponding to the solutions of~(\ref{eq:geo}) with
$\Ss_V(x)\le a$.  Here we assume that $\Ss_V$ is a Morse function 
and its gradient flow satisfies the Morse-Smale condition
(i.e. the stable and unstable manifolds intersect 
transversally, see~\cite{DAVIES} for the unstable manifold).  
As in the finite dimensional case
one can show that the Morse-Witten homology
$\HM^a_*(\Ll M,\Ss_V)$ is naturally isomorphic to 
the singular homology of the sublevel set
$$
    \Ll^aM = \left\{x\in\Ll M\,|\,\Ss_V(x)\le a\right\}.
$$
On the other hand one can use the $L^2$ gradient flow 
of $\Aa_V$ to construct Floer homology groups $\HF_*^a(T^*M,H_V)$.
Our main result is the following.

\begin{theorem}\label{thm:main}
Asssume $\Ss_V$ is Morse and $a$ is either a regular
value of $\Ss_V$ or is equal to infinity. 
Then there is a natural isomorphism
$$
\HF^a_*(T^*M,H_V;R)\cong\HM^a_*(\Ll M,\Ss_V;R)
$$
for every principal ideal domain $R$. 
If $M$ is not simply connected then there is a separate
isomorphism for each component of the loop space. 
The isomorphism commutes with the homomorphisms
$\HF^a_*(T^*M,H_V)\to\HF^b_*(T^*M,H_V)$
and $\HM^a_*(\Ll M,\Ss_V)\to\HM^b_*(\Ll M,\Ss_V)$ 
for $a<b$. 
\end{theorem}

\begin{corollary}\label{cor:main}
Let $\Ss_V$ and $a$ be as in Theorem~\ref{thm:main}.
Then there is a natural isomorphism
$$
\HF^a_*(T^*M,H_V;R)\cong\mathrm{H}_*(\Ll^aM;R)
$$
for every principal ideal domain $R$. 
If $M$ is not simply connected then there is a separate
isomorphism for each component of the loop space. 
The isomorphism commutes with the homomorphisms
$\HF^a_*(T^*M,H_V)\to\HF^b_*(T^*M,H_V)$
and $\mathrm{H}_*(\Ll^aM)\to\mathrm{H}_*(\Ll^bM)$ for $a<b$. 
\end{corollary}

\begin{proof}
Theorem~\ref{thm:main} and Theorem~\ref{thm:morse}
\end{proof}

Both the Morse-Witten homology $HM^a_*(\Ll M,\Ss_V)$
and the Floer homology $\HF_*^a(T^*M,H_V)$ are based on the 
same chain complex 
$
    C^a_*
$
which is generated by the solutions of~(\ref{eq:geo}) 
and graded by the Morse index (as critical points of $\Ss_V$).  
In~\cite{JOA1} it is shown that this Morse index agrees,  
up to a universal additive constant zero or
one, with minus the Conley-Zehnder index.   
Thus it remains to compare the boundary operators
and this will be done by considering an adiabatic limit
with a family of metrics on $T^*M$ which scales the 
vertical part down to zero. 
Another approach to Corollary~\ref{cor:main} 
is contained in Viterbo's paper~\cite{VITERBO}. 
Some recent applications of Corollary~\ref{cor:main} 
can be found in~\cite{WEBER1}; these applications
require the statement with action windows 
and fixed homotopy classes of loops.

\subsection*{The Floer chain complex and its adiabatic limit}

We assume throughout that $\Ss_V$ is a Morse function 
on the loop space, i.e. that the $1$-periodic solutions 
of~(\ref{eq:geo}) are all nondegenerate.
(For a proof that this holds for a generic potential 
$V$ see~\cite{JOA1}.)  Under this assumption the set
$$
    \Pp^a(V) := \left\{x\in\Pp(V)\,|\,\Ss_V(x)\le a\right\}
$$
is finite for every real number $a$.  
Moreover, each critical point $x\in\Pp(V)$ has 
well defined stable and unstable manifolds
with respect to the (negative) $L^2$ gradient 
flow (see for example Davies~\cite{DAVIES}).
Call $\Ss_V$ {\bf Morse--Smale} if it is a Morse function
and the unstable manifold $W^u(y)$ intersects the stable 
manifold $W^s(x)$ transversally for any two critical points
$x,y\in\Pp(V)$.  

Assume $\Ss_V$ is a Morse function and
consider the $\Z$-module
$$
     C^a =  
     C^a(V)
     = \bigoplus_{x\in\Pp^a(V)} \Z x.
$$
If $\Ss_V$ and $\Aa_V$ are Morse--Smale then
this module carries two boundary operators. The first 
is defined by counting the (negative) gradient flow lines of 
$\Ss_V$.  They are solutions $u:\R\times S^1\to M$
of the heat equation 
\begin{equation}\label{eq:heat}
     \p_su - \Nabla{t}\p_tu - \nabla V_t(u) = 0
\end{equation}
satisfying 
\begin{equation}\label{eq:heat-lim}
     \lim_{s\to\pm\infty} u(s,t) = x^\pm(t),\qquad
     \lim_{s\to\pm\infty}\p_su = 0,
\end{equation}
where $x^\pm\in\Pp(V)$.  The limits are uniform in $t$. 
The space of solutions of~(\ref{eq:heat}) and~(\ref{eq:heat-lim})
will be denoted by $\Mm^0(x^-,x^+;V)$. 
The Morse--Smale hypothesis guarantees that, for every
pair $x^\pm\in\Pp^a(V)$, the space $\Mm^0(x^-,x^+;V)$ is a 
smooth manifold whose dimension is equal to the difference
of the Morse indices. In the case of Morse index difference 
one it follows that the quotient $\Mm^0(x^-,x^+;V)/\R$
by the (free) time shift action is a finite set. 
Counting the number of solutions with appropriate signs 
gives rise to a boundary operator on $C^a(V)$.  
The homology $\HM^a_*(\Ll M,\Ss_V)$
of the resulting chain complex is naturally isomorphic 
to the singular homology of the loop space for every
regular value $a$ of $\Ss_V$:
$$
    \HM^a_*(\Ll M,\Ss_V) \cong\mathrm{H}_*(\Ll^aM;\Z),\qquad
    \Ll^aM := \left\{x\in\Ll M\,|\,\Ss_V(x)\le a\right\}.
$$
The details of this isomorphism will be established 
in a separate paper (see Appendix~\ref{app:MS} for 
a summary of the relevant results).
 
The second boundary operator is defined by counting the 
negative gradient flow lines of the symplectic
action functional $\Aa_V$.  These are the solutions 
$(u,v):\R\times S^1\to TM$ of the {\bf Floer equations} 
\begin{equation}\label{eq:floer}
\p_su - \Nabla{t}v - \nabla V_t(u) = 0,\qquad
\Nabla{s}v + \p_tu - v = 0,
\end{equation}
\begin{equation}\label{eq:floer-lim}
\lim_{s\to\pm\infty}u(s,t)=x^\pm(t),\qquad
\lim_{s\to\pm\infty}v(s,t)=\dot x^\pm(t).
\end{equation}
Here we also assume that $\p_su$ and $\Nabla{s}v$ 
converge to zero, uniformly in $t$, as $|s|$ tends to infinity. 
For notational simplicity we identify the tangent 
and cotangent bundles of $M$ via the metric.   
Counting the index-$1$ solutions of~(\ref{eq:floer})
and~(\ref{eq:floer-lim}) with appropriate signs
we obtain the Floer boundary operator.  
We wish to prove that the resulting 
Floer homology groups $\HF^a_*(T^*M,H_V)$
are naturally isomorphic to $\HM^a_*(\Ll M,\Ss_V)$. 
To construct this isomorphism we modify 
equation~(\ref{eq:floer}) by introducing a small 
parameter $\eps$ as follows
\begin{equation}\label{eq:floer-eps}
\p_su - \Nabla{t}v - \nabla V(t,u) = 0,\qquad
\Nabla{s}v + \eps^{-2}(\p_tu - v) = 0.
\end{equation}
The space of solutions of~(\ref{eq:floer-eps}) 
and~(\ref{eq:floer-lim}) will be denoted by $\Mm^\eps(x^-,x^+;V)$.
The Floer homology groups for different values 
of $\eps$ are isomorphic (see Remark~\ref{rmk:eps}
below). Thus the task at hand is to prove that, 
for $\eps>0$ sufficiently small, there is a one-to-one
correspondence between the solutions of~(\ref{eq:heat})
and those of~(\ref{eq:floer-eps}). 
A first indication, why one might expect such a 
correspondence, is the energy identity
\begin{equation*}
\begin{split}
  E^\eps(u,v)
   &=
    \frac12\int_{-\infty}^\infty\int_0^1\left(
    |\p_su|^2 + |\Nabla{t}v+\nabla V_t(u)|^2
    + \eps^2|\Nabla{s}v|^2
    + \eps^{-2}|\p_tu-v|^2
    \right) \\
   &=
   \Ss_V(x^-) - \Ss_V(x^+)
\end{split}
\end{equation*}
for the solutions of~(\ref{eq:floer-eps})
and~(\ref{eq:floer-lim}).
It shows that $\p_tu-v$ must converge to zero in the 
$L^2$ norm as $\eps$ tends to zero.  If $\p_tu=v$ then 
the first equation in~(\ref{eq:floer-eps}) is equivalent
to~(\ref{eq:heat}). 

\begin{remark}\label{rmk:eps}\rm
Let $M$ be a Riemannian manifold.  Then the tangent 
space of the cotangent bundle $T^*M$ at a point
$(x,y)$ with $y\in T_x^*M$ can be identified
with the direct sum $T_xM\oplus T_x^*M$.
The isomorphism takes the derivative $\dot z(t)$
of a curve $\R\to T^*M:t\mapsto z(t)=(x(t),y(t))$
to the pair $(\dot x(t),\Nabla{t}y(t))$.  With this 
identification the almost complex structure 
$J_\eps$ and the metric $G_\eps$ on $T^*M$,
given by
\begin{equation*}
     J_\eps = \begin{pmatrix}
     0 & -\eps g^{-1} \\ \eps^{-1}g & 0
     \end{pmatrix},\qquad
     G_\eps = \begin{pmatrix}
     \eps^{-1}g & 0 \\ 0 & \eps g^{-1}
     \end{pmatrix},
\end{equation*}
are compatible with the standard symplectic form
$\om$ on $T^*M$.  Here we denote by $g:TM\to T^*M$
the isomorphism induced by the metric.  
The case $\eps=1$ corresponds to the standard almost
complex structure.  The Floer equations for the almost 
complex structure $J_\eps$ and the Hamiltonian~(\ref{eq:H})
are
$$
\p_sw-J_\eps(w)(\p_tw-X_{H_t}(w)) = 0.
$$
If we write $w(s,t)=(u(s,t),v(s,t))$
with $v(s,t)\in T_{u(s,t)}^*M$
then this equation has the form 
\begin{equation}\label{eq:floerT*M}
\p_su-\eps g^{-1}\Nabla{t}v - \eps\nabla V_t(u) = 0,\qquad
\Nabla{s}v + \eps^{-1}g\p_tu - \eps^{-1}v = 0. 
\end{equation}
A function $w=(u,v)$ is a solution of~(\ref{eq:floerT*M})
if and only if the functions $\tu(s,t):=u(\eps^{-1}s,t)$
and $\tv(s,t):=g^{-1}v(\eps^{-1}s,t)$ satisfy~(\ref{eq:floer-eps}).  
In view of this discussion it follows from the Floer homotopy argument 
that the Floer homology defined with the solutions of~(\ref{eq:floer-eps})
is independent of the choice of $\eps>0$. 
\end{remark}

Assume $\Ss_V$ is Morse--Smale.  Then we shall prove that, 
for every $a\in\R$, there exists an $\eps_0>0$ such that,
for $0<\eps<\eps_0$ and every pair 
$x^+,x^-\in\Pp^a(V)$ with Morse index difference one,
there is a natural bijective correspondence between
the (shift equivalence classes of) solutions 
of~(\ref{eq:heat}), (\ref{eq:heat-lim})
and those of~(\ref{eq:floer-eps}), (\ref{eq:floer-lim}). 
This will follow from Theorems~\ref{thm:existence}
and~\ref{thm:onto} below. 

It is an open question if the function $\Ss_V$ is Morse--Smale 
(with respect to the $L^2$ metric on the loop space) 
for a generic potential $V$. However, it is easy to establish 
transversality for a general class of abstract perturbations
$\Vv:\Ll M\to\R$ (see Section~\ref{sec:perturbations}). 
We shall use these perturbations to prove 
Theorem~\ref{thm:main} in general. 

The general outline of the proof is similar to that
of the Atiyah--Floer conjecture in~\cite{DOSA} which
compares two elliptic PDEs via an adiabatic limit argument. 
By contrast our adiabatic limit theorem compares elliptic with 
parabolic equations.  This leads to new features in the
analysis that are related to the fact that the parabolic 
equation requires different scaling in space and time 
directions.

The present paper is organized as follows.  
The next section introduces a relevant class of abstract 
perturbations $\Vv:\Ll M\to\R$.
Section~\ref{sec:linear}  explains the relevant linearized 
operators and states the estimates for the right inverse.  
These are proved in 
Appendices~\ref{app:Lp} and~\ref{app:inverse}. 
In Section~\ref{sec:adiabatic} we construct a map
$\Tt^\eps:\Mm^0(x^-,x^+;\Vv)\to\Mm^\eps(x^-,x^+;\Vv)$ 
which assigns to every parabolic cylinder of index one
a nearby Floer connecting orbit for $\eps>0$ sufficiently
small. The existence of this map was established 
in the thesis of the second author~\cite{JOA}, 
where the results of Section~\ref{sec:linear},
Section~\ref{sec:adiabatic}, and Appendix~\ref{app:inverse}
were proved. Sections~\ref{sec:apriori}, 
\ref{sec:gradient}, and~\ref{sec:2nd-derivs}
are of preparatory nature and establish uniform estimates
for the solutions of~(\ref{eq:floer-eps}). 
Section~\ref{sec:exp-decay} deals with exponential decay,
Section~\ref{sec:time-shift} establishes local 
surjectivity of the map $\Tt^\eps$ by a time-shift 
argument, and in Section~\ref{sec:onto}
we prove that $\Tt^\eps$ is bijective. 
Things are put together in Section~\ref{sec:proof}
where we compare orientations and prove Theorem~\ref{thm:main}.
Appendix~\ref{app:MS} summarizes some results about 
the heat flow~(\ref{eq:heat}) which will be proved
in~\cite{WEBER}.  In Appendix~\ref{app:mean} we prove several
mean value inequalities that play a central role in our apriori
estimates of Sections~\ref{sec:apriori}, \ref{sec:gradient}, 
and~\ref{sec:2nd-derivs}.  


\section{Perturbations}\label{sec:perturbations}

In this section we introduce a class of perturbations of 
equations~(\ref{eq:heat}) and~(\ref{eq:floer-eps}) for which
transversality is easy to achieve.  The perturbations take the 
form of smooth maps 
$
\Vv:\Ll M\to\R.
$
For $x\in\Ll M$ let $\grad\Vv(x)\in\Om^0(S^1,x^*TM)$
denote the $L^2$-gradient of $\Vv$; it is defined by 
$$
\int_0^1\inner{\grad\Vv(u)}{\p_su}\,dt := \frac{d}{ds}\Vv(u)
$$
for every smooth path $\R\to\Ll M:s\mapsto u(s,\cdot)$. 
The {\bf covariant Hessian} of $\Vv$ at a loop $x:S^1\to M$
is the operator
$
\Hh_\Vv(x):\Om^0(S^1,x^*TM)\to\Om^0(S^1,x^*TM)
$
defined by 
$$
\Hh_\Vv(u)\p_su := \Nabla{s}\grad\Vv(u)
$$
for every smooth map $\R\to\Ll M:s\mapsto u(s,\cdot)$.
The axiom~$(V1)$ below asserts that this Hessian is a zeroth 
order operator. We impose the following conditions on $\Vv$;
here $\Abs{\cdot}$ denotes the pointwise absolute value 
at $(s,t)\in\R\times S^1$ and $\Norm{\cdot}_{L^p}$ denotes 
the $L^p$-norm over $S^1$ at time $s$. 

\smallbreak

\begin{description}
\item[(V0)]
$\Vv$ is continuous with respect to the $C^0$ topology
on $\Ll M$.  Moreover, there is a constant $C>0$ such that
$$
\sup_{x\in\Ll M}\left|\Vv(x)\right|
+\sup_{x\in\Ll M}\left\|\grad\Vv(x)\right\|_{L^\infty(S^1)}
\le C.
$$
\item[(V1)]
There is a constant $C>0$ such that
\begin{align*}
\Abs{\Nabla{s}\grad\Vv(u)}
&\le C\bigl(\left|\p_su\right|+\Norm{\p_su}_{L^1}\bigr), \\
\left|\Nabla{t}\grad\Vv(u)\right|
&\le C\Bigl(1+\left|\p_tu\right|\Bigr)
\end{align*}
for every smooth map $\R\to\Ll M:s\mapsto u(s,\cdot)$ 
and every $(s,t)\in\R\times S^1$. 
\item[(V2)]
There is a constant $C>0$ such that
\begin{align*}
\Abs{\Nabla{s}\Nabla{s}\grad\Vv(u)} 
&\le C\Bigl(\Abs{\Nabla{s}\p_su} 
+ \Norm{\Nabla{s}\p_su}_{L^1} 
+ \bigl(\Abs{\p_su} + \Norm{\p_su}_{L^2}\bigr)^2
\Bigr), \\
\Abs{\Nabla{t}\Nabla{s}\grad\Vv(u)} 
&\le C\Bigl(
\Abs{\Nabla{t}\p_su} 
+ \bigl(1+\Abs{\p_tu}\bigr)
\bigl(\Abs{\p_su} + \Norm{\p_su}_{L^1}\bigr)
\Bigr)
\end{align*}
and
$$
\Abs{\Nabla{s}\Nabla{s}\grad\Vv(u)
- \Hh_\Vv(u)\Nabla{s}\p_su}
\le C\bigl(\Abs{\p_su} + \Norm{\p_su}_{L^2}\bigr)^2
$$
for every smooth map $\R\to\Ll M:s\mapsto u(s,\cdot)$ 
and every $(s,t)\in\R\times S^1$.
\item[(V3)]
There is a constant $C>0$ such that
\begin{align*}
\Abs{\Nabla{s}\Nabla{s}\Nabla{s}\grad\Vv(u)} 
&\le C\Bigl(\Abs{\Nabla{s}\Nabla{s}\p_su} 
+ \Norm{\Nabla{s}\Nabla{s}\p_su}_{L^1} \\
&\quad
+ \bigl(\Abs{\Nabla{s}\p_su} + \Norm{\Nabla{s}\p_su}_{L^2}\bigr)
\bigl(\Abs{\p_su} + \Norm{\p_su}_{L^2}\bigr)\\
&\quad
+ \bigl(\Abs{\p_su} + \Norm{\p_su}_{L^\infty}\bigr)
\bigl(\Abs{\p_su} + \Norm{\p_su}_{L^2}\bigr)^2
\Bigr), \\
\Abs{\Nabla{t}\Nabla{s}\Nabla{s}\grad\Vv(u)} 
&\le C\Bigl(
\Abs{\Nabla{t}\Nabla{s}\p_su} 
+ \Abs{\Nabla{t}\p_su}
\bigl(\Abs{\p_su} + \Norm{\p_su}_{L^1}\bigr) \\
&\quad
+ \bigl(1+\Abs{\p_tu}\bigr)
\bigl(\Abs{\Nabla{s}\p_su} + \Norm{\Nabla{s}\p_su}_{L^1}\bigr) \\
&\quad
+ \bigl(1+\Abs{\p_tu}\bigr)\bigl(\Abs{\p_su} + \Norm{\p_su}_{L^2}\bigr)^2 
\Bigr), \\
\Abs{\Nabla{t}\Nabla{t}\Nabla{s}\grad\Vv(u)} 
&\le C\Bigl(
\Abs{\Nabla{t}\Nabla{t}\p_su}
+ \bigl(1+\Abs{\p_tu}\bigr)\Abs{\Nabla{t}\p_su} \\
&\quad
+ \bigl(1+\Abs{\p_tu}^2+\Abs{\Nabla{t}\p_tu}\bigr)
\bigl(\Abs{\p_su} + \Norm{\p_su}_{L^1}\bigr)
\Bigr)
\end{align*}
for every smooth map $\R\to\Ll M:s\mapsto u(s,\cdot)$ 
and every $(s,t)\in\R\times S^1$.
\item[(V4)]
For any two integers $k>0$ and $\ell\ge 0$
there is a constant $C=C(k,\ell)$ such that
$$
\Abs{\nabla_t^\ell\nabla_s^k\grad\Vv(u)}
\le C\sum_{k_j,\ell_j}
\left(
\prod_{j\atop\ell_j>0}
\Abs{\nabla_t^{\ell_j}\nabla_s^{k_j}u}
\right)
\prod_{j\atop\ell_j=0}
\Biggl(\Abs{\nabla_s^{k_j}u}
+\left\|\nabla_s^{k_j}u\right\|_{L^{p_j}}
\Biggr)
$$
for every smooth map $\R\to\Ll M:s\mapsto u(s,\cdot)$ 
and every $(s,t)\in\R\times S^1$;
here $p_j\ge 1$ and $\sum_{\ell_j=0}1/p_j=1$; 
the sum runs over all partitions $k_1+\cdots+k_m=k$
and $\ell_1+\cdots+\ell_m\le\ell$ such that $k_j+\ell_j\ge1$ 
for all $j$. For $k=0$ the same inequality holds
with an additional summand $C$ on the right. 
\end{description}

\begin{remark}\rm
The archetypal example of a perturbation is
$$
\Vv(x):= \rho\left(\left\|x-x_0\right\|_{L^2}^2\right)
\int_0^1V_t(x(t))\,dt,
$$
where $\rho:\R\to[0,1]$ is a smooth cutoff function,
$x_0:S^1\to M$ is a smooth loop, and $x-x_0$ denotes the 
difference in some ambient Euclidean space into which 
$M$ is (isometrically) embedded.  Any such perturbation 
satisfies~$(V0-V4)$.  
\end{remark}

\begin{remark}\label{rmk:HV}\rm
If 
$$
\Vv(x)=\int_0^1V_t(x(t))\,dt
$$ 
then 
$$
\grad\Vv(x) = \nabla V_t(x),\qquad
\Hh_\Vv(x)\xi = \Nabla{\xi}\nabla V_t(x),
$$
for $x\in\Ll M$ and  $\xi\in\Om^0(S^1,x^*TM)$. 
\end{remark}

With an abstract perturbation 
$\Vv$ the classical and symplectic action are given by
$$
\Ss_\Vv(x) = \frac12\int_0^1\Abs{\dot x(t)}^2\,dt - \Vv(x)
$$
and
$$
\Aa_\Vv(x,y) = \int_0^1\left(\inner{y(t)}{\dot x(t)}
-\frac12\Abs{y(t)}^2\right)\,dt - \Vv(x)
$$
for $x\in\Ll M$ and $y\in\Om^0(S^1,x^*T^*M)$. 
Equation~(\ref{eq:floer-eps}) has the form
\begin{equation}\label{eq:floer-V}
   \p_su - \Nabla{t}v - \grad\Vv(u) = 0,\qquad
   \Nabla{s}v + \eps^{-2}(\p_tu - v) = 0.
\end{equation}
and the limit equation is
\begin{equation}\label{eq:heat-V}
   \p_su - \Nabla{t}\p_tu - \grad\Vv(u) = 0.
\end{equation}
Here $\grad\Vv(u)$ denotes the value of $\grad\Vv$ on the
loop $t\mapsto u(s,t)$.  The relevant set of critical points 
consists of the loops $x:S^1\to M$ that satisfy the 
differential equation $\Nabla{t}\dot x=\grad\Vv(x)$ and will be
denoted by $\Pp(\Vv)$. The subset $\Pp^a(\Vv)\subset\Pp(\Vv)$ 
consists of all critical points $x$ with $\Ss_\Vv(x)\le a$.


\section{The linearized operators}\label{sec:linear}

Throughout this section we fix a perturbation
$\Vv$ that satisfies~$(V0-V4)$.
Linearizing the heat equation~(\ref{eq:heat-V})
gives rise to the operator
$$
    \Dd^0_u:\Om^0(\R\times S^1,u^*TM)\to\Om^0(\R\times S^1,u^*TM)
$$
given by
\begin{equation}\label{eq:D0}
     \Dd^0_u\xi 
     = \Nabla{s}\xi
       - \Nabla{t}\Nabla{t}\xi 
       - R(\xi,\p_tu)\p_tu
       - \Hh_\Vv(u)\xi,
\end{equation}
for every element $\xi$
of the set $\Om^0(\R\times S^1,u^*TM)$
of smooth vector fields along~$u$.
If $\Ss_\Vv$ is Morse then
this is a Fredholm operator between 
appropriate Sobolev completions.
More precisely, define 
$$
\Ll_u=\Ll_u^p,\qquad \Ww_u = \Ww_u^p
$$
as the completions of the space of smooth 
compactly supported sections of the pullback
tangent bundle $u^*TM\to\R\times S^1$
with respect to the norms
$$
     \left\|\xi\right\|_{\Ll}
     = \left(\int_{-\infty}^\infty\int_0^1 
     |\xi|^p\,dtds\right)^{1/p},
$$
$$
     \left\|\xi\right\|_{\Ww}
     = \left(\int_{-\infty}^\infty\int_0^1 
      |\xi|^p + |\Nabla{s}\xi|^p 
      + |\Nabla{t}\Nabla{t}\xi|^p
      \,dtds\right)^{1/p}.
$$
Then $\Dd^0_u:\Ww_u^p\to\Ll_u^p$ is a Fredholm operator
for $p>1$ (Theorem~\ref{thm:par-Fredholm}) with index
$$
     \INDEX\Dd^0_u = \IND_\Vv(x^-)
     -\IND_\Vv(x^+).
$$
Here $\IND_\Vv(x)$ denotes the Morse index,
i.e. the number of negative eigenvalues of the Hessian
of $\Ss_\Vv$.  This Hessian is given by
$$
     A^0(x)\xi = -\Nabla{t}\Nabla{t}\xi 
                - R(\xi,\dot x)\dot x
                - \Hh_\Vv(x)\xi,
$$
where $R$ denotes the Riemann curvature tensor
and $\Hh_\Vv$ denotes the covariant Hessian of $\Vv$
(see Section~\ref{sec:perturbations}). 
The Morse--Smale condition asserts that the operator 
$\Dd^0_u$ is surjective for every finite energy solution
of~(\ref{eq:heat}).  That this condition can
be achieved by a generic perturbation $\Vv$
is proved in~\cite{WEBER} (see Appendix~\ref{app:MS}).

Linearizing equation~(\ref{eq:floer-V})
gives rise to the first order differential operator
$$
    \Dd^\eps_{u,v}:
     W^{1,p}(\R\times S^1,u^*TM\oplus u^*TM)
     \to L^p(\R\times S^1,u^*TM\oplus u^*TM)
$$
given by 
\begin{equation}\label{eq:D-eps}
    \Dd^\eps_{u,v}\begin{pmatrix}
     \xi \\ \eta \end{pmatrix}
     = \begin{pmatrix}
       \Nabla{s}\xi -\Nabla{t}\eta - R(\xi,\p_tu)v - 
       \Hh_\Vv(u)\xi \\
       \Nabla{s}\eta + R(\xi,\p_su)v
       + \eps^{-2}(\Nabla{t}\xi-\eta)
       \end{pmatrix}
\end{equation}
for $(\xi,\eta)\in W^{1,p}(\R\times S^1,u^*TM\oplus u^*TM)$.

\begin{remark}\label{rmk:eps2}\rm
Assume $\Ss_\Vv$ is Morse and let $p>1$.
Then $\Dd^\eps_{u,v}$ is a Fredholm operator
for every pair $(u,v)$ that satisfies~(\ref{eq:floer-lim}) 
and its index is given by 
$$
      \INDEX\,\Dd^\eps_{u,v} 
      = \IND_\Vv(x^-)-\IND_\Vv(x^+).
$$
To see this rescale $u$ and $v$ as in
Remark~\ref{rmk:eps}.
Then the operator on the rescaled vector fields
$\tilde\xi(s,t):=\xi(\eps^{-1}s,t)$
and $\tilde\eta(s,t):=g^{-1}\eta(\eps^{-1}s,t)$
has the same form as in Floer's original
papers~\cite{FLOER5}
with the almost complex structure $J_\eps$ of
Remark~\ref{rmk:eps}.
That this operator is Fredholm was proved
in~\cite{FLOER2,SALZ,ROSA}
for $p=2$. An elegant proof of the
Fredholm property
for general $p>1$ was given by
Donaldson~\cite{DON} 
for the instanton case; it adapts easily to the
symplectic
case~\cite{Sa97}. The Fredholm index can be
expressed as 
a difference of the Conley--Zehnder
indices~\cite{SALZ,ROSA}.
That it agrees with the difference of the
Morse indices was 
proved in~\cite{JOA}.
\end{remark}

Let us now fix a solution $u$ of~(\ref{eq:heat})
and define
$
      v:=\p_tu.
$
For this pair $(u,v)$ we must prove that 
the operator $\Dd_u^\eps:=\Dd_{u,\p_tu}^\eps$ 
is onto for $\eps>0$
sufficiently small and prove an estimate
for the right
inverse which is independent of $\eps$.  
We will establish this under the
assumption that the 
operator $\Dd_u^0$ is onto.
To obtain uniform estimates
for the inverse with constants
independent of $\eps$
we must work with suitable $\eps$-dependent norms.
For compactly supported vector fields
$\zeta=(\xi,\eta)\in\Om^0(\R\times S^1,u^*TM\oplus u^*TM)$ 
define 
$$
     \left\|\zeta\right\|_{0,p,\eps}
     = \left(\int_{-\infty}^\infty\int_0^1
       \left(\left|\xi\right|^p 
       + \eps^p\left|\eta\right|^p\right)\,dtds
       \right)^{1/p},
$$
\begin{equation*}
\begin{split}
     \left\|\zeta\right\|_{1,p,\eps}
    &=
     \biggl(\int_{-\infty}^\infty\int_0^1
       \bigl(\left|\xi\right|^p 
       + \eps^p\left|\eta\right|^p
       + \eps^p\left|\Nabla{t}\xi\right|^p 
       + \eps^{2p}\left|\Nabla{t}\eta\right|^p \\
    &\quad
       +\, \eps^{2p}\left|\Nabla{s}\xi\right|^p 
       + \eps^{3p}\left|\Nabla{s}\eta\right|^p
       \bigr)\,dtds\biggr)^{1/p}.
\end{split}
\end{equation*}

\begin{theorem}\label{thm:elliptic-eps}
Let $(u,v):\R\times S^1\to TM$ be a smooth map
such that $v$ and the derivatives 
$\p_su,\p_tu,\Nabla{t}\p_su,\Nabla{t}\p_tu$ are bounded
and $\lim_{s\to\pm\infty}u(s,t)$ exists, uniformly in $t$.
Then, for every $p>1$,
there are positive constants $c$ and 
$\eps_0$ such that, for every $\eps\in(0,\eps_0)$
and every
$
     \zeta=(\xi,\eta)\in W^{1,p}(\R\times S^1,
     u^*TM\oplus u^*TM),
$
we have
\begin{equation}\label{eq:elliptic-eps}
\begin{split}
     &\eps^{-1}\left\|\Nabla{t}\xi-\eta\right\|_{L^p}
      + \left\|\Nabla{t}\eta\right\|_{L^p}
      + \left\|\Nabla{s}\xi\right\|_{L^p}
      + \eps\left\|\Nabla{s}\eta\right\|_{L^p} \\
     &\le 
      c\left(
      \left\|\Dd_{u,v}^\eps\zeta\right\|_{0,p,\eps}
      + \left\|\xi\right\|_{L^p}
      + \eps^2\left\|\eta\right\|_{L^p}
      \right).
\end{split}
\end{equation}
The formal adjoint operator $(\Dd_{u,v}^\eps)^*$
defined below 
satisfies the same estimate.
Moreover, the constants $c$ and $\eps_0$
are invariant under $s$-shifts of $u$.
\end{theorem}

The formal adjoint operator 
$$
     (\Dd^\eps_{u,v})^*:
     W^{2,p}(\R\times S^1,u^*TM\oplus u^*TM)
     \to W^{1,p}(\R\times S^1,u^*TM\oplus u^*TM)
$$
with respect to the $(0,2,\eps)$-inner
product associated 
to the $(0,2,\eps)$-norm has the form
$$
    (\Dd^\eps_{u,v})^*\begin{pmatrix}
     \xi \\ \eta \end{pmatrix}
     = \begin{pmatrix}
       -\Nabla{s}\xi -\Nabla{t}\eta - R(\xi,v)\p_tu 
       - \Hh_\Vv(u)\xi
       + \eps^2R(\eta,v)\p_su \\
       -\Nabla{s}\eta 
       + \eps^{-2}(\Nabla{t}\xi-\eta)
       \end{pmatrix}
$$
for $\xi,\eta\in W^{1,p}(\R\times S^1,u^*TM)$.
We shall also use the projection operator 
$$
     \pi_\eps:L^p(S^1,x^*TM)\times L^p(S^1,x^*TM)\to 
              W^{1,p}(S^1,x^*TM)
$$
given by 
$$
     \pi_\eps(\xi,\eta) 
     = (\1-\eps\Nabla{t}\Nabla{t})^{-1}
       (\xi-\eps^2\Nabla{t}\eta)
$$
for $x\in\Ll M$ and $\xi,\eta\in\Om^0(S^1,x^*TM)$.
This operator, for the loop $x(t)=u(s,t)$, will be applied 
to the pair $(\xi(s,\cdot),\eta(s,\cdot))$. 

\begin{theorem}\label{thm:inverse}
Assume $\Ss_\Vv$ is Morse-Smale and
let $u\in\Mm^0(x^-,x^+;\Vv)$. Then, for every $p>1$, 
there are positive constants $c$ and 
$\eps_0$ (invariant under $s$-shifts of $u$)
such that, for every 
$\eps\in(0,\eps_0)$ the following are true.
The operator $\Dd_u^\eps:=\Dd^\eps_{u,\p_tu}$
is onto and for every pair
$$
\zeta:=(\xi,\eta)\in \im\,(\Dd^\eps_u)^*
\subset W^{1,p}(\R\times S^1,u^*TM\oplus u^*TM)
$$
we have
\begin{equation}\label{eq:inverse-xieta}
      \left\|\xi\right\|_{L^p}
      + \eps^{1/2}\left\|\eta\right\|_{L^p}
      + \eps^{1/2}\left\|\Nabla{t}\xi\right\|_{L^p}
      \le c\left(
       \eps\left\|\Dd_u^\eps\zeta\right\|_{0,p,\eps}
       + \left\|\pi_\eps(\Dd_u^\eps\zeta)\right\|_{L^p}
       \right),
\end{equation}
\begin{equation}\label{eq:inverse-zeta}
      \left\|\zeta\right\|_{1,p,\eps}
      \le c\left(
       \eps\left\|\Dd_u^\eps\zeta\right\|_{0,p,\eps}
       + \left\|\pi_\eps(\Dd_u^\eps\zeta)\right\|_{L^p}
       \right). 
\end{equation}
\end{theorem}

The proofs of Theorems~\ref{thm:elliptic-eps} 
and~\ref{thm:inverse} are given in
Appendix~\ref{app:inverse}. 
They are based on a simplified form of 
Theorem~\ref{thm:elliptic-eps} for flat manifolds
with $\Vv=0$ which is proved in Appendix~\ref{app:Lp}.
In particular, Corollary~\ref{cor:CZ-eps} shows that the 
$\eps$-weights on the left hand side of 
equation~(\ref{eq:elliptic-eps}) appear
in a natural manner by a rescaling argument
and, for $p=2$, these terms can be interpreted as a 
linearized version of the energy.  
This was in fact the motivation for 
introducing the above $\eps$-dependent norms. 
The proof of Theorem~\ref{thm:inverse} is based on 
Theorem~\ref{thm:elliptic-eps} and a comparison 
of the operators $\Dd^0_u$ and $\Dd^\eps_u$. 

To construct a solution of~(\ref{eq:floer-eps})
near a parabolic cylinder it is useful to combine
Theorems~\ref{thm:elliptic-eps} and~\ref{thm:inverse}
into the following corollary.  This corollary involves
an $\eps$-dependent norm which at first glance appears
to be somewhat less natural but plays a useful role
for technical reasons. 

Given a smooth map
$u:\R\times S^1\to M$ and a compactly 
supported pair of vector fields 
$\zeta=(\xi,\eta)\in\Om^0(\R\times S^1,u^*TM\oplus u^*TM)$
we define 
\begin{equation}\label{eq:NORM}
\begin{split}
\NORM{\zeta}_\eps 
&:= \Norm{\xi}_p 
    + \eps^{1/2}\Norm{\eta}_p 
    + \eps^{1/2}\Norm{\Nabla{t}\xi}_p 
    + \Norm{\eta-\Nabla{t}\xi}_p 
    + \eps^2\Norm{\Nabla{s}\eta}_p \\
&\quad
    + \eps\Norm{\Nabla{t}\eta}_p
    + \eps\Norm{\Nabla{s}\xi}_p
    + \eps^{3/2p}\Norm{\xi}_\infty
    + \eps^{1/2+2/p}\Norm{\eta}_\infty.
\end{split}
\end{equation}
For small $\eps$ this norm is much bigger than the 
$(1,p,\eps)$-norm.

\begin{corollary}\label{cor:inverse}
Assume $\Ss_\Vv$ is Morse-Smale and
let $u\in\Mm^0(x^-,x^+;\Vv)$.
Then, for every $p>1$, there are positive constants 
$c$ and $\eps_0$ such that, for every
$\eps\in(0,\eps_0)$ the following holds. If
$$
    \zeta=(\xi,\eta)\in\im\,(\Dd_u^\eps)^*,\qquad
    \zeta'=(\xi',\eta'):=\Dd^\eps_u\zeta,
$$
then
\begin{equation}\label{eq:NORM1}
    \NORM{\zeta}_\eps\le c\left(
    \Norm{\xi'}_p+\eps^{3/2}\Norm{\eta'}_p\right).
\end{equation}
\end{corollary}

\begin{proof}
Let $c_2$ be the constant of Theorem~\ref{thm:elliptic-eps}
and $c_3$ be the constant of Theorem~\ref{thm:inverse}.
Then, by Theorem~\ref{thm:inverse},
\begin{equation*}
\begin{split}
&
    \Norm{\xi}_p + \eps^{1/2}\Norm{\eta}_p 
    + \eps^{1/2}\Norm{\Nabla{t}\xi}_p \\
&\le 
    c_3\left(
    \eps\Norm{\xi'}_p+\eps^2\Norm{\eta'}_p
    + \Norm{(\1-\eps\Nabla{t}\Nabla{t})^{-1}
       (\xi'-\eps^2\Nabla{t}\eta')}_p
    \right)  \\
&\le 
    c_3\left(
    (1+\eps)\Norm{\xi'}_p+
    (\eps^2+\kappa_p\eps^{3/2})\Norm{\eta'}_p
    \right) \\
&\le
    c_4\left(\Norm{\xi'}_p+\eps^{3/2}\Norm{\eta'}_p
    \right).
\end{split}
\end{equation*}
Here the second step follows from Lemma~\ref{le:eat-eps}. 
Combining the last estimate with Theorem~\ref{thm:elliptic-eps}
we obtain
\begin{equation*}
\begin{split}
&
    \Norm{\eta-\Nabla{t}\xi}_p 
    + \eps\Norm{\Nabla{t}\eta}_p
    + \eps\Norm{\Nabla{s}\xi}_p
    + \eps^2\Norm{\Nabla{s}\eta}_p \\
&\le 
    c_2\eps\left(
    \Norm{\xi'}_p + \eps\Norm{\eta'}_p
    + \Norm{\xi}_p + \eps^2\Norm{\eta}_p
    \right) \\
&\le 
    c_2\left(
    \eps\Norm{\xi'}_p + \eps^2\Norm{\eta'}_p
    +c_4\eps\left(\Norm{\xi'}_p+\eps^{3/2}\Norm{\eta'}_p\right)
    \right) \\
&\le 
    c_2(1+c_4)
    \left(\eps\Norm{\xi'}_p+\eps^2\Norm{\eta'}_p
    \right).
\end{split}
\end{equation*}
Now let $c_5$ be the constant of Lemma~\ref{le:balanced} below.
Then 
\begin{equation}\label{eq:infty}
\begin{split}
    \eps^{3/2p}\Norm{\xi}_\infty
    &\le c_5\left(\Norm{\xi}_p+\eps^{1/2}\Norm{\Nabla{t}\xi}_p
        +\eps\Norm{\Nabla{s}\xi}_p\right), \\
    \eps^{1/2+2/p}\Norm{\eta}_\infty
    &\le c_5\left(\eps^{1/2}\Norm{\eta}_p+\eps\Norm{\Nabla{t}\eta}_p
        +\eps^2\Norm{\Nabla{s}\eta}_p\right).
\end{split}
\end{equation}
(Here we used the cases $(\beta_1,\beta_2)=(1/2,1)$ and
$(\beta_1,\beta_2)=(1/2,3/2)$.)
Combining these four estimates we obtain~(\ref{eq:NORM1}).  
\end{proof}

The second estimate in the proof of Corollary~\ref{cor:inverse}
shows that one can obtain a stronger estimate than~(\ref{eq:NORM1}) 
from Theorems~\ref{thm:elliptic-eps} and~\ref{thm:inverse}.
Namely, (\ref{eq:NORM1}) continues to hold if $\NORM{\zeta}_\eps$
is replaced by the stronger norm where the $L^p$ norms 
of $\Nabla{t}\xi-\eta$, $\Nabla{t}\eta$, $\Nabla{s}\xi$, 
and $\Nabla{s}\eta$ are multiplied by an additional 
factor $\eps^{-1/2}$.   The reason for not using this stronger 
norm lies in the proof of Theorem~\ref{thm:existence}.  
In the first step of the iteration we solve an equation of 
the form $\Dd^\eps_u\zeta_0=\zeta'=(0,\eta')$
where $\eta'$ is bounded (in $L^p$) with all its derivatives. 
Our goal in this first step is to obtain the sharpest possible 
estimate for $\zeta_0$ and its first derivatives.  We shall see that
this estimate has the form $\NORM{\zeta_0}_\eps\le c\eps^2$
and that such an estimate in terms of $\eps^2$ cannot be 
obtained with the stronger norm indicated above. 

\begin{lemma}\label{le:balanced}
Let $u\in\Cinf(\R \times S^1,M)$ such that
$\|\p_s u\|_\infty$ and $\|\p_t u\|_\infty$
are finite and
$
     \lim_{s\to\pm\infty}u(s,t)
$
exists, uniformly in $t$. Then, for every $p>2$,
there is a constant $c>0$ such that
$$
     \left\|\xi\right\|_\infty
     \le c \eps^{-(\beta_1+\beta_2)/p}
     \left(\left\|\xi\right\|_p
     +\eps^{\beta_1}\left\|\Nabla{t}\xi\right\|_p
     +\eps^{\beta_2}\left\|\Nabla{s}\xi\right\|_p
     \right)
$$
for every $\eps\in(0,1]$, every pair of nonnegative
real numbers $\beta_1$ and $\beta_2$,
and every compactly supported vector field
$\xi\in \Omega^0(\R\times S^1,u^*TM)$.
\end{lemma}

\begin{proof}
Define
$
     \tilde{u}:Z_\eps
     :=\R\times\left(\R/
     \eps^{-\beta_1}\Z\right)
     \to M
$
and $\tilde{\xi}\in \Omega^0(Z_\eps,\tilde{u}^*TM)$
by
$$
     \tilde{u}(s,t)
     :=u(\eps^{\beta_2}s,
     \eps^{\beta_1}t),\qquad
     \tilde{\xi}(s,t)
     :=\xi(\eps^{\beta_2}s,
     \eps^{\beta_1}t).
$$
The estimate is equivalent to the Sobolev
inequality
$$
     \bigl\|\tilde{\xi}\bigr\|_\infty
     \le c \left(\bigl\|\tilde{\xi}\bigr\|_p
     +\bigl\|\Nabla{t}\tilde{\xi}\bigr\|_p
     +\bigl\|\Nabla{s}\tilde{\xi}\bigr\|_p
     \right)
$$
with a uniform constant $c=c(p,\Norm{\p_su}_\infty,
\Norm{\p_tu}_\infty)$ that is independent
of $\eps\in(0,1]$.
(To see how the $L^\infty$ bounds on $\p_su$ and
$\p_tu$ enter the estimate, embedd $M$
into some euclidean space and use
the Gauss-Weingarten formula.)
\end{proof}


\section{Existence and uniqueness}\label{sec:adiabatic}

Throughout this section we fix a perturbation
$\Vv$ that satisfies~$(V0-V4)$.
In the next theorem we denote by 
$$
     \Phi(x,\xi):T_xM\to T_{\exp_x(\xi)}M
$$
parallel transport along the geodesic $\tau\mapsto\exp_x(\tau\xi)$. 

\begin{theorem}[Existence]\label{thm:existence}
Assume $\Ss_\Vv$ is Morse--Smale and fix 
two constants $a\in\R$ and $p>2$. 
Then there are positive constants $c$ and 
$\eps_0$ such that the following holds.
For every $\eps\in(0,\eps_0)$,
every pair $x^\pm\in\Pp^a(\Vv)$ of index difference one,
and every $u\in\Mm^0(x^-,x^+;\Vv)$, there exists a pair 
$(u^\eps,v^\eps)\in\Mm^\eps(x^-,x^+;\Vv)$ of the form
$$
      u^\eps = \exp_u(\xi),\qquad
      v^\eps = \Phi(u,\xi)(\p_tu+\eta),\qquad
      (\xi,\eta)\in\im\,(\Dd_u^\eps)^*,
$$
where $\xi$ and $\eta$ satisfy the inequalities
\begin{equation}\label{eq:exist1}
\begin{split}
    &\left\|\Nabla{t}\xi-\eta\right\|_{L^p}
     +\left\|\xi\right\|_{L^p}
     +\eps^{1/2}\left\|\eta\right\|_{L^p}
     +\eps^{1/2}\left\|\Nabla{t}\xi\right\|_{L^p} \\
    &+\eps\left\|\Nabla{t}\eta\right\|_{L^p}
     +\eps\left\|\Nabla{s}\xi\right\|_{L^p}
     +\eps^2\left\|\Nabla{s}\eta\right\|_{L^p} 
     \le c\eps^2
\end{split}
\end{equation}
and
\begin{equation}\label{eq:exist2}
     \left\|\xi\right\|_{L^\infty}
     \le c\eps^{2-3/2p}, \qquad
     \left\|\eta\right\|_{L^\infty}
     \le c\eps^{3/2-2/p}.
\end{equation}
\end{theorem}

\begin{remark}\rm
The estimates~(\ref{eq:exist1}) and (\ref{eq:exist2})
can be summarized in the form
$$
     \NORM{\zeta}_\eps\le c\eps^2
$$
(with a larger constant $c$).
\end{remark}

\begin{theorem}[Uniqueness]\label{thm:unique}
Assume $\Ss_\Vv$ is Morse--Smale and fix 
two constants $a\in\R$ and $C>0$. 
Then there are positive constants $\delta$ and 
$\eps_0$ such that, for every $\eps\in(0,\eps_0)$,
every pair $x^\pm\in\Pp^a(\Vv)$ of index difference one,
and every $u\in\Mm^0(x^-,x^+;\Vv)$ the following holds.
If 
\begin{equation}\label{eq:unique}
     (\xi_i,\eta_i)\in\im\,(\Dd_u^\eps)^*,\qquad
     \left\|\xi_i\right\|_{L^\infty}\le\delta\eps^{1/2},\qquad
     \left\|\eta_i\right\|_{L^\infty}\le C,
\end{equation}
for $i=1,2$ and the pairs
$$
      u_i^\eps := \exp_u(\xi_i),\qquad
      v_i^\eps := \Phi(u,\xi_i)(\p_tu+\eta_i),
$$
belong to the moduli space $\Mm^\eps(x^-,x^+;\Vv)$,
then $(u_1^\eps,v_1^\eps)=(u_2^\eps,v_2^\eps)$.
\end{theorem}

In the hypotheses of Theorem~\ref{thm:unique} we did not
specify the Sobolev space to which $\zeta_i=(\xi_i,\eta_i)$
is required to belong.  The reason is that $\zeta_i$ is smooth 
and, by exponential decay, belongs to the Sobolev space 
$W^{k,p}(\R\times S^1,u^*TM\oplus u^*TM)$ for every integer
$k\ge0$ and every $p\ge1$.

\begin{definition}\label{def:T}
Assume $\Ss_\Vv$ is Morse--Smale and fix 
three constants $a\in\R$, $C>0$, and $p>2$.
Choose positive constants $\eps_0$, $\delta$, and $c$ such that
the assertions of Theorem~\ref{thm:existence} and~\ref{thm:unique}
hold with these constants. Shrink $\eps_0$ so that
$c\eps_0^{1/2}<\delta$ and $c\eps_0^{1/2}\le C$. 
Define the map
$$
      \Tt^\eps:\Mm^0(x^-,x^+;\Vv)\to\Mm^\eps(x^-,x^+;\Vv)
$$
by
$$
      \Tt^\eps(u):=(u^\eps,v^\eps),\qquad
      u^\eps := \exp_u(\xi),\qquad
      v^\eps := \Phi(u,\xi)(\p_tu+\eta),
$$
where the pair
$
      (\xi,\eta)\in\im\,(\Dd_u^\eps)^*
$
is chosen such that~(\ref{eq:exist1}) and~(\ref{eq:exist2}) are 
satisfied and $(\exp_u(\xi),\Phi(u,\xi)(\p_tu+\eta))\in\Mm^\eps(x^-,x^+;\Vv)$.
Such a pair $(\xi,\eta)$ exists, by Theorem~\ref{thm:existence},
and is unique, by Theorem~\ref{thm:unique}. 
The map $\Tt^\eps$ is shift equivariant. 
 \end{definition}

The proof of Theorem~\ref{thm:existence}
is based on the Newton--Picard iteration method
to detect a zero of a map near an approximate zero.
The first step is to define
a suitable map between Banach spaces.
In order to do so let $(u,v):\R \times S^1 \to TM$ be a smooth
map and consider the map $\Ff_{u,v}^\eps:
W^{1,p}(\R \times S^1,u^*TM \oplus u^*TM)
\to L^p(\R \times S^1,u^*TM \oplus u^*TM)$
given by
\begin{equation}\label{eq:Ff_(u,v)^eps}
     \Ff_{u,v}^\eps \begin{pmatrix} \xi \\ \eta \end{pmatrix}
     :=
     \begin{pmatrix} \Phi(u,\xi)^{-1} & 0 \\
     0 & \Phi(u,\xi)^{-1} \end{pmatrix}
     \Ff_\eps \begin{pmatrix} \exp_u \xi \\
     \Phi(u,\xi) (v + \eta) \end{pmatrix},
\end{equation}
where
\begin{equation}\label{eq:Ff_eps}
     \Ff_\eps \begin{pmatrix} u^\eps \\ v^\eps \end{pmatrix}
     := \begin{pmatrix} \p_s u^\eps -\Nabla{t} v^\eps -\grad\Vv(u^\eps) \\
     \Nabla{s} v^\eps + \eps^{-2} (\p_t u^\eps -v^\eps) \end{pmatrix} .
\end{equation}
Thus, abbreviating $\Phi:=\Phi(u,\xi)$, we have
$$
     \Ff_{u,v}^\eps \begin{pmatrix} \xi \\ \eta \end{pmatrix}
     :=
     \begin{pmatrix} \Phi^{-1}
     \left(\p_s\exp_u(\xi) - \Nabla{t}(\Phi(v + \eta))
     -\grad\Vv(\exp_u(\xi))\right) \\
     \Phi^{-1}\left(\Nabla{s}(\Phi(v + \eta))
     + \eps^{-2}\p_t\exp_u(\xi)\right)-\eps^{-2}(v + \eta)
      \end{pmatrix}.
$$
Moreover, 
the differential of $\Ff_{u,v}^\eps$ at the origin
is given by
$
     d\Ff_{u,v}^\eps (0,0)=\Dd_{u,v}^\eps
$
(see~\cite[Appendix A.3]{JOA}).

One of the key ingredients in the iteration
is to have control over the variation
of derivatives. This is provided
by the following quadratic estimates.

\begin{proposition} \label{prop:quadest1}
There exists a constant $\delta>0$ with the following 
significance.  For every $p>1$ and every $c_0>0$ there
is a constant $c>0$ such that the following is true.
Let $(u,v):\R\times S^1\to TM$ be a smooth map
and $Z=(X,Y),\zeta =(\xi , \eta ) \in
\Omega^0(\R \times S^1, u^*TM \oplus u^*TM)$
be two pairs of vector fields along $u$ 
such that
$$
     \|\p_su\|_\infty+\|\p_tu\|_\infty+\|v\|_\infty\le c_0,\quad
     \|\xi\|_\infty+\|X\|_\infty \le \delta,\quad
     \|\eta\|_\infty+\|Y\|_\infty \le c_0.
$$
Then the vector fields $F_1$, $F_2$ along $u$, defined by
$$
     \Ff^\eps_{u,v}(Z+\zeta) 
     -\Ff^\eps_{u,v}(Z) 
     -d\Ff^\eps_{u,v} (Z) \zeta
     =:\begin{pmatrix} F_1 \\ F_2 \end{pmatrix} ,
$$
satisfy the inequalities
\begin{equation*}
\begin{split}
      \| F_1 \|_p 
     &\le c \|\xi\|_\infty \Bigl( \|\xi\|_p
      +\|\eta\|_p+\|\Nabla{t} \xi\|_p
      +\|\Nabla{s} \xi\|_p \|\xi\|_\infty \Bigr) \\
     &\quad+ c\Bigl(\|\Nabla{t}X\|_p
       +\|\Nabla{s}X\|_p\Bigr)\|\xi\|_\infty^2
      +c\|\Nabla{t} X\|_p\|\xi\|_\infty
      \|\eta\|_\infty \\
     &\quad + c \|X\|_\infty \Bigl(  
      \|\Nabla{s} \xi\|_p \|\xi\|_\infty
      + \|\Nabla{t} \xi\|_p \|\eta\|_\infty \Bigr), \\
      \| F_2 \|_p 
     &\le c \|\xi\|_\infty \Bigl(\eps^{-2}\|\xi\|_p 
      +\|\eta\|_p +\|\Nabla{s}\xi\|_p
      +\eps^{-2} \|\Nabla{t} \xi\|_p \|\xi\|_\infty
      \Bigr) \\
     &\quad+ c\Bigl(\|\Nabla{s}X\|_p
      +\eps^{-2}\|\Nabla{t}X\|_p\Bigr)\|\xi\|_\infty^2
       +c\|\Nabla{s}X\|_p\|\xi\|_\infty\|\eta\|_\infty\\
     &\quad + c \|X\|_\infty \Bigl( 
      \eps^{-2} \|\Nabla{t} \xi\|_p \|\xi\|_\infty
       + \|\Nabla{s} \xi\|_p \|\eta\|_\infty \Bigr).
\end{split}
\end{equation*}
\end{proposition}

\begin{proposition} \label{prop:quadest2}
There exists a constant $\delta>0$ with the following 
significance.  For every $p>1$ and every $c_0>0$ there
is a constant $c>0$ such that the following is true.
Let $(u,v):\R\times S^1\to TM$ be a smooth map
and $Z=(X,Y),\zeta =(\xi , \eta ) \in
\Omega^0(\R \times S^1, u^*TM \oplus u^*TM)$
be two pairs of vector fields along $u$ 
such that
$$
     \|\p_su\|_\infty+\|\p_tu\|_\infty+\|v\|_\infty\le c_0,\qquad
     \|X\|_\infty \le \delta,\qquad
     \|Y\|_\infty \le c_0.
$$
Then the vector fields $F_1$, $F_2$ along $u$, defined by
$$
     d\Ff^\eps_{u,v} (Z) \zeta
     -d\Ff^\eps_{u,v} (0) \zeta
     =:\begin{pmatrix} F_1 \\ F_2 \end{pmatrix},
$$
satisfy the inequalities
\begin{equation*}
\begin{split}
      \| F_1 \|_p
     &\le c \| \xi \|_\infty \Bigl( \|X\|_p 
      +\|Y\|_p + \|\Nabla{t} X\|_p
      + \|\Nabla{s} X\|_p \|X\|_\infty \Bigr) \\
     &\quad +c \|X\|_\infty \Bigl( \|\eta\|_p 
      + \|\Nabla{t} \xi\|_p + \|\Nabla{s} \xi\|_p \|X\|_\infty 
      + \|\Nabla{t} X\|_p \|\eta\|_\infty \Bigr), \\
      \| F_2 \|_p
     &\le c \| \xi \|_\infty \Bigl(
       \eps^{-2}\|X\|_p
       +\eps^{-2}\|\Nabla{t} X\|_p \|X\|_\infty
      +\|Y\|_p+\|\Nabla{s}X\|_p \Bigr) \\
     &\quad +c \|X\|_\infty \Bigl( \eps^{-2} 
      \|\Nabla{t} \xi\|_p \|X\|_\infty + \| \eta \|_p 
      + \|\Nabla{s} \xi\|_p 
      + \| \Nabla{s} X \|_p \| \eta \|_\infty \Bigr).
\end{split}
\end{equation*}
\end{proposition}

For the proof of Propositions~\ref{prop:quadest1}
and~\ref{prop:quadest2} we refer
to~\cite[Chapter~5]{JOA}.
To understand the estimate of
Proposition~\ref{prop:quadest2}
note that $\eta$ and $Y$ appear only as
zeroth order terms, that $\Nabla{s}\xi$ and
$\Nabla{s} X$ appear only in cubic terms
in $F_1$, and that $\Nabla{t}\xi$ and
$\Nabla{t} X$ appear only in cubic terms
in $F_2$. This follows from the fact
that the first component of $\Ff_\eps$
is linear in $\p_su$ and the second
component is linear in $\p_tu$.
In Proposition~\ref{prop:quadest1}
we have included cubic terms that arise
when the derivative hits $X$.
In this case we must use the $L^\infty$ norms
on the factors $\xi$ and $\eta$ and can
profit from the fact that $\Nabla{s}X$
and $\Nabla{t}X$ will be small in $L^p$.
The constant $\delta$ appears as a condition
for the pointwise quadratic estimates in suitable
coordinate charts on $M$. 

We now reformulate the quadratic estimates 
in terms of the norm~(\ref{eq:NORM}).

\begin{corollary}\label{cor:quad}
There exists a constant $\delta>0$ with the following 
significance.  For every $p>1$ and every $c_0>0$ there
is a constant $c>0$ such that the following holds.
If $(u,v)$, $Z=(X,Y)$ and $\zeta=(\xi,\eta)$
satisfy the hypotheses of Proposition~\ref{prop:quadest1}
then 
\begin{equation*}
\begin{split}
&
    \Norm{\Ff^\eps_{u,v}(Z+\zeta)-\Ff^\eps_{u,v}(Z)
    -d\Ff^\eps_{u,v}(Z)\zeta}_{0,p,\eps^{3/2}} \\
&\le 
     c\NORM{\zeta}_\eps\left(
     \eps^{-1/2}\Norm{\xi}_\infty+\eps^{-1}\Norm{\xi}_\infty^2
     \right) 
     + c\eps^{-1-3/2p}\NORM{Z}_\eps\NORM{\zeta}_\eps
        \Bigl(\Norm{\xi}_\infty+\eps^{1/2}\Norm{\eta}_\infty\Bigr).
\end{split}
\end{equation*}
If $(u,v)$, $Z=(X,Y)$ and $\zeta=(\xi,\eta)$ satisfy 
the hypotheses of Proposition~\ref{prop:quadest2} then 
\begin{equation*}
     \Norm{d\Ff^\eps_{u,v}(Z)\zeta
     -d\Ff^\eps_{u,v}(0)\zeta}_{0,p,\eps^{3/2}} 
     \le c\left(\eps^{-1/2-3/2p}\NORM{Z}_\eps
     +\eps^{-1-7/2p}\NORM{Z}_\eps^2\right)
      \NORM{\zeta}_\eps.
\end{equation*}
\end{corollary}

\begin{proof}
The result follows from Propositions~\ref{prop:quadest2}
and~\ref{prop:quadest1} via term by term inspection.
In particular, we must use the inequalities
$$
     \Norm{\xi}_\infty+\eps\Norm{\eta}_\infty\le 
     c\eps^{-3/p}\Norm{\zeta}_{1,p,\eps},\qquad 
     \Norm{X}_\infty\le\eps^{-3/2p}\NORM{Z}_\eps 
$$
at various places.
The first follows from Lemma~\ref{le:balanced}
with $(\beta_1,\beta_2)=(1,2)$ and the second
from the definition of the norm in~(\ref{eq:NORM}).
\end{proof}

\begin{proof}[Proof of
Theorem~\ref{thm:existence}.]
Given $u\in \Mm^0(x^-,x^+;\Vv)$
with $x^\pm\in \Pp^a(\Vv)$
we aim to detect an element of $\Mm^\eps(x^-,x^+;\Vv)$
near $u$. We set $v:=\p_tu$ and carry out 
the Newton--Picard iteration method for the map
$\Ff_u^\eps:=\Ff_{u,\p_tu}^\eps$.
Key ingredients are a small initial value, 
a uniformly bounded right inverse 
and control over the variation
of derivatives (which is provided
by the quadratic estimates above).
Because $\Ss_\Vv$ is Morse-Smale, the sets
$\Pp^a(\Vv)$ and $\Mm^0(x^-,x^+;\Vv)/\R$
are finite (the latter in addition relies on 
the assumption of index difference one).
All constants appearing below turn out to
be invariant under $s$-shifts of $u$.
Hence they can be chosen to depend on $a$ only.

Since $u\in \Mm^0(x^-,x^+;\Vv)$ it follows
from Theorems~\ref{thm:par-apriori}
and~\ref{thm:par-exp-decay}
that there is a constant $c_0>0$ such that
\begin{equation}\label{eq:aprior-bound1}
     \left\|\p_su\right\|_\infty
     +\left\|\p_tu\right\|_\infty
     +\left\|\Nabla{t}\p_tu\right\|_\infty
     \le c_0
\end{equation}
and
\begin{equation}\label{eq:aprior-bound2}
     \left\|\Nabla{t}\p_su\right\|_\infty
     +\left\|\Nabla{t}\p_su\right\|_p
    +\left\|\Nabla{t}\Nabla{t}\p_su\right\|_p
     \le c_0.
\end{equation}
Thus the assumptions in Theorem~\ref{thm:elliptic-eps},
Theorem~\ref{thm:inverse}, Proposition~\ref{prop:quadest1},
Proposition~\ref{prop:quadest2} and Lemma~\ref{le:balanced} 
are satisfied. Moreover, by~(\ref{eq:aprior-bound2}) 
the value of the initial point $Z_0:=0$
is indeed small with respect to the
$(0,p,\eps)$-norm:
\begin{equation}\label{eq:initial-value}
     \left\|\Ff_u^\eps(0)\right\|_{0,p,\eps}
     =\left\|\Ff^\eps(u,\p_tu)\right\|_{0,p,\eps}
     =\left\|\begin{pmatrix} 0 \\ \Nabla{s} \p_tu
     \end{pmatrix} \right\|_{0,p,\eps}
     \le c_0 \eps .
\end{equation}
Here we used in addition (\ref{eq:Ff_(u,v)^eps}), 
(\ref{eq:Ff_eps}) and the parabolic equations.
Define the initial correction term
$\zeta_0=(\xi_0,\eta_0)$ by
$$
     \zeta_0
   :=- {\Dd_u^\eps}^* 
     (\Dd_u^\eps{\Dd_u^\eps}^*)^{-1}
     \Ff_u^\eps (0) .
$$
Recursively, for $\nu\in\N$, define
the sequence of correction terms 
$\zeta_\nu=(\xi_\nu,\eta_\nu)$
by
\begin{equation}\label{eq:zeta-nu}
     \zeta_\nu
   :=- {\Dd_u^\eps}^* 
     (\Dd_u^\eps{\Dd_u^\eps}^*)^{-1}
     \Ff_u^\eps (Z_\nu),\qquad
     Z_\nu =(X_\nu,Y_\nu)
     :=\sum_{\ell=0}^{\nu-1} \zeta_\ell.
\end{equation}
We prove by induction that there is a constant $c>0$
such that 
$$
\NORM{\zeta_\nu}_\eps
\le \frac{c}{2^\nu}\eps^2,\qquad
\left\|\Ff_u^\eps(Z_{\nu+1})\right\|_{0,p,\eps^{3/2}}
     \le\frac{c}{2^\nu}\eps^{7/2-3/2p}.
\eqno{(H_\nu)}
$$

\medskip
\noindent
{\bf Initial Step: \boldmath$\nu =0$\unboldmath.}
By definition of $\zeta_0$ we have
$$
     \Dd_u^\eps \zeta_0
     =-\Ff_u^\eps (0)
     =\begin{pmatrix}0\\-\Nabla{s}\p_tu
     \end{pmatrix}.
$$
Thus, by Theorem~\ref{thm:inverse} 
(with constant $c_1>0$), 
\begin{equation*}
\begin{split}
      \|\xi_0\|_p+\eps^{1/2}\|\eta_0\|_p
      +\eps^{1/2}\|\Nabla{t}\xi_0\|_p
     &\le c_1\left(\eps\|(0,\Nabla{s}\p_tu)\|
      _{0,p,\eps}
      +\|\pi_\eps(0,\Nabla{s}\p_t u)\|_p \right) \\
     &\le c_1\left(\eps^2\|\Nabla{s}\p_tu\|_p
      +\eps^2\|\Nabla{t}\Nabla{s}\p_t u\|_p
      \right)\\
     &\le c_0c_1\eps^2.
\end{split} 
\end{equation*}
Here the second inequality follows from
Lemma~\ref{le:eat-eps} and the last 
from~(\ref{eq:aprior-bound2}).
By Theorem~\ref{thm:elliptic-eps} 
(with constant $c_2>0$),
\begin{equation*}
 \begin{split}
    &\|\Nabla{t}\xi_0-\eta_0\|_p
    +\eps\|\Nabla{t}\eta_0\|_p+\eps\|\Nabla{s}\xi_0\|_p
     +\eps^2\|\Nabla{s}\eta_0\|_p \\
    &\le c_2\eps\left(\Norm{(0,\Nabla{s}\p_t u)}_{0,p,\eps}
     +\|\xi_0\|_p+\eps^2\|\eta_0\|_p\right) \\
    &\le c_2\eps\left(\eps\|\Nabla{s}\p_t u\|_p
     + c_0c_1\eps^2\right) \\
    &\le c_0c_2(1+c_1\eps)\eps^2.
\end{split}
\end{equation*}
The last inequality follows again from~(\ref{eq:aprior-bound2}).
Combining these two estimates with~(\ref{eq:infty})
we obtain 
\begin{equation}\label{eq:zeta0}
     \eps^{3/2p}\Norm{\xi_0}_\infty
     + \eps^{1/2+2/p}\Norm{\eta_0}_\infty
     \le \NORM{\zeta_0}_\eps\le 
     c\eps^2.
\end{equation}
with a suitable constant $c>0$ (depending only on $c_0,c_1,c_2$
and the constant of Lemma~\ref{le:balanced}).  
This proves the first estimate in~($H_\nu$)
for $\nu=0$.  To prove the second estimate we observe
that $Z_1=\zeta_0$ and hence, by Proposition~\ref{prop:quadest1}
(with constant $c_3>0$),
\begin{equation*} 
\begin{split}
     &\left\|\Ff_u^\eps(Z_1)\right\|_{0,p,\eps^{3/2}}\\
    &=\left\|\Ff_u^\eps(\zeta_0)
     -\Ff_u^\eps(0)
     -\Dd_u^\eps\zeta_0\right\|_{0,p,\eps^{3/2}} \\
    &\le c_3\|\xi_0\|_\infty\Bigl(
     \|\xi_0\|_p+\|\eta_0\|_p+\|\Nabla{t} \xi_0\|_p
     +\|\Nabla{s}\xi_0\|_p \|\xi_0\|_\infty
     \Bigr) \\
    &\quad+c_3\eps^{3/2}\|\xi_0\|_\infty
     \Bigl(\eps^{-2}\|\xi_0\|_p
     +\|\eta_0\|_p+\|\Nabla{s}\xi_0\|_p
     +\eps^{-2}\|\Nabla{t}\xi_0\|_p\|\xi_0\|_\infty
     \Bigr) \\
    &\le c\eps^{7/2-3/2p}.
\end{split} 
\end{equation*}
with a suitable constant $c>0$ (depending only on $c_0,c_1,c_2$
and the constant of Lemma~\ref{le:balanced}).  
Thus we have proved~($H_\nu$) for $\nu=0$. 
{From} now on we fix the constant $c$ for which the 
estimate~($H_0$) has been established. 

\medskip
\noindent
{\bf Induction step:
\boldmath$\nu-1\Rightarrow\nu$\unboldmath.}
Let $\nu\ge 1$ and assume that 
$(H_0),\dots,(H_{\nu-1})$ are true.
Then
$$
  \NORM{Z_\nu}_\eps
  \le\sum_{\ell=0}^{\nu-1}\NORM{\zeta_\ell}_\eps
  \le c\eps^2\sum_{\ell=0}^{\nu-1} 2^{-\ell}
  \le 2c\eps^2,
$$
$$
   \left\|\Ff_u^\eps(Z_\nu)\right\|_{0,p,\eps^{3/2}}
   \le\frac{c}{2^{\nu-1}}\eps^{7/2-3/2p}.
$$
By~(\ref{eq:zeta-nu}) we have 
$$
   \Dd^\eps_u\zeta_\nu=-\Ff_u^\eps(Z_\nu),\qquad   
   \zeta_\nu\in\im(\Dd^\eps_u)^*.
$$
Hence, by Corollary~\ref{cor:inverse},
(with constant $c_4>0$), 
\begin{equation}\label{eq:zetanu}
\NORM{\zeta_\nu}_\eps
\le c_4\Norm{\Ff^\eps_u(Z_\nu)}_{0,p,\eps^{3/2}}
\le \frac{cc_4}{2^{\nu-1}}\eps^{7/2-3/2p}
\le \frac{c}{2^\nu}\eps^2.
\end{equation}
The last inequality holds whenever 
$c_4\eps^{3/2-3/2p}\le1/2$. 

By what we have just proved the vector fields $Z_\nu$
and $\zeta_\nu$ satisfy the requirements of 
Corollary~\ref{cor:quad} (with the constant $c_5>0$).
Hence 
\begin{equation*}
\begin{split}
    \Norm{\Ff_u^\eps(Z_{\nu+1})}_{0,p,\eps^{3/2}}
    &\le\Norm{\Ff_u^\eps(Z_\nu+\zeta_\nu)
     -\Ff_u^\eps(Z_\nu)-d\Ff_u^\eps(Z_\nu)\zeta_\nu}_{0,p,\eps^{3/2}} \\
    &\quad+\Norm{d\Ff_u^\eps(Z_\nu)\zeta_\nu
     -\Dd_u^\eps\zeta_\nu}_{0,p,\eps^{3/2}} \\
    &\le
     c_5\left(\eps^{-1/2}\Norm{\xi_\nu}_\infty+\eps^{-1}
         \Norm{\xi_\nu}_\infty^2\right)
     \NORM{\zeta_\nu}_\eps \\
    &\quad
     + c_5\eps^{-1-3/2p}\NORM{Z_\nu}_\eps
        \Bigl(\Norm{\xi_\nu}_\infty+\eps^{1/2}\Norm{\eta_\nu}_\infty\Bigr)
        \NORM{\zeta_\nu}_\eps \\
    &\quad
     + c_5\eps^{-1/2-3/2p}\NORM{Z_\nu}_\eps\NORM{\zeta_\nu}_\eps 
     +c_5\eps^{-1-7/2p}\NORM{Z_\nu}_\eps^2\NORM{\zeta_\nu}_\eps \\
    &\le
     c_5\left(c\eps^{3/2-3/2p}
     + c^2\eps^{3-3/p}\right)\NORM{\zeta_\nu}_\eps
     + 2c^2c_5\eps^{3-7/2p}\NORM{\zeta_\nu}_\eps \\
    &\quad
     + 2cc_5\eps^{3/2-3/2p}\NORM{\zeta_\nu}_\eps
     + 4c^2c_5\eps^{3-7/2p}\NORM{\zeta_\nu}_\eps \\
    &\le 
     \frac{1}{2c_4}\NORM{\zeta_\nu}_\eps \\
    &\le 
     \frac{c}{2^\nu}\eps^{7/2-3/2p}.
\end{split}
\end{equation*}
In the third step we have used the inequalities
$$
     \Norm{\xi_\nu}_\infty\le 
     \eps^{-3/2p}\NORM{\zeta_\nu}_\eps\le c\eps^{2-3/2p}
$$
and
$$
     \Norm{\xi_\nu}_\infty+\eps^{1/2}\Norm{\eta_\nu}_\infty
     \le \eps^{-2/p}\NORM{\zeta_\nu}_\eps\le c\eps^{2-2/p}
$$
as well as $\NORM{Z}_\nu\le 2c\eps^2$.
The fourth step holds for $\eps$ sufficiently small,
and the last step follows from~(\ref{eq:zetanu}). 
This completes the induction and proves~($H_\nu$)
for every $\nu$. 

It follows from~($H_\nu$) that 
$Z_\nu$ is a Cauchy sequence with respect 
to $\NORM{\cdot}_\eps$.  
Denote its limit by
$$
    \zeta := \lim_{\nu\to\infty}Z_\nu = \sum_{\nu=0}^\infty\zeta_\nu.
$$
By construction and by~($H_\nu$), the limit satisfies
$$
    \NORM{\zeta}_\eps \le 2c\eps^2,\qquad
    \Ff^\eps_\nu(\zeta)=0,\qquad \zeta\in\im\,(\Dd^\eps_u)^*.
$$
Hence, by~(\ref{eq:Ff_(u,v)^eps}), the pair
$$
     (u^\eps,v^\eps)
     :=\left(exp_u(\xi),\Phi(u,\xi)
     (\p_tu+\eta)\right)
$$
is a solution of~(\ref{eq:floer-eps}). 
Since $\NORM{\zeta}_\eps$ is finite it follows 
that $(\p_su^\eps,\Nabla{s}v^\eps)$ is bounded.
Hence, by the standard elliptic bootstrapping arguments
for pseudoholomorphic curves, the shifted functions
$u^\eps(s+\cdot,\cdot),v^\eps(s+\cdot,\cdot)$
converge in the $\Cinf$ topology on every
compact set as $s$ tends to $\pm\infty$.
Since $\zeta\in W^{1,p}$, the limits must be the 
periodic orbits $x^\pm$ and, moreover, the pair
$(\p_su^\eps(s,t),\Nabla{s}v^\eps(s,t))$
converges to zero, uniformly in $t$, as 
$s$ tends to $\pm\infty$.  Hence 
$(u^\eps,v^\eps)\in\Mm^\eps(x^-,x^+;\Vv)$.
Evidently, each step in the iteration including the 
constants in the estimates is invariant
under time shift. This proves the theorem.
\end{proof}

\begin{proof}[Proof of Theorem~\ref{thm:unique}.]
Fix a constant $p>2$ and an index one
parabolic cylinder $u\in \Mm^0(x^-,x^+;\Vv)$.
Denote $v:=\p_tu$ and $\Ff_u^\eps:=\Ff_{u,\p_tu}^\eps$.
As in the proof of Theorem~\ref{thm:existence},
the map $u$ satisfies the estimates~(\ref{eq:aprior-bound1})
and~(\ref{eq:aprior-bound2}). Denote by
$$
     \Tt^\eps(u)
     =\left(\exp_u(X),
     \Phi(u,X)(\p_tu+Y)\right)
$$
the solution of~(\ref{eq:floer-eps})
constructed in  Theorem~\ref{thm:existence}. 
Then
$$
     Z\in \im\,(\Dd_u^\eps)^*,\qquad
     \Ff^\eps_u(Z)=0,\qquad
     \NORM{Z}_\eps\le c\eps^2
$$
for a suitable constant $c>0$. 
Now suppose $(u^\eps,v^\eps)\in\Mm^\eps(x^-,x^+;\Vv)$
satisfies the hypotheses of the theorem. 
This means that there is a pair
$$
     \zeta=(\xi,\eta)\in W^{1,p}(\R\times S^1,u^*TM\oplus u^*TM)
$$
such that
$$
     \zeta\in \im\,(\Dd_u^\eps)^*,\qquad
     \Ff^\eps_u(\zeta)=0,\qquad
     \Norm{\xi}_\infty\le\delta\eps^{1/2},\qquad
     \Norm{\eta}_\infty\le C.
$$
The difference 
$$
     \zeta':=(\xi',\eta'):=\zeta-Z
$$
satisfies the inequalities
$$
     \Norm{\xi'}_\infty\le\delta\eps^{1/2}
     +c\eps^{2-3/2p}\le 2\delta\eps^{1/2},\qquad
     \Norm{\eta'}_\infty\le C+c\eps^{3/2-2/p}\le 2C,
$$
provided that $\eps$ is sufficiently small.
Hence, by Corollary~\ref{cor:inverse} 
(with a constant $c_1>0$) and 
Corollary~\ref{cor:quad}
(with a constant $c_2>0$),
we have
\begin{equation*}
\begin{split}
     \NORM{\zeta'}_\eps
&\le 
     c_1\Norm{\Dd^\eps_u\zeta'}_{0,p,\eps^{3/2}}  \\
&\le
     c_1\Norm{\Ff_u^\eps(Z +\zeta')-\Ff_u^\eps(Z)
      -d{\Ff_u^\eps}(Z)\zeta'}_{0,p,\eps^{3/2}} \\
&\quad
     + c_1\Norm{d\Ff_u^\eps(Z)\zeta'
      -d{\Ff_u^\eps}(0)\zeta'}_{0,p,\eps^{3/2}} \\
&\le
     c_1c_2\left(
     \eps^{-1/2}\Norm{\xi'}_\infty+\eps^{-1}\Norm{\xi'}_\infty^2
     \right)
     \NORM{\zeta'}_\eps \\
&\quad
     + c_1c_2\eps^{-1-3/2p}\NORM{Z}_\eps
        \Bigl(\Norm{\xi'}_\infty+\eps^{1/2}\Norm{\eta'}_\infty\Bigr)
        \NORM{\zeta'}_\eps \\
&\quad
     + c_1c_2\eps^{-1/2-3/2p}\NORM{Z}_\eps\NORM{\zeta'}_\eps 
     + c_1c_2\eps^{-1-7/2p}\NORM{Z}_\eps^2\NORM{\zeta'}_\eps \\
&\le
     c_1c_2\left(2\delta+4\delta^2\right)
     \NORM{\zeta'}_\eps 
     + cc_1c_2\eps^{3/2-3/2p}
        \Bigl(2\delta+2C\Bigr)\NORM{\zeta'}_\eps \\
&\quad
     + cc_1c_2\eps^{3/2-3/2p}\NORM{\zeta'}_\eps 
     + c^2c_1c_2\eps^{3-7/2p}\NORM{\zeta'}_\eps \\
&\le
     \frac12\NORM{\zeta'}_\eps.
\end{split}
\end{equation*}
The last inequality holds when $\delta$ and $\eps$ are sufficiently 
small.  It follows that $\zeta'=0$ and this proves the theorem.
\end{proof}


\section{An apriori estimate} \label{sec:apriori}

\begin{theorem}\label{thm:apriori}
Fix a constant $c_0>0$ and a perturbation
$\Vv:\Ll M\to\R$ that satisfies~$(V0)$ and~$(V1)$. 
Then there is a constant $C=C(c_0,\Vv)>0$ such that
the following holds.  If $0<\eps\le1$ and 
$(u,v):\R\times S^1\to TM$ is a solution 
of~(\ref{eq:floer-V}) such that
\begin{equation}\label{eq:EA}
      E^\eps(u,v)\le c_0,\qquad
      \sup_{s\in\R}\Aa_\Vv(u(s,\cdot),v(s,\cdot))\le c_0
\end{equation}
then
$
     \|v\|_\infty \le C.
$
\end{theorem}

For $\eps=1$ and $\Vv(x)=\int_0^1V_t(x(t))\,dt$ 
this result was proved by Cieliebak~\cite[Theorem~5.4]{Ci94}.
His proof combines the 2-dimensional maximum
principle and the Krein-Rutman theorem.
Our proof is based on the following $L^2$-estimate.

\begin{proposition}\label{prop:L-infty-v-ok}
Fix a constant $c_0>0$ and a perturbation
$\Vv:\Ll M\to\R$ that satisfies~$(V0)$ and~$(V1)$. 
Then there is a constant $c=c(c_0,\Vv)>0$ such that
the following holds.  If $0<\eps\le1$ and 
$(u,v):\R\times S^1\to TM$ is a solution 
of~(\ref{eq:floer-eps}) that satisfies~(\ref{eq:EA})
then 
$$
     \sup_{s\in\R}\int_0^1\Abs{v(s,t)}^2\:dt\le c.
$$
\end{proposition}

\begin{proof}
Define $F:\R\to\R$ by 
$$
F(s):=\int_0^1\Abs{v(s,t)}^2\,dt.
$$
We prove that there is a constant $\mu=\mu(\Vv)>0$
such that 
\begin{equation}\label{eq:muf}
\eps^2F''-F'+\mu F + 1\ge 0. 
\end{equation}
To see this we abbreviate
$$
L_\eps := \eps^2\p_s^2+\p_t^2-\p_s,\qquad
\Ll_\eps := \eps^2\Nabla{s}\Nabla{s}+\Nabla{t}\Nabla{t}-\Nabla{s}.
$$
By~(\ref{eq:floer-V}), we have
\begin{equation}\label{eq:veps}
\Ll_\eps v = -\Nabla{t}\grad\Vv(u)
\end{equation}
and hence 
\begin{equation*}
\begin{split}
L_\eps\frac{|v|^2}{2}
&=\eps^2\Abs{\Nabla{s}v}^2+\Abs{\Nabla{t}v}^2
     +\inner{\Ll_\eps v}{v} \\
&=\eps^2\Abs{\Nabla{s}v}^2+\Abs{\Nabla{t}v}^2
     - \inner{\Nabla{t}\grad\Vv(u)}{v} \\
&\ge \eps^2\Abs{\Nabla{s}v}^2
     +\Abs{\Nabla{t}v}^2 
     - C\bigl(1+\Abs{\p_tu}\bigr)\Abs{v}\\
&\ge \eps^2\Abs{\Nabla{s}v}^2
     +\Abs{\Nabla{t}v}^2 
     - C\bigl(1+\Abs{v}+\eps^2\Abs{\Nabla{s}v}\bigr)\Abs{v}\\
&\ge \frac{\eps^2}{2}\Abs{\Nabla{s}v}^2
     +\Abs{\Nabla{t}v}^2
     - \left(\frac{C^2}{2}+C+\frac{\eps^2C^2}{2}\right)\Abs{v}^2 
     - \frac{1}{2} \\
&\ge - \left(C+C^2\right)\Abs{v}^2 
     - \frac{1}{2}.
\end{split}
\end{equation*}
Here $C$ is the constant in~$(V1)$. 
Integrating this inequality over the interval
$0\le t\le 1$ gives~(\ref{eq:muf}) with
$\mu:=2C+2C^2$. It follows from~(\ref{eq:muf}) 
and Lemma~\ref{le:apriori-1eps} with $f$ replaced by 
$f+1/\mu$ and $r:=1/2$ that
\begin{equation}\label{eq:fint}
F(s)\le 
F(s)+\frac{1}{\mu}
\le 16 c_2 e^{\mu/4}
\int_{s-1}^{s+1}\left(F(\sigma)+\frac{1}{\mu}\right)\,d\sigma
\end{equation}
for every $s\in\R$. 

Next we observe that, by~(\ref{eq:EA}), we have
\begin{equation*}
\begin{split}
c_0 &\ge \Aa_\Vv(u(s,\cdot),v(s,\cdot)) \\
&=\int_0^1\left(\inner{v(s,t)}{\p_tu(s,t)}
     -\frac{\Abs{v(s,t)}^2}{2}\right)dt 
     - \Vv(u(s,\cdot))\\
&=\int_0^1\left(\frac{\Abs{v(s,t)}^2}{2}
     -\eps^2\inner{v(s,t)}{\Nabla{s}v(s,t)}\right)dt 
     - \Vv(u(s,\cdot)) \\
&\ge\int_0^1\left(\frac{\Abs{v(s,t)}^2}{4}
     -\eps^4\Abs{\Nabla{s}v(s,t)}^2\right)dt - C.
\end{split}
\end{equation*}
Here $C$ is the constant in~$(V0)$ and we have used 
the fact that $\p_tu=v-\eps^2\Nabla{s}v$. 
This implies 
$$
F(s)\le 4\left(c_0+C
+ \int_0^1\eps^2\Abs{\Nabla{s}v(s,t)}^2\,dt\right)
$$
for every $s\in\R$. Integrating this inequality we obtain
$$
\int_{s-1}^{s+1}F(\sigma)\,d\sigma
\le 8c_0+8C+8E^\eps(u,v)\le 16c_0+8C.
$$
Now the assertion follows from~(\ref{eq:fint}).
\end{proof}

\begin{proof}[Proof of Theorem~\ref{thm:apriori}]
In the proof of Proposition~\ref{prop:L-infty-v-ok}
we have seen that there is a constant $\mu=\mu(\Vv)>0$
such that every solution $(u,v)$ of~(\ref{eq:floer-eps})
with $0<\eps\le1$ satisfies the inequality
\begin{equation}\label{eq:v}
L_\eps\Abs{v}^2 \ge-\mu\Abs{v}^2-1.
\end{equation}
Now let $(s_0,t_0) \in \R \times S^1$ and apply 
Lemma~\ref{le:apriori-basic-eps}
with $r=1$ to the function 
$
     w:\R\times\R\supset P_1^\eps\to\R,
$
given by
$
     w(s,t):=\Abs{v(s+s_0,t+t_0)}^2+1/\mu
$:
\begin{equation*}
\begin{split}
     \Abs{v(s_0,t_0)}^2
     &\le2c_2e^\mu
      \int_{-1-\eps}^\eps\int_{-1}^1
      \left(\Abs{v(s+s_0,t+t_0)}^2+\frac{1}{\mu}\right)\,dtds \\
     &\le12c_2e^\mu\left(\frac{1}{\mu}+
      \sup_{s\in\R}\int_0^1\Abs{v(s,t)}^2\,dt\right).
\end{split}
\end{equation*}
Hence the result follows from Proposition~\ref{prop:L-infty-v-ok}.
\end{proof}


\section{Gradient bounds}\label{sec:gradient}

\begin{theorem}\label{thm:gradient}
Fix a constant $c_0>0$ and a perturbation
$\Vv:\Ll M\to\R$ that satisfies~$(V0-V3)$.  
Then there is a constant $C=C(c_0,\Vv)>0$ 
such that the following holds. If $0<\eps\le1$ and 
$(u,v):\R\times S^1\to TM$ is a solution 
of~(\ref{eq:floer-V}) that satisfies~(\ref{eq:EA}),
i.e. $E^\eps(u,v)\le c_0$ and 
$\sup_{s\in\R}\Aa_\Vv(u(s,\cdot),v(s,\cdot))\le c_0$,
then 
\begin{equation}\label{eq:uvE}
\begin{split}
&\Abs{\p_su(s,t)}^2+\Abs{\Nabla{s}v(s,t)}^2 \\
&
+\int_{s-1/2}^{s+1/2}\int_0^1\Bigl(
\Abs{\Nabla{t}\p_su}^2 
+ \Abs{\Nabla{s}\p_su}^2
+ \Abs{\Nabla{t}\Nabla{s}v}^2 
+ \eps^2\Abs{\Nabla{s}\Nabla{s}v}^2
\Bigl)  \\
&\le CE_{[s-1,s+1]}^\eps(u,v)
\end{split}
\end{equation}
for all $s$ and $t$. Here $E_I^\eps(u,v)$ denotes the 
energy of $(u,v)$ over the domain $I\times S^1$.
\end{theorem}

\begin{remark}\rm
Note that~(\ref{eq:uvE}) implies the estimate
$$
\Norm{\p_tu-v}_{L^\infty}\le \eps^2\sqrt{CE^\eps(u,v)}
$$
for every solution $(u,v):\R\times S^1\to TM$ 
of~(\ref{eq:floer-V}) that satisfies~(\ref{eq:EA}).
\end{remark}

The proof of Theorem~\ref{thm:gradient} has five steps.  
The first step is a bubbling argument and establishes 
a weak form of the required $L^\infty$ estimate
(with $\p_su$ replaced by $\eps^2\p_su$ and $\Nabla{s}v$
replaced by $\eps^3\Nabla{s}v$). The second step establishes
an $L^2$-version of the estimate for 
$\Norm{\p_su(s,\cdot)}_{L^2(S^1)}+\eps\Norm{\Nabla{s}v(s,\cdot)}_{L^2(S^1)}$.
The third step is an auxiliary result of the same type 
for the second derivatives. The fourth step establishes
the $L^\infty$ bound with $\Nabla{s}v$ replaced by $\eps\Nabla{s}v$.
The final step then proves the theorem in full. 

\begin{lemma}\label{le:weak-L-infty}
Fix a constant $c_0>0$ and a perturbation
$\Vv:\Ll M\to\R$ that satisfies~$(V0-V1)$.  
Then the following holds.
\begin{enumerate}
\item[\rm\bfseries(i)]
For every $\delta>0$ there is an $\eps_0>0$
such that every solution $(u,v):\R\times S^1\to M$
of~(\ref{eq:floer-V}) and~(\ref{eq:EA}) with
$0<\eps\le\eps_0$ satisfies the inequality
\begin{equation}\label{eq:weak-L-infty}
     \eps^2\left\|\p_su\right\|_\infty
      + \eps^3\left\|\Nabla{s}v\right\|_\infty
     \le\delta.
\end{equation}
\item[\rm\bfseries(ii)]
For every $\eps_0>0$ there is a constant $c>0$
such that every solution $(u,v):\R\times S^1\to M$
of~(\ref{eq:floer-V}) and~(\ref{eq:EA}) with
$\eps_0\le\eps\le1$ satisfies
$$
      \left\|\p_su\right\|_\infty
      + \left\|\Nabla{s}v\right\|_\infty
      \le c.
$$
\end{enumerate}
\end{lemma}

\begin{proof}
We prove~(i).  Suppose, by contradiction, that
the result is false. Then there is a sequence of solutions
$(u_\nu,v_\nu):\R\times S^1\to M$ of~(\ref{eq:floer-V})
with $\eps_\nu>0$ satisfying
$$
E^{\eps_\nu}(u_\nu,v_\nu)\le c_0,\qquad
\sup_{s\in\R}\Aa_\Vv(u_\nu(s,\cdot),v_\nu(s,\cdot))\le c_0,\qquad
\lim_{\nu\to\infty}\eps_\nu=0,
$$
and
$$
     \eps_\nu^2\left\|\p_su_\nu\right\|_\infty
     + \eps_\nu^3\left\|\Nabla{s}v_\nu\right\|_\infty
     \ge 2\delta
$$
for suitable constants $c_0>0$ and $\delta>0$. 
Since $(u_\nu,v_\nu)$ has finite energy the functions
$|\p_su_\nu(s,t)|$ and $|\Nabla{s}v_\nu(s,t)|$ converge
to zero as $|s|$ tends to infinity.  Hence the function 
$|\p_su_\nu|+\eps_\nu|\Nabla{s}v_\nu|$ takes on its maximum 
at some point $z_\nu=s_\nu+it_\nu$, i.e.
$$
c_\nu := \sup_{\R\times S^1}
\left(\Abs{\p_su_\nu}+\eps_\nu\Abs{\Nabla{s}v_\nu}\right)
= \Abs{\p_su_\nu(s_\nu,t_\nu)}
+ \eps_\nu\Abs{\Nabla{s}v_\nu(s_\nu,t_\nu)}
$$
and 
\begin{equation}\label{eq:ceps}
\eps_\nu^2c_\nu \ge \delta.
\end{equation}
Applying a time shift and using the periodicity in $t$
we may assume without loss of generality that $s_\nu=0$ 
and $0\le t_\nu\le1$.

Now consider the sequence
$$
     \tilde w_\nu=(\tilde u_\nu,\tilde v_\nu)
     :\R^2\to TM
$$
defined by
$$
     \tilde{u}_\nu(s,t)
     := u_\nu\left(\frac{s}{c_\nu},
     t_\nu+\frac{t}{\eps_\nu c_\nu}\right),
     \qquad
     \tilde{v}_\nu(s,t)
     :=\eps_\nu v_\nu\left(\frac{s}{c_\nu},
     t_\nu+\frac{t}{\eps_\nu c_\nu}\right).
$$
This sequence satisfies the partial differential equation
\begin{equation}\label{eq:floerw}
     \p_s\tilde{u}_\nu-\Nabla{t}\tilde{v}_\nu
     =\frac{1}{c_\nu}\xi_\nu,\qquad
     \Nabla{s}\tilde{v}_\nu+\p_t\tilde{u}_\nu
     =\frac{1}{{\eps_\nu}^2c_\nu}\tilde{v}_\nu,
\end{equation}
where 
$$
\xi_\nu(s,t) := \grad\Vv(u_\nu(s/c_\nu,\cdot))(t_\nu+t/\eps_\nu c_\nu)
\in T_{\tilde{u}_\nu(s,t)}M.
$$
By definition of $c_\nu$ we have
\begin{equation}\label{eq:nonconstant}
\Abs{\p_s\tilde{u}_\nu(0,t_\nu)}
+ \Abs{\Nabla{s}\tilde{v}_\nu(0,t_\nu)} = 1
\end{equation}
and 
$$
\Abs{\p_s\tilde{u}_\nu(s,t)}+\Abs{\Nabla{s}\tilde{v}_\nu(s,t)}\le1
$$
for all $s$ and $t$. Since $\Abs{\tilde{v}_\nu}$ is uniformly
bounded, by Theorem~\ref{thm:apriori}, and $\Abs{\xi_\nu}$ is uniformly 
bounded, by axiom~$(V0)$, it then follows from~(\ref{eq:floerw}) 
that $\tilde{u}_\nu$ and $\tilde{v}_\nu$ are uniformly 
bounded in $C^1$. Moreover, it follows from~$(V1)$ that
$$
\Abs{\Nabla{t}\xi_\nu(s,t)}
\le \frac{C}{\eps_\nu c_\nu}
\Bigl(1+\Abs{\p_tu_\nu(s/c_\nu,t_\nu+t/\eps_\nu c_\nu)}\Bigr)
= C\left(\frac{1}{\eps_\nu c_\nu}
+ \Abs{\p_t\tilde{u}_\nu(s,t)}\right)
$$
and
$$
\Abs{\Nabla{s}\xi_\nu(s,t)}
\le \frac{C}{c_\nu}
\Abs{\p_su_\nu(s/c_\nu,t_\nu+t/\eps_\nu c_\nu)}
= C\Abs{\p_s\tilde{u}_\nu(s,t)}.
$$
Since the sequence $1/\eps_\nu^2c_\nu$ is bounded, 
by~(\ref{eq:ceps}), it now follows from~(\ref{eq:floerw}) 
that $\p_s\tilde{u}_\nu-\Nabla{t}\tilde{v}_\nu$ and 
$\Nabla{s}\tilde{v}_\nu+\p_t\tilde{u}_\nu$ are uniformly
bounded in $C^1$, and hence in $W^{1,p}$ for any $p>2$ 
and on any compact subset of $\R^2$. Since $\tilde{u}_\nu$
and $\tilde{v}_\nu$ are uniformly bounded in $C^1$, this implies
that they are also uniformly bounded in $W^{2,p}$ over every 
compact subset of $\R^2$, by the standard elliptic bootstrapping
techniques for $J$-holomorphic curves (see~\cite[Appendix~B]{MS}).
Hence, by the Arz\'ela--Ascoli theorem, there is a subsequence
that converges in the $C^1$ topology to a solution
$(\tilde{u},\tilde{v})$ of the partial differential
equation 
$$
\p_s\tilde{u}-\Nabla{t}\tilde{v} = 0,\qquad
\Nabla{s}\tilde{v}+\p_t\tilde{u} = \lambda\tilde{v},
$$
where $\lambda=\lim_{\nu\to\infty}1/\eps_\nu^2c_\nu$. 
Since $v_\nu$ is uniformly bounded and $\eps_\nu\to0$ 
we have $\tilde{v}\equiv0$ and so $\tilde{u}$ is constant. 
On the other hand it follows from~(\ref{eq:nonconstant}) 
that $(\tilde{u},\tilde{v})$ is nonconstant; contradiction. 
This proves~(i). 

The proof of~(ii) is almost word by word the same,
except that $\eps_\nu$ no longer converges to zero
while $c_\nu$ still diverges to infinity. So the limit 
$\tilde{w}=(\tilde{u},\tilde{v}):\C\to TM\cong T^*M$
is a $J$-holomorphic curve with finite energy and, 
by removal of singularities, extends to a 
nonconstant $J$-holomorphic sphere $\tilde w:S^2\to T^*M$,
which cannot exist since the symplectic form 
on $T^*M$ is exact.  
\end{proof}

The second step in the proof of Theorem~\ref{thm:gradient} is
to prove an integrated version of the estimate with $\Nabla{s}v$
replaced by $\eps\Nabla{s}v$.

\begin{lemma}\label{le:gradient-L2}
Fix a constant $c_0>0$ and a perturbation
$\Vv:\Ll M\to\R$ that satisfies~$(V0-V2)$.  
Then there is a constant $C=C(c_0,\Vv)>0$ 
such that the following holds. If $0<\eps\le1$ and 
$(u,v):\R\times S^1\to TM$ is a solution 
of~(\ref{eq:floer-V}) that satisfies~(\ref{eq:EA})
then, for every $s\in\R$,
\begin{equation}\label{eq:uvE-L2}
\begin{split}
&\int_0^1\Bigl(
\Abs{\p_su(s,t)}^2+\eps^2\Abs{\Nabla{s}v(s,t)}^2
\Bigr)\,dt \\
&
+\int_{s-1/4}^{s+1/4}\int_0^1\Bigl(
\Abs{\Nabla{t}\p_su}^2 
+ \eps^2\Abs{\Nabla{s}\p_su}^2
+ \eps^2\Abs{\Nabla{t}\Nabla{s}v}^2 
+ \eps^4\Abs{\Nabla{s}\Nabla{s}v}^2
\Bigl)  \\
&\le CE_{[s-1/2,s+1/2]}^\eps(u,v).
\end{split}
\end{equation}
\end{lemma}

\begin{corollary}\label{cor:gradient-L2}
Fix a constant $c_0>0$ and a perturbation
$\Vv:\Ll M\to\R$ that satisfies~$(V0-V2)$.  
Then there is a constant $C=C(c_0,\Vv)>0$ 
such that the following holds. If $0<\eps\le1$ and 
$(u,v):\R\times S^1\to TM$ is a solution 
of~(\ref{eq:floer-V}) that satisfies~(\ref{eq:EA}),
then 
$$
\int_{s-1/4}^{s+1/4}\int_0^1\Abs{\Nabla{s}v}^2
\le CE_{[s-1/2,s+1/2]}^\eps(u,v)
$$
for every $s\in\R$. 
\end{corollary}

\begin{proof}
Since $\Nabla{s}v=\Nabla{t}\p_su+\eps^2\Nabla{s}\Nabla{s}v$
this estimate follows immediately from 
Lemma~\ref{le:gradient-L2}.
\end{proof}

\begin{proof}[Proof of Lemma~\ref{le:gradient-L2}.]
Define the functions $f,g:\R\times S^1\to\R$ by
$$
f := \frac{1}{2}\Bigl(\Abs{\p_su}^2
+\eps^2\Abs{\Nabla{s}v}^2\Bigr)
$$
and
$$
g := \frac{1}{2}\Bigl(\Abs{\Nabla{t}\p_su}^2 
+ \eps^2\Abs{\Nabla{s}\p_su}^2
+ \eps^2\Abs{\Nabla{t}\Nabla{s}v}^2 
+ \eps^4\Abs{\Nabla{s}\Nabla{s}v}^2
\Bigr),
$$
and abbreviate
$$
F(s) := \int_0^1f(s,t)\,dt,\qquad
G(s) := \int_0^1g(s,t)\,dt.
$$
Recall the definition of $L_\eps := \eps^2\p_s^2+\p_t^2-\p_s$
and $\Ll_\eps:=\eps^2\Nabla{s}\Nabla{s}+\Nabla{t}\Nabla{t}-\Nabla{s}$
in the proof of Proposition~\ref{prop:L-infty-v-ok}.
Then 
\begin{equation}\label{eq:Lfg}
L_\eps f = 2g + U + \eps^2V,\qquad
U:=\inner{\p_su}{\Ll_\eps\p_su},\qquad
V:=\inner{\Nabla{s}v}{\Ll_\eps\Nabla{s}v}.
\end{equation}
We shall prove that $U$ and $V$ satisfy the 
pointwise inequality
\begin{equation}\label{eq:UV}
|U|+\eps^2|V|\le \mu f 
+ \frac12\left(
g + \Norm{\p_su}_{L^2(S^1)}^2
+ \eps^4\Norm{\Nabla{s}\p_su}_{L^2(S^1)}^2
\right)
\end{equation}
for a suitable constant $\mu>0$.
Inserting this inequality in~(\ref{eq:Lfg}) gives
$$
L_\eps f + \mu f+F 
\ge g + \frac12(g-G).
$$
Now integrate over the interval $0\le t\le 1$ to obtain
$$
\eps^2F''-F'+ (\mu+1) F\ge G.
$$
With this understood the result follows from
Lemmas~\ref{le:apriori-1eps} and~\ref{le:apriori-1-L2}. 

To prove~(\ref{eq:UV}) we observe that, 
by~(\ref{eq:floer-V}),
\begin{equation}\label{eq:Lf-1}
\begin{split}
\Ll_\eps\p_su 
&= \eps^2\Nabla{s}\Nabla{s}
\left(\Nabla{t}v+\grad\Vv(u)\right)
+\Nabla{t}\Nabla{s}\left(v-\eps^2\Nabla{s}v\right) \\
&\quad
-\Nabla{s}\left(\Nabla{t}v+\grad\Vv(u)\right) \\
&=\eps^2\left[\Nabla{s}\Nabla{s},\Nabla{t}\right]v
+\left[\Nabla{t},\Nabla{s}\right]v
-\Nabla{s}\grad\Vv(u)
+\eps^2\Nabla{s}\Nabla{s}\grad\Vv(u) \\
&=2\eps^2R(\p_su,\p_tu)\Nabla{s}v
+ \eps^2\left(\Nabla{\p_su}R\right)(\p_su,\p_tu)v
- R(\p_su,\p_tu)v \\
&\quad 
+\eps^2R(\Nabla{s}\p_su,\p_tu)v
+\eps^2R(\p_su,\Nabla{s}\p_tu)v \\
&\quad
-\Nabla{s}\grad\Vv(u)
+\eps^2\Nabla{s}\Nabla{s}\grad\Vv(u).
\end{split}
\end{equation}
Now fix a sufficiently small constant $\delta>0$ and choose
$\eps_0>0$ such that the assertion of Lemma~\ref{le:weak-L-infty}~(i)
holds. Choose $C>0$ such that the assertion of Theorem~\ref{thm:apriori}
holds and assume $0<\eps\le\eps_0\le\delta/C$.  Then, by 
Theorem~\ref{thm:apriori} and Lemma~\ref{le:weak-L-infty}, 
we have
\begin{equation}\label{eq:uvs}
     \eps^2\left\|\p_su\right\|_\infty\le\delta,\qquad 
     \eps^3\left\|\Nabla{s}v\right\|_\infty\le\delta,\qquad
     \left\|v\right\|_\infty\le C,\qquad
     \eps\left\|\p_tu\right\|_\infty\le2\delta.
\end{equation}
The last estimate uses the identity $\p_tu=v-\eps^2\Nabla{s}v$.
Now take the pointwise inner product of~(\ref{eq:Lf-1}) with $\p_su$
and estimate the resulting seven expressions separately. 
By~(\ref{eq:uvs}) and~$(V1)$, the terms four, five, and six 
are bounded by the right hand side of~(\ref{eq:UV}). 
For the last term we find, by~$(V2)$, 
\begin{align*}
\eps^2\Abs{\inner{\p_su}{\Nabla{s}\Nabla{s}\grad\Vv(u)}}
&
\le \eps^2C\Abs{\p_su}
\Bigl(
\Abs{\Nabla{s}\p_su}
+ \Norm{\Nabla{s}\p_su}_{L^2(S^1)}\Bigr) \\
&\quad
+ \eps^2C\Abs{\p_su}
\left(\Abs{\p_su}+ \Norm{\p_su}_{L^2(S^1)}\right)^2 \\
&
\le \eps^2C\Abs{\p_su}
\Bigl(
\Abs{\Nabla{s}\p_su}
+ \Norm{\Nabla{s}\p_su}_{L^2(S^1)}\Bigr) \\
&\quad
+ 2C\delta\Abs{\p_su}
\left(\Abs{\p_su}+ \Norm{\p_su}_{L^2(S^1)}\right) \\
&\le
\mu f + \frac18\left(g+\Norm{\p_su}_{L^2(S^1)}^2
+ \eps^4\Norm{\Nabla{s}\p_su}_{L^2(S^1)}^2\right).
\end{align*}
For the first three terms on the right in~(\ref{eq:Lf-1}) 
we argue as follows. Differentiate the equation 
$v=\p_tu+\eps^2\Nabla{s}v$ covariantly with respect 
to $s$ to obtain
\begin{equation}\label{eq:uv}
     \Nabla{s}v
     =\Nabla{s}\p_tu+\eps^2\Nabla{s}\Nabla{s}v,\qquad
     \p_tu = v -\eps^2\Nabla{t}\p_su -\eps^4\Nabla{s}\Nabla{s}v.
\end{equation}
Now express half the first term on the right 
in~(\ref{eq:Lf-1}) in the form
\begin{equation*}
\begin{split}
   &\eps^2\bigl\langle\p_su,
    R(\p_su,\p_tu)\Nabla{s}v
    \bigr\rangle \\
   &=\eps^2\bigl\langle\p_su,
    R(\p_su,v)\Nabla{t}\p_su\bigr\rangle
    +\eps^4\bigl\langle\p_su,
    R(\p_su,v)\Nabla{s}\Nabla{s}v
    \bigr\rangle \\
   &\quad
    -\eps^4\bigl\langle\p_su,
    R(\p_su,\Nabla{t}\p_su)\Nabla{t}\p_su
    \bigr\rangle
    -\eps^6\bigl\langle\p_su,
    R(\p_su,\Nabla{t}\p_su)\Nabla{s}\Nabla{s}v
    \bigr\rangle \\
   &\quad
    -\eps^6\bigl\langle\p_su,
    R(\p_su,\Nabla{s}\Nabla{s}v)\Nabla{t}\p_su
    \bigr\rangle
    -\eps^8\bigl\langle\p_su,
    R(\p_su,\Nabla{s}\Nabla{s}v)\Nabla{s}\Nabla{s}v
    \bigr\rangle.
\end{split}
\end{equation*}
Here we have replaced $\p_tu$ and $\Nabla{s}v$ by the expressions
in~(\ref{eq:uv}). In the first two terms we eliminate one of
the factors $\p_su$ by using the inequality $\eps^2|\p_su|\le\delta$
and in the last four terms we eliminate both factors $\p_su$ by
the same inequality. The next two terms in our expression for $U$
have the form 
$$
    \eps^2\inner{\p_su}{\left(\Nabla{\p_su}R\right)(\p_su,\p_tu)v}
    - \inner{\p_su}{R(\p_su,\p_tu)v}.
$$
Replace $\p_tu$ by the expression in~(\ref{eq:uv}) and elimate 
in each of the resulting summands one or two of the factors
$\eps^2\p_su$ as above.  This proves the required estimate 
for $U$ and $0<\eps\le\eps_0$. 

To estimate $V$ we observe that, by~(\ref{eq:floer-V}),
\begin{equation}\label{eq:Lf-2}
\begin{split}
\Ll_\eps\Nabla{s}v 
&=\Nabla{s}\Nabla{s}(\eps^2\Nabla{s}v-v)
+\Nabla{t}\Nabla{s}\Nabla{t}v
+\Nabla{t}([\Nabla{t},\Nabla{s}]v) \\
&= -\Nabla{s}\Nabla{t}\p_su
+\Nabla{t}\Nabla{s}(\p_su-\grad\Vv(u)) 
-\Nabla{t}(R(\p_su,\p_tu)v) \\
&= - R(\p_su,\p_tu)\p_su 
+ R(\p_su,\p_tu)\grad\Vv(u)
-\Nabla{t}(R(\p_su,\p_tu)v) \\       
&\quad
- \Nabla{s}\Nabla{t}\grad\Vv(u) \\
&=-2R(\p_su,\p_tu)\p_su
+ R(\p_su,\p_tu)\grad\Vv(u) \\
&\quad-\left(\Nabla{\p_tu}R\right)(\p_su,\p_tu)v
-R(\Nabla{t}\p_su,\p_tu)v
-R(\p_su,\Nabla{t}\p_tu)v \\
&\quad
- \Nabla{t}\Nabla{s}\grad\Vv(u).
\end{split}
\end{equation}
The last step uses the identity
$\Nabla{t}v=\p_su-\grad\Vv(u)$.
Now take the pointwise inner product with $\eps^2\Nabla{s}v$.
Then the first term has the same form as the one dicussed above. 
In the second and fourth term we estimate $\eps|\p_tu|$ 
by $2\delta$ and we use~$(V0)$. For the last term we find, by~$(V2)$,
\begin{align*}
\eps^2\Abs{\inner{\Nabla{s}v}{\Nabla{t}\Nabla{s}\grad\Vv(u)}}
&
\le \eps^2C\Abs{\Nabla{s}v}
\Bigl(
\Abs{\Nabla{t}\p_su} + \Abs{\p_su}+ \Norm{\p_su}_{L^2(S^1)} \\
&\qquad\qquad\qquad
+ \Abs{\p_tu}\bigl(\Abs{\p_su}+ \Norm{\p_su}_{L^2(S^1)}\bigr) 
\Bigr)\\
&
\le \eps^2C\Abs{\Nabla{s}v}
\Bigl(
\Abs{\Nabla{t}\p_su}
+ \Abs{\p_su}+ \Norm{\p_su}_{L^2(S^1)}\Bigr) \\
&\quad
+ 2\eps C\delta\Abs{\Nabla{s}v}
\left(\Abs{\p_su}+ \Norm{\p_su}_{L^2(S^1)}\right) \\
&\le
\mu f + \frac18\left(g + \Norm{\p_su}_{L^2(S^1)}^2\right).
\end{align*}
This leaves the terms three and five. 
In the third term we estimate $\eps^2|\p_tu|^2$
by $4\delta^2$ and use the identity 
$$
\Nabla{s}v=\Nabla{s}\p_tu+\eps^2\Nabla{s}\Nabla{s}v
$$ 
of~(\ref{eq:uv}).  For term five we use the identity
$$
\Nabla{t}\p_tu
= \Nabla{t}(v-\eps^2\Nabla{s}v)
=\p_su-\grad\Vv(u) - \eps^2\Nabla{t}\Nabla{s}v
$$
to obtain the expression
\begin{equation*}
\begin{split}
&\eps^2\inner{\Nabla{s}v}{R(\p_su,\p_su-\grad\Vv(u)
-\eps^2\Nabla{t}\Nabla{s}v)v} \\
&= -\eps^2\inner{\Nabla{s}v}{R(\p_su,\grad\Vv(u))v}
-\eps^4\inner{\Nabla{s}v}{R(\p_su,\Nabla{t}\Nabla{s}v)v}.
\end{split}
\end{equation*}
In the last summand we use the estimate $\eps^2|\p_su|\le\delta$.
This proves~(\ref{eq:UV}) for $0<\eps\le\eps_0$.
For $\eps_0\le\eps\le1$ the estimate~(\ref{eq:UV})
follows immediately from~(\ref{eq:Lf-1}), (\ref{eq:Lf-2}),
and Lemma~\ref{le:weak-L-infty}~(ii).
\end{proof}

The third step in the proof of Theorem~\ref{thm:gradient} is
to estimate the summand $\eps^4\Norm{\Nabla{s}\p_su}_{L^2(S^1)}^2$
in~(\ref{eq:UV}) in terms of the energy. This is the content of the 
following lemma. 

\begin{lemma}\label{le:gradient-L22}
Fix a constant $c_0>0$ and a perturbation
$\Vv:\Ll M\to\R$ that satisfies~$(V0-V3)$.  
Then there is a constant $C=C(c_0,\Vv)>0$ 
such that the following holds. If $0<\eps\le1$ and 
$(u,v):\R\times S^1\to TM$ is a solution 
of~(\ref{eq:floer-V}) that satisfies~(\ref{eq:EA})
then, for every $s\in\R$, 
\begin{equation}\label{eq:uvE-L22}
\begin{split}
&\int_0^1\Bigl(
\eps^2\Abs{\Nabla{t}\p_su}^2 
+ \eps^4\Abs{\Nabla{s}\p_su}^2
+ \eps^4\Abs{\Nabla{t}\Nabla{s}v}^2 
+ \eps^6\Abs{\Nabla{s}\Nabla{s}v}^2
\Bigl)  \\
&\le CE_{[s-1/2,s+1/2]}^\eps(u,v).
\end{split}
\end{equation}
\end{lemma}

\begin{proof}
Define $f_1$ and $g_1$ by
$$
2f_1
:= \Abs{\p_su}^2
+\eps^2\Abs{\Nabla{s}v}^2
+ \eps^2\Abs{\Nabla{t}\p_su}^2 
+ \eps^4\Abs{\Nabla{s}\p_su}^2
+ \eps^4\Abs{\Nabla{t}\Nabla{s}v}^2,
$$
\begin{align*}
2g_1 &:= 
\Abs{\Nabla{t}\p_su}^2 
+ \eps^2\Abs{\Nabla{s}\p_su}^2
+ \eps^2\Abs{\Nabla{t}\Nabla{s}v}^2
+ \eps^4\Abs{\Nabla{s}\Nabla{s}v}^2 \\
&\quad
+ \eps^2\Abs{\Nabla{t}\Nabla{t}\p_su}^2 
+ \eps^4\Abs{\Nabla{s}\Nabla{t}\p_su}^2 
+ \eps^4\Abs{\Nabla{t}\Nabla{s}\p_su}^2 \\
&\quad
+ \eps^6\Abs{\Nabla{s}\Nabla{s}\p_su}^2 
+ \eps^4\Abs{\Nabla{t}\Nabla{t}\Nabla{s}v}^2 
+ \eps^6\Abs{\Nabla{s}\Nabla{t}\Nabla{s}v}^2 
\end{align*}
and abbreviate $F_1(s):=\int_0^1f_1(s,t)\,dt$ and
$G_1(s):=\int_0^1g_1(s,t)\,dt$. Then 
\begin{equation}\label{eq:Lfst}
L_\eps f_1 = 2g_1 + U + \eps^2 V
+ \eps^2U_t + \eps^4U_s + \eps^4V_t 
\end{equation}
where $U:=\inner{\p_su}{\Ll_\eps\p_su}$ and 
$V:=\inner{\Nabla{s}v}{\Ll_\eps\Nabla{s}v}$ 
as in Lemma~\ref{le:gradient-L2} and
$$
U_t:=\inner{\Nabla{t}\p_su}{\Ll_\eps\Nabla{t}\p_su},\qquad
U_s:=\inner{\Nabla{s}\p_su}{\Ll_\eps\Nabla{s}\p_su},
$$
$$
V_t:=\inner{\Nabla{t}\Nabla{s}v}{\Ll_\eps\Nabla{t}\Nabla{s}v}.
$$
We shall prove the estimate
\begin{equation}\label{eq:UVst}
\begin{split}
&\Abs{U}+\eps^2\Abs{V}+\eps^2\Abs{U_t} + \eps^4\Abs{U_s+V_t} \\
&\le \mu f_1 
+ \frac12\left(g_1 
+ \Norm{\p_su}_{L^2(S^1)}^2
+ \eps^4\Norm{\Nabla{s}\p_su}_{L^2(S^1)}^2
+ \eps^8\Norm{\Nabla{s}\Nabla{s}\p_su}_{L^2(S^1)}^2
\right) \\
&\le \mu f_1 + F_1 + g_1 + G_1 
\end{split}
\end{equation}
for a suitable constant $\mu>0$.
By~(\ref{eq:Lfst}) and~(\ref{eq:UVst}), 
$
L_\eps f_1 + \mu f_1 +F_1 \ge  g_1 - G_1.
$
Integrating this inequality over the interval $0\le t\le1$
gives 
$$
\eps^2F_1''-F_1'+(\mu+1)F_1\ge 0.
$$
Hence it follows from Lemma~\ref{le:apriori-1eps} 
with $r:=1/5$ that
$$
F_1(s) \le c\int_{s-1/4}^{s+1/4}F_1(\sigma)\,d\sigma
\le c\left(1+\frac{C\eps^2}{2}\right)
E_{[s-1/2,s+1/2]}^\eps(u,v).
$$
Here $c:=250 c_2e^{(\mu+1)/25}$, where $c_2$ is the constant of
Lemma~\ref{le:apriori-1eps}, and the second inequality
follows from Lemma~\ref{le:gradient-L2}.  Now use 
Lemma~\ref{le:gradient-L2} again and the identity 
$\eps^2\Nabla{s}\Nabla{s}v=\Nabla{s}v-\Nabla{s}\p_tu$
to estimate the term $\eps^6\Abs{\Nabla{s}\Nabla{s}v}^2$. 
 
It remains to prove~(\ref{eq:UVst}). 
For the terms $\Abs{U}+\eps^2\Abs{V}$ the estimate 
was established in~(\ref{eq:UV}).
To estimate the term $\eps^2\Abs{U_t}$ write
\begin{equation}\label{eq:Ut}
\begin{split}
\Ll_\eps\Nabla{t}\p_su
&= \Nabla{t}\Ll_\eps\p_su
+ \eps^2[\Nabla{s}\Nabla{s},\Nabla{t}]\p_su 
- [\Nabla{s},\Nabla{t}]\p_su \\
&= \Nabla{t}\Ll_\eps\p_su
+ \eps^2\Nabla{s}\left(R(\p_su,\p_tu)\p_su)\right)  \\
&\quad
+ \eps^2R(\p_su,\p_tu)\Nabla{s}\p_su
- R(\p_su,\p_tu)\p_su \\
&=2\eps^2\Nabla{t}\left(R(\p_su,\p_tu)\Nabla{s}v\right)
+ \eps^2\Nabla{t}\left(\left(\Nabla{\p_su}R\right)(\p_su,\p_tu)v\right) \\
&\quad
- \Nabla{t}\left(R(\p_su,\p_tu)v\right)   \\
&\quad
+\eps^2\Nabla{t}\left(R(\Nabla{s}\p_su,\p_tu)v\right)
+\eps^2\Nabla{t}\left(R(\p_su,\Nabla{s}\p_tu)v\right) \\
&\quad
-\Nabla{t}\Nabla{s}\grad\Vv(u)
+\eps^2\Nabla{t}\Nabla{s}\Nabla{s}\grad\Vv(u) \\
&\quad
+ \eps^2\Nabla{s}\left(R(\p_su,\p_tu)\p_su)\right)  \\
&\quad
+ \eps^2R(\p_su,\p_tu)\Nabla{s}\p_su
- R(\p_su,\p_tu)\p_su. 
\end{split}
\end{equation}
The last equation follows from~(\ref{eq:Lf-1}). 
Now take the pointwise inner product with $\eps^2\Nabla{t}\p_su$.
We begin by explaining how to estimate the first term. 
We encounter an expression of the form 
$\eps^4\inner{\Nabla{t}\p_su}{(\Nabla{\p_tu}R)(\p_su,\p_tu)\Nabla{s}v}$.
Here we can use the identity 
$$
\Nabla{s}v = \Nabla{s}\p_tu + \eps^2\Nabla{s}\Nabla{s}v
$$
to obtain an inequality
$$
\Abs{\Nabla{s}v}\Abs{\Nabla{t}\p_su}\le 3g_1
$$
By~(\ref{eq:uvs})
we can estimate the product $\eps^4\Abs{\p_su}\Abs{\p_tu}^2$
by a small constant. Another expression we encounter is 
$\eps^4\inner{\Nabla{t}\p_su}{R(\Nabla{t}\p_su,\p_tu)\Nabla{s}v}$;
by~(\ref{eq:uvs}), we have 
$\eps^4\Abs{\p_tu}\Abs{\Nabla{s}v}\le 2\delta^2$
and so the expression can be estimated by a small constant times $g_1$. 
Then we encounter the expression
$\eps^4\inner{\Nabla{t}\p_su}{R(\p_su,\Nabla{t}\p_tu)\Nabla{s}v}$;
here we use the identity
$$
\Nabla{t}\p_tu = \Nabla{t}(v-\eps^2\Nabla{s}v)
= \p_su - \grad\Vv(u) - \eps^2\Nabla{t}\Nabla{s}v;
$$
the crucial observation is that the summand $\p_su$
can be dropped when inserting this formula in
$R(\p_su,\Nabla{t}\p_tu)$; in the summand
$\eps^4\inner{\Nabla{t}\p_su}{R(\p_su,\grad\Vv(u))\Nabla{s}v}$
we use~$(V0)$ and $\eps^2\Abs{\p_su}\le\delta$; 
for the summand 
$\eps^6\inner{\Nabla{t}\p_su}{R(\p_su,\Nabla{t}\Nabla{s}v)\Nabla{s}v}$
we use $\eps^5\Abs{\p_su}\Abs{\Nabla{s}v}\le\delta^2$
and $\eps\Abs{\Nabla{t}\p_su}\Abs{\Nabla{t}\Nabla{s}v}\le Cg_1$.
The last expression we encounter is 
$\eps^4\inner{\Nabla{t}\p_su}{R(\p_su,\p_tu)\Nabla{t}\Nabla{s}v}$;
here we use $\eps^3\Abs{\p_su}\Abs{\p_tu}\le2\delta^2$, by~(\ref{eq:uvs}),
and again $\eps\Abs{\Nabla{t}\p_su}\Abs{\Nabla{t}\Nabla{s}v}\le Cg_1$.
This deals with the first term; the next two terms can be estimated
by the same method.

In the fourth term we encounter the expression
$\eps^4\inner{\Nabla{t}\p_su}{R(\Nabla{t}\Nabla{s}\p_su,\p_tu)v}$;
here we use $\eps\Abs{\p_tu}\le2\delta$ and
$\eps^2\Abs{\Nabla{t}\p_su}\Abs{\Nabla{t}\Nabla{s}\p_su}\le Cg_1$.
Another expression is 
$\eps^4\inner{\Nabla{t}\p_su}{R(\Nabla{s}\p_su,\p_tu)\Nabla{t}v}$;
here we use $\Nabla{t}v=\p_su-\grad\Vv(u)$ and the inequalities
$\eps^3\Abs{\p_su}\Abs{\p_tu}\le2\delta^2$ and
$\eps\Abs{\Nabla{t}\p_su}\Abs{\Nabla{s}\p_su}\le Cg_1$.
A third expression is
$\eps^4\inner{\Nabla{t}\p_su}{R(\Nabla{s}\p_su,\Nabla{t}\p_tu)v}$;
here we use the formula
\begin{equation}\label{eq:Lu}
\begin{split}
\eps^2\Nabla{s}\p_su + \Nabla{t}\p_tu
&= \eps^2\Nabla{s}(\Nabla{t}v+\grad\Vv(u))
+ \Nabla{t}(v-\eps^2\Nabla{s}v)  \\
&= \p_su + \eps^2R(\p_su,\p_tu)v
   - \grad\Vv(u) + \eps^2\Nabla{s}\grad\Vv(u);
\end{split}
\end{equation}
so the curvature term can be estimated by 
\begin{equation}\label{eq:sstt}
\Abs{R(\Nabla{s}\p_su,\Nabla{t}\p_tu)}
\le C\Abs{\Nabla{s}\p_su}\left(1+\Abs{\p_su} 
+ \eps^2\Abs{\p_su}\Abs{\p_tu}\right).
\end{equation}
This completes the discussion of the fourth term.
The fifth term is similar, except that the 
cubic expression in the second derivatives vanishes.
The last three terms can be disposed off similarly;
the only new expression that appears is 
$\eps^4\inner{\Nabla{t}\p_su}{(\Nabla{\p_su}R)(\p_su,\p_tu)\p_su}$;
here we use $\p_tu=v-\eps^2\Nabla{s}v$ and the inequalities
$\eps^2\Abs{\p_su}\le\delta$ as well as
$\Abs{\Nabla{t}\p_su}\Abs{\p_su}\le g_1+f_1$ and
$\Abs{\Nabla{t}\p_su}\Abs{\Nabla{s}v}\le 3g_1$.

This leaves the terms involving $\grad\Vv$.
For $\eps^2\inner{\Nabla{t}\p_su}{\Nabla{t}\Nabla{s}\grad\Vv(u)}$
we use~$(V2)$ and for
$\eps^4\inner{\Nabla{t}\p_su}{\Nabla{t}\Nabla{s}\Nabla{s}\grad\Vv(u)}$
we use~$(V3)$. Both terms can be estimated by 
$C\eps(f_1+g_1+\Norm{\p_su}_{L^2(S^1)}^2
+\eps^4\Norm{\Nabla{s}\p_su}_{L^2(S^1)}^2)$.  
This completes the estimate of $\eps^2\Abs{U_t}$. 

To estimate the term $\eps^4\Abs{U_s+V_t}$ write
\begin{equation}\label{eq:Us}
\begin{split}
\Ll_\eps\Nabla{s}\p_su
&= \Nabla{s}\Ll_\eps\p_su
+ [\Nabla{t}\Nabla{t},\Nabla{s}]\p_su \\
&= \Nabla{s}\Ll_\eps\p_su \\
&\quad
- \Nabla{t}\left(R(\p_su,\p_tu)\p_su\right) 
- R(\p_su,\p_tu)\Nabla{t}\p_su \\
&=2\eps^2\Nabla{s}\left(R(\p_su,\p_tu)\Nabla{s}v\right)
+ \eps^2\Nabla{s}\left(\left(\Nabla{\p_su}R\right)(\p_su,\p_tu)v\right) \\
&\quad
- \Nabla{s}\left(R(\p_su,\p_tu)v\right)   \\
&\quad
+\eps^2\Nabla{s}\left(R(\Nabla{s}\p_su,\p_tu)v\right)
+\eps^2\Nabla{s}\left(R(\p_su,\Nabla{s}\p_tu)v\right)  \\
&\quad
-\Nabla{s}\Nabla{s}\grad\Vv(u)
+\eps^2\Nabla{s}\Nabla{s}\Nabla{s}\grad\Vv(u) \\
&\quad
- \Nabla{t}\left(R(\p_su,\p_tu)\p_su\right) 
- R(\p_su,\p_tu)\Nabla{t}\p_su
\end{split}
\end{equation}
(where the last equation follows from~(\ref{eq:Lf-1}))
and
\begin{equation}\label{eq:Vt}
\begin{split}
\Ll_\eps\Nabla{t}\Nabla{s}v
&= \Nabla{t}\Ll_\eps\Nabla{s}v 
+ \eps^2[\Nabla{s}\Nabla{s},\Nabla{t}]\Nabla{s}v
- [\Nabla{s},\Nabla{t}]\Nabla{s}v \\
&= \Nabla{t}\Ll_\eps\Nabla{s}v
+ \eps^2\Nabla{s}\left(R(\p_su,\p_tu)\Nabla{s}v)\right)  \\
&\quad
+ \eps^2R(\p_su,\p_tu)\Nabla{s}\Nabla{s}v
- R(\p_su,\p_tu)\Nabla{s}v \\
&=-2\Nabla{t}\left(R(\p_su,\p_tu)\p_su\right)
+ \Nabla{t}\left(R(\p_su,\p_tu)\grad\Vv(u)\right) \\
&\quad
-\Nabla{t}\left(\left(\Nabla{\p_tu}R\right)(\p_su,\p_tu)v\right) \\
&\quad
-\Nabla{t}\left(R(\Nabla{t}\p_su,\p_tu)v\right)
-\Nabla{t}\left(R(\p_su,\Nabla{t}\p_tu)v\right) \\
&\quad
- \Nabla{t}\Nabla{t}\Nabla{s}\grad\Vv(u) \\
&\quad
+ \eps^2\Nabla{s}\left(R(\p_su,\p_tu)\Nabla{s}v)\right)  \\
&\quad
+ \eps^2R(\p_su,\p_tu)\Nabla{s}\Nabla{s}v
- R(\p_su,\p_tu)\Nabla{s}v
\end{split}
\end{equation}
(where the last equation follows from~(\ref{eq:Lf-2})).
The terms that require special attention are those involving $\grad\Vv$
and the cubic terms in the second derivatives.
The cubic terms in the second derivatives are
\begin{equation*}
\begin{split}
U_{s0} &:= 2\eps^6\inner{\Nabla{s}\p_su}
{R(\Nabla{s}\p_su,\Nabla{s}\p_tu)v}, \\
V_{t0} &:= 2\eps^4\inner{\Nabla{t}\Nabla{s}v}
{R(\Nabla{t}\p_tu,\Nabla{t}\p_su)v}.
\end{split}
\end{equation*}
Now insert 
$$
\Nabla{s}\p_su = \Nabla{s}\left(\Nabla{t}v+\grad\Vv(u)\right),\qquad
\Nabla{t}\p_tu = \Nabla{t}\left(v-\eps^2\Nabla{s}v\right)
$$
into $U_{s0}$ and $V_{t0}$, respectively. Then the only difficult 
remaining terms are the ones involving again three second derivatives.
After replacing $\Nabla{s}\Nabla{t}v$ by 
$\Nabla{t}\Nabla{s}v+R(\p_su,\p_tu)v$ we obtain
\begin{equation*}
\begin{split}
U_{s1} &:= 2\eps^6\inner{\Nabla{s}\p_su}
{R(\Nabla{t}\Nabla{s}v,\Nabla{s}\p_tu)v},\\
V_{t1} &:= -2\eps^6\inner{\Nabla{t}\Nabla{s}v}
{R(\Nabla{t}\Nabla{s}v,\Nabla{t}\p_su)v}.
\end{split}
\end{equation*}
The sum is
\begin{equation*}
\begin{split}
U_{s1}+V_{t1}
&= 2\eps^6\inner{\Nabla{s}\p_su-\Nabla{t}\Nabla{s}v}
{R(\Nabla{t}\Nabla{s}v,\Nabla{s}\p_tu)v} \\
&= 2\eps^6\inner{\Nabla{s}(\p_su-\Nabla{t}v) + R(\p_su,\p_tu)v}
{R(\Nabla{t}\Nabla{s}v,\Nabla{s}\p_tu)v} \\
&= 2\eps^6\inner{\Nabla{s}\grad\Vv(u) + R(\p_su,\p_tu)v}
{R(\Nabla{t}\Nabla{s}v,\Nabla{s}\p_tu)v}
\end{split}
\end{equation*}
and can be estimated in the required fashion.

The terms involving $\grad\Vv$ can be estimated by 
\begin{align*}
&\eps^6\Abs{\inner{\Nabla{s}\p_su}
{\Nabla{s}\Nabla{s}\Nabla{s}\grad\Vv(u)}} 
+ \eps^4\Abs{\inner{\Nabla{s}\p_su}
{\Nabla{s}\Nabla{s}\grad\Vv(u)}}  \\
&
+ \eps^4\Abs{\inner{\Nabla{t}\Nabla{s}v}
{R(\p_su,\p_tu)\Nabla{t}\grad\Vv(u)}}
+ \eps^4\Abs{\inner{\Nabla{t}\Nabla{s}v}
{\Nabla{t}\Nabla{t}\Nabla{s}\grad\Vv(u)}} \\
&\le
C\eps^2\Abs{\Nabla{s}\p_su}
\Bigl(\eps^4\Abs{\Nabla{s}\Nabla{s}\p_su}
+\eps^2\Abs{\Nabla{s}\p_su} 
+ \Abs{\p_su}\Bigr) \\
&\quad
+ C\eps^2\Abs{\Nabla{s}\p_su}
\Bigl(\eps^4\Norm{\Nabla{s}\Nabla{s}\p_su}_{L^2(S^1)}
+\eps^2\Norm{\Nabla{s}\p_su}_{L^2(S^1)}
+\Norm{\p_su}_{L^2(S^1)}
\Bigr) \\
&\quad
+ C\eps^2\Abs{\Nabla{t}\Nabla{s}v}
\Bigl(\eps^2\Abs{\Nabla{t}\Nabla{t}\p_su}
+ \eps\Abs{\Nabla{t}\p_su}
+ \Abs{\p_su}
+ \Norm{\p_su}_{L^2(S^1)}
\Bigr) \\
&\quad
+ C\eps^4\Abs{\Nabla{t}\Nabla{s}v}^2 \\
&\le 
\mu f_1+
\frac18\left(g_1+\Norm{\p_su}_{L^2(S^1)}^2
+\eps^4\Norm{\Nabla{s}\p_su}_{L^2(S^1)}^2  
+\eps^8\Norm{\Nabla{s}\Nabla{s}\p_su}_{L^2(S^1)}^2\right).
\end{align*}
Here the first inequality follows from~$(V1-3)$;
it also uses the identity 
$\Nabla{t}\p_tu=\p_su-\Nabla{s}\grad\Vv(u)
-\eps^2\Nabla{t}\Nabla{s}v$ 
and the fact that $\eps^2\Abs{\p_su}$ and $\eps\Abs{\p_tu}$
are uniformly bounded (Lemma~\ref{le:weak-L-infty}).
All the other summands appearing in our expression for
$\eps^4\Abs{U_s+V_t}$ can be estimated by the same arguments 
as for $\eps^2\Abs{U_t}$. This implies~(\ref{eq:UVst}) 
and completes the proof of Lemma~\ref{le:gradient-L22}.
\end{proof}

The fourth step in the proof of Theorem~\ref{thm:gradient}
is to establish the $L^\infty$ estimate with $\Nabla{s}v$ 
replaced by $\eps\Nabla{s}v$. 

\begin{lemma}\label{le:gradient-eps}
Fix a constant $c_0>0$ and a perturbation
$\Vv:\Ll M\to\R$ that satisfies~$(V0-V3)$.  
Then there is a constant $C=C(c_0,\Vv)>0$ 
such that the following holds. If $0<\eps\le1$ and 
$(u,v):\R\times S^1\to TM$ is a solution 
of~(\ref{eq:floer-V}) that satisfies~(\ref{eq:EA}),
i.e. $E^\eps(u,v)\le c_0$ and 
$\sup_{s\in\R}\Aa_\Vv(u(s,\cdot),v(s,\cdot))\le c_0$,
then 
\begin{equation}\label{eq:uvepsE}
\Abs{\p_su(s,t)}^2+\eps^2\Abs{\Nabla{s}v(s,t)}^2
\le CE_{[s-1,s+1]}^\eps(u,v)
\le Cc_0
\end{equation}
for all $s$ and $t$. 
\end{lemma}

\begin{proof}
Let $f$, $g$, $F$, $G$, $U$, $V$, and $\mu$ be as 
in the proof of Lemma~\ref{le:gradient-L2}. 
Choose a constant $C>0$ such that the assertions of 
Lemmas~\ref{le:gradient-L2} and~\ref{le:gradient-L22} 
hold with this constant. Then, by~(\ref{eq:Lfg}) 
and~(\ref{eq:UV}), we have
\begin{align*}
L_\eps f
&= 2g + U + \eps^2V \\
&\ge -\mu f - \frac12\left(
\Norm{\p_su}_{L^2(S^1)}^2+\eps^4\Norm{\Nabla{s}\p_su}_{L^2(S^1)}^2
\right) \\
&\ge -\mu f - CE_{[s-1/2,s+1/2]}^\eps(u,v)
\end{align*}
for all $(s,t)\in\R\times S^1$.  Let $s_0\in\R$ 
and denote 
$
a:= \frac{C}{\mu}
E_{[s_0-1,s_0+1]}^\eps(u,v).
$
Then 
$
L_\eps(f+a) + \mu(f+a)\ge 0
$
for $s_0-1/2\le s\le s_0+1/2$. 
Hence we may apply Lemma~\ref{le:apriori-basic-eps}
with $r=1/3$ to the function $w(s,t):=f(s_0+s,t_0+t)+a$:
\begin{equation*}
\begin{split}
     f(s_0,t_0)
    &\le54c_2e^{\mu/9}
     \int_{s_0-1/9-\eps/3}^{s_0+\eps/3}
     \int_0^1(f(s,t)+a)\,dtds \\
    &\le54c_2e^{\mu/9}
     \int_{s_0-1/2}^{s_0+1/2}\int_0^1
     \left(\frac12\Abs{\p_su(s,t)}^2
     +\frac{\eps^2}{2}\Abs{\Nabla{s}v(s,t)}^2
     + a\right)\,dtds \\
    &\le54c_2e^{\mu/9}
     \Bigl(E_{[s_0-1,s_0+1]}^\eps(u,v)+a\Bigr) \\
    &=54c_2e^{\mu/9}\left(1 + \frac{C}{\mu}\right)
      E_{[s_0-1,s_0+1]}^\eps(u,v). 
\end{split}
\end{equation*}
This proves the lemma. 
\end{proof}

\begin{proof}[Proof of Theorem~\ref{thm:gradient}.]
Define $f_2$ and $g_2$ by
\begin{align*}
2f_2 &:= \Abs{\p_su}^2 + \Abs{\Nabla{s}v}^2,\\
2g_2 &:= 
\Abs{\Nabla{t}\p_su}^2 
+ \eps^2\Abs{\Nabla{s}\p_su}^2
+ \Abs{\Nabla{t}\Nabla{s}v}^2
+ \eps^2\Abs{\Nabla{s}\Nabla{s}v}^2 
\end{align*}
and abbreviate 
$
F_2(s):=\int_0^1f_2(s,t)\,dt
$
and
$
G_2(s):=\int_0^1g_2(s,t)\,dt.
$
Then
\begin{equation}\label{eq:Lf2}
L_\eps f_2 = 2g_2 + U + V
\end{equation}
where $U$ and $V$ are as in Lemma~\ref{le:gradient-L2}.
These functions satisfy the estimate
\begin{equation}\label{eq:UV2}
\Abs{U}+\Abs{V} \le \mu f_2 
+ \frac12\left(g_2
+ \Norm{\p_su}_{L^2}^2
+ \eps^4\Norm{\Nabla{s}\p_su}_{L^2}^2\right)
\end{equation}
for a suitable constant $\mu>0$; here $\Norm{\cdot}_{L^2}$
denotes the $L^2$-norm over the circle at time $s$.
This follows from~(\ref{eq:Lf-1}) and~(\ref{eq:Lf-2}) 
via term by term inspection. (We use the fact that
$\Abs{\p_su}$, $\eps\Abs{\Nabla{s}v}$, and 
$\Abs{\p_tu}$ are uniformly bounded, 
by Lemma~\ref{le:gradient-eps}.)  

By Lemmas~\ref{le:gradient-L2} 
and~\ref{le:gradient-L22}, we have
$$
\int_0^1\left(\Abs{\p_su(s,t)}^2
+ \eps^4\Abs{\Nabla{s}\p_su(s,t)}^2\right)\,dt 
\le CE_{[s-1/2,s+1/2]}^\eps(u,v)
$$
for a suitable constant $C$ and every $s\in\R$.
Hence it follows from~(\ref{eq:Lf2})
and~(\ref{eq:UV2}) that
$$
L_\eps f_2(s,t) \ge -\mu f_2(s,t)
-CE_{[s-1/2,s+1/2]}^\eps(u,v)
$$
for all $(s,t)\in\R\times S^1$. Fix a number 
$s_0$ and abbreviate
$
a:= \frac{C}{\mu}
E_{[s_0-1,s_0+1]}^\eps(u,v).
$
Then 
$
L_\eps(f_2+a) + \mu(f_2+a)\ge 0
$
for $s_0-1/2\le s\le s_0+1/2$. 
Hence we may apply Lemma~\ref{le:apriori-basic-eps}
with $r=1/3$ to the function $w(s,t):=f_2(s_0+s,t_0+t)+a$:
\begin{equation*}
\begin{split}
     f_2(s_0,t_0)
    &\le54c_2e^{\mu/9}
     \int_{s_0-1/9-\eps/3}^{s_0+\eps/3}
     \int_0^1(f_2(s,t)+a)\,dtds \\
    &\le54c_2e^{\mu/9}
     \int_{s_0-1/2}^{s_0+1/2}\int_0^1
     \left(\frac12\Abs{\p_su(s,t)}^2
     +\frac{1}{2}\Abs{\Nabla{s}v(s,t)}^2
     + a\right)\,dtds \\
    &\le c_3
     \Bigl(E_{[s_0-1,s_0+1]}^\eps(u,v)+a\Bigr) \\
    &= c_3\left(1 + \frac{C}{\mu}\right)
      E_{[s_0-1,s_0+1]}^\eps(u,v). 
\end{split}
\end{equation*}
Here  the third inequality, with a suitable constant $c_3$,  
follows from Corollary~\ref{cor:gradient-L2}. 
This proves the pointwise estimate. 

To prove the $L^2$-estimate integrate~(\ref{eq:Lf2}) 
and~(\ref{eq:UV2}) over $0\le t\le1$ to obtain
$$
\eps^2F_2''-F_2'+(\mu+1)F_2\ge G_2
$$
for every $s\in\R$. Hence, by Lemma~\ref{le:apriori-1-L2}
with suitable choices of $R$ and $r$, we have 
$$
\int_{-1/2}^{1/2}G_2(s)\,ds
\le c_4\int_{-3/4}^{3/4}F_2(s)\,ds
$$
for every $s\in\R$ and a constant $c_4>0$ that depends
only on $R$, $r$, and $\mu$. Now it follows from
Corollary~\ref{cor:gradient-L2} that
$$
\int_{-3/4}^{3/4}F_2(s)\,ds\le c_5E^\eps_{[s-1,s+1]}(u,v)
$$
for every $s>0$ and some constant $c_5=c_5(c_0,\Vv)>0$. 
Hence
$$
\int_{-1/2}^{1/2}\int_0^1\left(
\Abs{\Nabla{t}\p_su}^2 
+ \Abs{\Nabla{t}\Nabla{s}v}^2
+ \eps^2\Abs{\Nabla{s}\Nabla{s}v}^2 
\right)\,dtds
\le 2c_4c_5E^\eps_{[s-1,s+1]}(u,v).
$$
The estimate for $\Nabla{s}\p_su$ now follows from the identity
$$
\Nabla{s}\p_su 
= \Nabla{s}\Nabla{t}v+\Nabla{s}\grad\Vv(u)
= \Nabla{t}\Nabla{s}v + R(\p_su,\p_tu)v + \Nabla{s}\grad\Vv(u). 
$$ 
This proves Theorem~\ref{thm:gradient}. 
\end{proof}


\section{Estimates of the second derivatives}\label{sec:2nd-derivs}

\begin{theorem}\label{thm:2nd-derivs}
Fix a constant $c_0>0$ and a perturbation
$\Vv:\Ll M\to\R$ that satisfies~$(V0-V4)$.  
Then there is a constant $C=C(c_0,\Vv)>0$ 
such that the following holds. If $0<\eps\le1$ 
and $(u,v):\R\times S^1\to TM$ is a solution
of~(\ref{eq:floer-V}) that satisfies~(\ref{eq:EA})
then
\begin{equation}\label{eq:2D}
\begin{split}
&\Norm{\Nabla{t}\p_su}_{L^p([-T,T]\times S^1)}
+\Norm{\Nabla{s}\p_su}_{L^p([-T,T]\times S^1)} \\
&+\Norm{\Nabla{t}\Nabla{s}v}_{L^p([-T,T]\times S^1)}
+\Norm{\Nabla{s}\Nabla{s}v}_{L^p([-T,T]\times S^1)} \\
&\le c\sqrt{E_{[-T-1,T+1]}^\eps(u,v)}
\end{split}
\end{equation}
for $T>1$ and $2\le p\le\infty$.
\end{theorem}

For $p=2$ the estimate, with $\Nabla{s}\Nabla{s}v$ replaced by 
$\eps\Nabla{s}\Nabla{s}v$,  was established in Theorem~\ref{thm:gradient}.   
The strategy is to prove the estimate for $p=\infty$ and, as a byproduct, 
to get rid of the factor $\eps$ for $p=2$ (see Corollary~\ref{cor:2D} below).
The result for general $p$ then follows by interpolation. 

\begin{lemma}\label{le:2D}
Fix a constant $c_0>0$ and a perturbation
$\Vv:\Ll M\to\R$ that satisfies~$(V0-V3)$.  
Then there is a constant $C=C(c_0,\Vv)>0$ 
such that the following holds. If $0<\eps\le1$ and 
$(u,v):\R\times S^1\to TM$ is a solution 
of~(\ref{eq:floer-V}) that satisfies~(\ref{eq:EA}),
then 
\begin{align*}
&\int_0^1
\Bigl(
\Abs{\Nabla{t}\p_su(s,t)}^2
+\Abs{\Nabla{s}\p_su(s,t)}^2
+\Abs{\Nabla{t}\Nabla{s}v(s,t)}^2
+\eps^2\Abs{\Nabla{s}\Nabla{s}v(s,t)}^2
\Bigr)\,dt \\
&+\int_{s-1/4}^{s+1/4}\int_0^1\Bigl(
\Abs{\Nabla{t}\Nabla{t}\p_su}^2 
+ \Abs{\Nabla{t}\Nabla{s}\p_su}^2 
+ \eps^2\Abs{\Nabla{s}\Nabla{s}\p_su}^2
\Bigl)  \\
&+\int_{s-1/4}^{s+1/4}\int_0^1\Bigl(
\Abs{\Nabla{t}\Nabla{t}\Nabla{s}v}^2 
+ \eps^2\Abs{\Nabla{t}\Nabla{s}\Nabla{s}v}^2 
+ \eps^4\Abs{\Nabla{s}\Nabla{s}\Nabla{s}v}^2
\Bigl)  \\
&\le CE_{[s-1/2,s+1/2]}^\eps(u,v)
\end{align*}
for every $s\in\R$. 
\end{lemma}

\begin{corollary}\label{cor:2D}
Fix a constant $c_0>0$ and a perturbation
$\Vv:\Ll M\to\R$ that satisfies~$(V0-V3)$.  
Then there is a constant $C=C(c_0,\Vv)>0$ 
such that the following holds. If $0<\eps\le1$ and 
$(u,v):\R\times S^1\to TM$ is a solution 
of~(\ref{eq:floer-V}) that satisfies~(\ref{eq:EA}),
then 
$$
\int_{s-1/4}^{s+1/4}\int_0^1\Abs{\Nabla{s}\Nabla{s}v}^2
\le CE_{[s-1/2,s+1/2]}^\eps(u,v)
$$
for every $s\in\R$. 
\end{corollary}

\begin{proof}
Since 
$$
\Nabla{s}\Nabla{s}v
=\Nabla{s}\Nabla{s}\p_tu+\eps^2\Nabla{s}\Nabla{s}\Nabla{s}v
=R(\p_su,\p_tu)\p_su+\Nabla{t}\Nabla{s}\p_su
+\eps^2\Nabla{s}\Nabla{s}\Nabla{s}v
$$
this estimate follows immediately from Lemma~\ref{le:2D}.
\end{proof}

\begin{proof}[Proof of Lemma~\ref{le:2D}.]
Define $f_3$ and $g_3$ by
$$
2f_3 := \Abs{\p_su}^2 + \Abs{\Nabla{s}v}^2
+ \Abs{\Nabla{t}\p_su}^2
+ \Abs{\Nabla{s}\p_su}^2 
+ \Abs{\Nabla{t}\Nabla{s}v}^2
+ \eps^2\Abs{\Nabla{s}\Nabla{s}v}^2
$$
and
\begin{align*}
2g_3
&:= \Abs{\Nabla{t}\p_su}^2 
+ \eps^2\Abs{\Nabla{s}\p_su}^2
+ \Abs{\Nabla{t}\Nabla{s}v}^2
+ \eps^2\Abs{\Nabla{s}\Nabla{s}v}^2 \\
&\quad
+ \Abs{\Nabla{t}\Nabla{t}\p_su}^2 
+ \eps^2\Abs{\Nabla{s}\Nabla{t}\p_su}^2 
+ \Abs{\Nabla{t}\Nabla{s}\p_su}^2 
+ \eps^2\Abs{\Nabla{s}\Nabla{s}\p_su}^2  \\
&\quad
+ \Abs{\Nabla{t}\Nabla{t}\Nabla{s}v}^2 
+ \eps^2\Abs{\Nabla{s}\Nabla{t}\Nabla{s}v}^2
+ \eps^2\Abs{\Nabla{t}\Nabla{s}\Nabla{s}v}^2 
+ \eps^4\Abs{\Nabla{s}\Nabla{s}\Nabla{s}v}^2
\end{align*}
and abbreviate 
$$
F_3(s):=\int_0^1f_3(s,t)\,dt,\qquad
G_3(s):=\int_0^1g_3(s,t)\,dt. 
$$
Then
\begin{equation}\label{eq:Lf3}
L_\eps f_3 = 2g_3 + U + V + U_t + U_s + V_t + \eps^2V_s
\end{equation}
where $U$, $V$, $U_t$, $U_s$, $V_t$, and $V_s$ 
are as in Lemma~\ref{le:gradient-L22}.
These functions satisfy the estimate
\begin{equation}\label{eq:UV3}
\begin{split}
&\Abs{U}+\Abs{V} + \Abs{U_t} + \Abs{U_s+V_t} + \eps^2\Abs{V_s}  \\
&\le \mu f_3 + \frac12\left(g_3
+ \Norm{\p_su}_{L^2}^2
+ \Norm{\Nabla{s}\p_su}_{L^2}^2
+ \eps^4\Norm{\Nabla{s}\Nabla{s}\p_su}_{L^2}^2
\right)  \\
&\le \mu f_3 + F_3 + \frac{1}{2}(g_3 + G_3)
\end{split}
\end{equation}
for a suitable constant $\mu>0$; here $\Norm{\cdot}_{L^2}$
denotes the $L^2$-norm over the circle at time $s$.
For $U$ and $V$ this follows from~(\ref{eq:UV3}) in the 
proof of Theorem~\ref{thm:gradient}. For $U_t$ this follows
from~(\ref{eq:Ut}) and for $U_s+V_t$ from~(\ref{eq:Us}) 
and~(\ref{eq:Vt}) by the same arguments as in the proof
of Lemma~\ref{le:gradient-L22}.
The improved estimate~(\ref{eq:UV3}) follows by combining 
these arguments with Theorem~\ref{thm:gradient}. 
For $V_s$ we use the formula
\begin{equation}\label{eq:Vs}
\begin{split}
\Ll_\eps\Nabla{s}\Nabla{s}v
&= \Nabla{s}\Ll_\eps\Nabla{s}v 
+ [\Nabla{t}\Nabla{t},\Nabla{s}]\Nabla{s}v \\
&= \Nabla{s}\Ll_\eps\Nabla{s}v \\
&\quad
- \Nabla{t}\left(R(\p_su,\p_tu)\Nabla{s}v)\right) 
- R(\p_su,\p_tu)\Nabla{t}\Nabla{s}v \\
&=-2\Nabla{s}\left(R(\p_su,\p_tu)\p_su\right)
+ \Nabla{s}\left(R(\p_su,\p_tu)\grad\Vv(u)\right) \\
&\quad
-\Nabla{s}\left(\left(\Nabla{\p_tu}R\right)(\p_su,\p_tu)v\right) \\
&\quad
-\Nabla{s}\left(R(\Nabla{t}\p_su,\p_tu)v\right)
-\Nabla{s}\left(R(\p_su,\Nabla{t}\p_tu)v\right) \\
&\quad
- \Nabla{s}\Nabla{t}\Nabla{s}\grad\Vv(u) \\
&\quad
- \Nabla{t}\left(R(\p_su,\p_tu)\Nabla{s}v)\right) 
- R(\p_su,\p_tu)\Nabla{t}\Nabla{s}v.
\end{split}
\end{equation}
(The last equation uses~(\ref{eq:Lf-2}).) 
The desired estimate now follows from a term by term inspection;
since all the first derivatives are uniformly bounded, 
by Theorem~\ref{thm:gradient}, we only need to examine the 
second and third derivatives; in particular, the cubic term 
$\eps^2\inner{\Nabla{s}\Nabla{s}v}{R(\Nabla{s}\p_su,\Nabla{t}\p_tu)v}$
can be estimated by 
$C\eps^2\Abs{\Nabla{s}\Nabla{s}v}\Abs{\Nabla{s}\p_su}$
(see~(\ref{eq:sstt}) in the proof of Lemma~\ref{le:gradient-L22}).

It follows from~(\ref{eq:Lf3}) and~(\ref{eq:UV3}) that
$$
L_\eps f_3 + \mu f_3+ F_3 \ge g_3  + \frac12(g_3 - G_3).
$$
Integrating this inequality over the interval $0\le t\le 1$
gives 
$$
\eps^2F_3''-F_3'+ (\mu+1) F_3\ge G_3.
$$
By Theorem~\ref{thm:gradient} and Corollary~\ref{cor:gradient-L2} 
we have 
$$
\int_{s-1/4}^{s+1/4}F_3(s)\,ds
\le CE^\eps_{[s-1/2,s+1/2]}(u,v)
$$
for a suitable constant $C=C(c_0,\Vv)>0$. Hence the 
estimate for the second derivatives follows from 
Lemma~\ref{le:apriori-1eps} with $r:=1/5$.
The estimate for the third derivatives 
follows from Lemma~\ref{le:apriori-1-L2}.  
\end{proof}

\begin{lemma}\label{le:2Deps}
Fix a constant $c_0>0$ and a perturbation
$\Vv:\Ll M\to\R$ that satisfies~$(V0-V3)$.  
Then there is a constant $c=c(c_0,\Vv)>0$ 
such that the following holds. If $0<\eps\le1$ and 
$(u,v):\R\times S^1\to TM$ is a solution 
of~(\ref{eq:floer-V}) that satisfies~(\ref{eq:EA}),
then 
\begin{align*}
&\Norm{\Nabla{t}\p_tu}_{L^\infty}
+\eps\Norm{\Nabla{s}\p_tu}_{L^\infty}
+\eps^2\Norm{\Nabla{s}\p_su}_{L^\infty}  \\
&+\eps\Norm{\Nabla{t}\Nabla{t}v}_{L^\infty}
+\eps^2\Norm{\Nabla{t}\Nabla{s}v}_{L^\infty}
+\eps^3\Norm{\Nabla{s}\Nabla{s}v}_{L^\infty} 
\le c.
\end{align*}
\end{lemma}

\begin{proof}
For every solution $(u,v)$ of~(\ref{eq:floer-V}) define
$$
\tu(s,t) := u(\eps s,t),\qquad \tv(s,t) := \eps v(\eps s,t).
$$
Then 
\begin{equation}\label{eq:tutv}
\p_s\tu - \Nabla{t}\tv = \eps\grad\Vv(\tu),\qquad
\Nabla{s}\tv + \p_t\tu = \frac{\tv}{\eps}.
\end{equation}
By Theorem~\ref{thm:gradient}, Lemma~\ref{le:2D}, 
and~$(V0-V3)$, the function $\tw:=(\tu,\tv)$ and the vector 
field 
$$
\zeta(s,t):=\left(\eps\Vv(u(\eps s,\cdot))(t),v(\eps s,t)\right)
$$ 
along $\tw$ are both uniformly bounded in $W^{3,2}$ 
(under the assumption~(\ref{eq:EA})); here we use the identities
\begin{align*}
\Nabla{t}\p_tu 
&=  \p_su - \grad\Vv(u) - \eps^2\Nabla{t}\Nabla{s}v, \\
\Nabla{t}\Nabla{t}\p_tu 
&=  \Nabla{t}\p_su - \Nabla{t}\grad\Vv(u) 
- \eps^2\Nabla{t}\Nabla{t}\Nabla{s}v, \\
\Nabla{t}\Nabla{t}v 
&=  \Nabla{s}v - \eps^2\Nabla{s}\Nabla{s}v 
- \Nabla{t}\grad\Vv(u), \\
\Nabla{t}\Nabla{t}\Nabla{t}v 
&=  \Nabla{t}\Nabla{s}v 
- \eps^2\Nabla{t}\Nabla{s}\Nabla{s}v 
- \Nabla{t}\Nabla{t}\grad\Vv(u).
\end{align*}
It follows that $\tw$ and $\zeta$ are both
uniformly bounded in $W^{2,p}$ for any $p>2$.  Since 
$$
\p_s\tw + J(\tw)\p_t\tw = \zeta
$$
it follows from~\cite[Proposition~B.4.9]{MS} that $\tu$ and $\tv$ 
are uniformly bounded in $W^{3,p}$ and hence in $C^2$. 
This proves the lemma. 
\end{proof}

\begin{lemma}\label{le:3D}
Fix a constant $c_0>0$ and a perturbation
$\Vv:\Ll M\to\R$ that satisfies~$(V0-V4)$.  
Then there is a constant $C=C(c_0,\Vv)>0$ 
such that the following holds. If $0<\eps\le1$ and 
$(u,v):\R\times S^1\to TM$ is a solution 
of~(\ref{eq:floer-V}) that satisfies~(\ref{eq:EA}),
then 
$$
\int_0^1
\eps^4\Abs{\Nabla{s}\Nabla{s}\p_su(s,t)}^2\,dt
\le CE_{[s-1/2,s+1/2]}^\eps(u,v)
$$
for every $s\in\R$. 
\end{lemma}

\begin{proof}
Define $f_4$ and $g_4$ by 
$$
2f_4 := \Abs{\p_su}^2 + \Abs{\Nabla{s}v}^2
+ \Abs{\Nabla{t}\p_su}^2
+ \Abs{\Nabla{t}\Nabla{t}\p_su}^2,
$$
\begin{align*}
2g_4
&:= \Abs{\Nabla{t}\p_su}^2 
+ \eps^2\Abs{\Nabla{s}\p_su}^2
+ \Abs{\Nabla{t}\Nabla{s}v}^2
+ \eps^2\Abs{\Nabla{s}\Nabla{s}v}^2 \\
&\quad
+ \Abs{\Nabla{t}\Nabla{t}\p_su}^2 
+ \eps^2\Abs{\Nabla{s}\Nabla{t}\p_su}^2 
+ \Abs{\Nabla{t}\Nabla{t}\Nabla{t}\p_su}^2 
+ \eps^2\Abs{\Nabla{s}\Nabla{t}\Nabla{t}\p_su}^2,
\end{align*}
and abbreviate
$
F_4(s) := \int_0^1f_4(s,t)\,dt
$
and
$
G_4(s) := \int_0^1g_4(s,t)\,dt.
$
Then 
\begin{equation}\label{eq:Lf4}
L_\eps f_4 = 2g_4 + U + V + U_t + U_{tt},
\end{equation}
where $U$, $V$, $U_t$ are as in Lemma~\ref{le:gradient-L22}
and
$
U_{tt} := \inner{\Nabla{t}\Nabla{t}\p_su}{\Ll_\eps\Nabla{t}\Nabla{t}\p_su}.
$
We shall prove that there is a constant $\mu>0$ such that
\begin{equation}\label{eq:UV4}
\Abs{U} + \Abs{V} + \Abs{U_t} + \Abs{U_{tt}}
\le \mu f_4 + \frac12\left(g_4 
+ \Norm{\p_su}_{L^2(S^1)}^2 
+ \eps^2\Norm{\Nabla{s}\p_su}_{L^2(S^1)}^2 
\right)
\end{equation}
It follows from~(\ref{eq:Lf4}) and~(\ref{eq:UV4}) that
$$
L_\eps f_4 + \mu f_4 + F_4 \ge g_4 + \frac12(g_4-G_4). 
$$
Integrating this inequality over the interval $0\le t\le 1$ gives
$$
\eps^2F_4'' - F_4' + (\mu+1)F_4\ge 0.
$$
By Theorem~\ref{thm:gradient} and Lemma~\ref{le:2D}, we have
$$
\int_{s-1/4}^{s+1/4}F_4(\sigma)\,d\sigma 
\le cE^\eps_{[s-1/2,s+1/2]}(u,v)
$$
for a suitable constant $c=c(c_0,\Vv)$.  Hence, 
by Lemma~\ref{le:apriori-1eps} with $r=1/5$, 
there is a constant $C=C(c_0,\Vv)$ such that 
$F_4(s)\le CE^\eps_{[s-1/2,s+1/2]}(u,v)$ 
for every $s\in\R$; this gives
$$
\int_0^1
\Abs{\Nabla{t}\Nabla{t}\p_su(s,t)}^2\,dt
\le CE_{[s-1/2,s+1/2]}^\eps(u,v).
$$
Now use~(\ref{eq:Lf-1}) and
$$
\eps^2\Nabla{s}\Nabla{s}\p_su 
= \Ll_\eps\p_su - \Nabla{t}\Nabla{t}\p_su + \Nabla{s}\p_su
$$
to get the required estimate for $\eps^4\Abs{\Nabla{s}\Nabla{s}\p_su}$. 

\smallbreak

For $U$ and $V$ the estimate~(\ref{eq:UV4}) was established in the 
proof of Theorem~\ref{thm:gradient}; for $U_t$ it follows 
from~(\ref{eq:Ut}) via the arguments used in the 
proof of Lemma~\ref{le:gradient-L22}. For $U_{tt}$ we use the identity
\begin{equation}\label{eq:Utt}
\begin{split}
\Ll_\eps\Nabla{t}\Nabla{t}\p_su
&= \Nabla{t}\Ll_\eps\Nabla{t}\p_su
+ \eps^2[\Nabla{s}\Nabla{s},\Nabla{t}]\Nabla{t}\p_su 
- [\Nabla{s},\Nabla{t}]\Nabla{t}\p_su \\
&= \Nabla{t}\Ll_\eps\p_su
+ \eps^2\Nabla{s}\left(R(\p_su,\p_tu)\Nabla{t}\p_su)\right)  \\
&\quad
+ \eps^2R(\p_su,\p_tu)\Nabla{s}\Nabla{t}\p_su
- R(\p_su,\p_tu)\Nabla{t}\p_su \\
&=2\eps^2\Nabla{t}\Nabla{t}\left(R(\p_su,\p_tu)\Nabla{s}v\right)
+ \eps^2\Nabla{t}\Nabla{t}\left(\left(\Nabla{\p_su}R\right)(\p_su,\p_tu)v\right) \\
&\quad
- \Nabla{t}\Nabla{t}\left(R(\p_su,\p_tu)v\right)   \\
&\quad
+\eps^2\Nabla{t}\Nabla{t}\left(R(\Nabla{s}\p_su,\p_tu)v\right)
+\eps^2\Nabla{t}\Nabla{t}\left(R(\p_su,\Nabla{t}\p_su)v\right) \\
&\quad
-\Nabla{t}\Nabla{t}\Nabla{s}\grad\Vv(u)
+\eps^2\Nabla{t}\Nabla{t}\Nabla{s}\Nabla{s}\grad\Vv(u) \\
&\quad
+ \eps^2\Nabla{t}\Nabla{s}\left(R(\p_su,\p_tu)\p_su)\right)  \\
&\quad
+ \eps^2\Nabla{t}\left(R(\p_su,\p_tu)\Nabla{s}\p_su\right)
- \Nabla{t}\left(R(\p_su,\p_tu)\p_su\right) \\
&\quad
+ \eps^2\Nabla{s}\left(R(\p_su,\p_tu)\Nabla{t}\p_su)\right)  \\
&\quad
+ \eps^2R(\p_su,\p_tu)\Nabla{s}\Nabla{t}\p_su
- R(\p_su,\p_tu)\Nabla{t}\p_su.
\end{split}
\end{equation}
Here the last equation follows from~(\ref{eq:Ut}). 
To establish~(\ref{eq:UV4}) we now use the 
pointwise estimates on the first derivatives in 
Theorem~\ref{thm:gradient} and the pointwise estimates 
on the second derivatives in Lemma~\ref{le:2Deps}.  
The term by term analysis shows that all the second,
third, and fourth order factors appear with the appropriate
powers of $\eps$. This proves~(\ref{eq:UV4}) and the lemma. 
\end{proof}

\begin{proof}[Proof of Theorem~\ref{thm:2nd-derivs}.]
For $p=2$ the estimate~(\ref{eq:2D}) follows from 
Theorem~\ref{thm:gradient} and Corollary~\ref{cor:2D}.  
To prove it for $p=\infty$ define $f_5$ and $g_5$ by
$$
2f_5 := \Abs{\p_su}^2 + \Abs{\Nabla{s}v}^2
+ \Abs{\Nabla{t}\p_su}^2
+ \Abs{\Nabla{s}\p_su}^2 
+ \Abs{\Nabla{t}\Nabla{s}v}^2
+ \Abs{\Nabla{s}\Nabla{s}v}^2
$$
and
\begin{align*}
2g_3
&:= \Abs{\Nabla{t}\p_su}^2 
+ \eps^2\Abs{\Nabla{s}\p_su}^2
+ \Abs{\Nabla{t}\Nabla{s}v}^2
+ \eps^2\Abs{\Nabla{s}\Nabla{s}v}^2 \\
&\quad
+ \Abs{\Nabla{t}\Nabla{t}\p_su}^2 
+ \eps^2\Abs{\Nabla{s}\Nabla{t}\p_su}^2 
+ \Abs{\Nabla{t}\Nabla{s}\p_su}^2 
+ \eps^2\Abs{\Nabla{s}\Nabla{s}\p_su}^2  \\
&\quad
+ \Abs{\Nabla{t}\Nabla{t}\Nabla{s}v}^2 
+ \eps^2\Abs{\Nabla{s}\Nabla{t}\Nabla{s}v}^2
+ \Abs{\Nabla{t}\Nabla{s}\Nabla{s}v}^2 
+ \eps^2\Abs{\Nabla{s}\Nabla{s}\Nabla{s}v}^2.
\end{align*}
Then
\begin{equation}\label{eq:Lf5}
L_\eps f_5 = 2g_5 + U + V + U_t + U_s + V_t + V_s
\end{equation}
where $U$, $V$, $U_t$, $U_s$, $V_t$, and $V_s$ 
are as in Lemma~\ref{le:gradient-L22}.
These functions satisfy the estimate
\begin{equation}\label{eq:UV5}
\begin{split}
&\Abs{U}+\Abs{V} + \Abs{U_t} + \Abs{U_s+V_t} + \Abs{V_s}  \\
&\le \mu f_5 + g_5
+ \Norm{\p_su}_{L^2}^2
+ \Norm{\Nabla{s}\p_su}_{L^2}^2
+ \eps^4\Norm{\Nabla{s}\Nabla{s}\p_su}_{L^2}^2
\end{split}
\end{equation}
for all $(s,t)\in\R\times S^1$ and a suitable constant $\mu>0$. 
To see this one argues as in the proof of Lemma~\ref{le:2D} and 
notices that the factor $\eps^2$ in front of $\Abs{V_s}$ is no longer
needed.  (It can now be dropped since, by Corollary~\ref{cor:2D},
the $L^2$-norm of $f_5$ is controlled by the energy.)

By~(\ref{eq:Lf5}) and~(\ref{eq:UV5}), we have
\begin{align*}
L_\eps f_5 +\mu f_5 
&\ge -\Norm{\p_su}_{L^2(S^1)}^2 
-\Norm{\Nabla{s}\p_su}_{L^2(S^1)}^2 
-\eps^4\Norm{\Nabla{s}\Nabla{s}\p_su}_{L^2(S^1)}^2 \\
&\ge -CE^\eps_{[s-1/2,s+1/2]}(u,v)
\end{align*}
for every $s\in\R$ and suitable positive constants
$\mu$ and $C$. Here the last inequality follows from 
Lemmas~\ref{le:2D} and~\ref{le:3D}.  Let $s_0\in\R$ and denote 
$$
a:=\frac{C}{\mu}E^\eps_{[s_0-1,s_0+1]}(u,v).
$$
Then 
$$
L_\eps(f_5+a)+\mu(f_5+a) \ge 0
$$
for $s_0-1/2\le s\le s_0+1/2$. 
Hence we may apply Lemma~\ref{le:apriori-basic-eps}
with $r=1/3$ to the function $w(s,t):=f_6(s_0+s,t_0+t)+a$:
\begin{equation*}
\begin{split}
f_5(s_0,t_0)
&\le54c_2e^{\mu/9}
\int_{s_0-1/9-\eps/3}^{s_0+\eps/3}
\int_0^1\bigl(f_5(s,t)+a\bigr)\,dtds \\
&\le54c_2e^{\mu/9}
\int_{s_0-1/2}^{s_0+1/2}\int_0^1
\bigl(f_5(s) + a\bigr)\,dtds \\
&\le c_3
\Bigl(E_{[s_0-1,s_0+1]}^\eps(u,v)+a\Bigr) \\
&= c_3\left(1 + \frac{C}{\mu}\right)E_{[s_0-1,s_0+1]}^\eps(u,v). 
\end{split}
\end{equation*}
Here the third inequality, with a suitable constant 
$c_3=c_3(c_0,\Vv)>0$, follows from Theorem~\ref{thm:gradient}
and Corollaries~\ref{cor:gradient-L2} and~\ref{cor:2D}.
This proves~(\ref{eq:2D}) for $p=\infty$. 
To prove the result for general $p$ we apply the interpolation inequality
$$
\|\xi\|_{L^p}\le\|\xi\|_{L^\infty}^{1-2/p}\|\xi\|_{L^2}^{2/p}
$$
to the terms on the left hand side of the estimate and use the 
results for $p=2$ and $p=\infty$.  This proves the theorem. 
\end{proof}


\section{Uniform exponential decay}
\label{sec:exp-decay}

\begin{theorem}\label{thm:exp-decay}
Fix a perturbation $\Vv:\Ll M\to\R$ that satisfies~$(V0-V3)$.  
Suppose $\Ss_\Vv$ is Morse and let $a\in\R$ be 
a regular value of $\Ss_\Vv$.
Then there exist positive constants
$\delta,c,\rho$ such that the following holds.
If $x^\pm\in\Pp^a(\Vv)$, $0<\eps\le 1$, $T_0>0$, 
and $(u,v)\in\Mm^\eps(x^-,x^+;\Vv)$ satisfies
\begin{equation}\label{eq:small-energy}
     E_{\R\setminus[-T_0,T_0]}^\eps(u,v)
     <\delta,
\end{equation}
then
$$
     E_{\R\setminus[-T,T]}^\eps(u,v)
     \le ce^{-\rho(T-T_0)}
     E_{\R\setminus[-T_0,T_0]}^\eps(u,v)
$$
for every $T\ge T_0+1$.
\end{theorem}

\begin{corollary}\label{cor:exp-decay}
Fix a perturbation $\Vv:\Ll M\to\R$ that satisfies~$(V0-V3)$.  
Suppose $\Ss_\Vv$ is Morse and let $x^\pm\in\Pp(\Vv)$.
Then there exist positive constants
$\delta,c,\rho$ such that the following holds.
If $0<\eps\le 1$, $T_0>0$, and $(u,v)\in\Mm^\eps(x^-,x^+;\Vv)$ 
satisfies~(\ref{eq:small-energy}) then
\begin{equation}\label{eq:exp-decay}
     \Abs{\p_su(s,t)}^2
     + \Abs{\Nabla{s} v(s,t)}^2
     \le ce^{-\rho|s|}E_{\R\setminus[-T_0,T_0]}^\eps(u,v)
\end{equation}
for every $|s|\ge T_0+2$.
\end{corollary}

\begin{proof}
Theorem~\ref{thm:gradient} and
Theorem~\ref{thm:exp-decay}.
\end{proof}

The proof of Theorem~\ref{thm:exp-decay}
is based on the following two lemmas.

\begin{lemma}[The Hessian]\label{le:hessian}
Fix a perturbation $\Vv:\Ll M\to\R$ that satisfies~$(V0-V2)$.  
Suppose $\Ss_\Vv$ is Morse and fix $a\in\R$.
Then there are positive constants $\delta_0$ and $c$ 
such that the following is true. If $x_0\in\Pp^a(\Vv)$ and
$(x,y)\in\Cinf(S^1,TM)$ satisfy
$$
    x=\exp_{x_0}(\xi_0),\qquad
    y=\Phi(x_0,\xi_0)(\p_tx_0+\eta_0),\qquad
    \Norm{\xi_0}_{W^{1,2}} + \Norm{\eta_0}_\infty
    \le\delta_0,
$$
then
\begin{equation*}
\begin{split}
    &\Norm{\xi}_2 + \Norm{\Nabla{t}\xi}_2
     + \Norm{\eta}_2+\Norm{\Nabla{t}\eta}_2
  \\&
     \le c\left(\Norm{\Nabla{t}\eta+R(\xi,\p_tx)y
     +\Hh_\Vv(x)\xi}_2
     + \Norm{\Nabla{t}\xi-\eta}_2\right)
\end{split}
\end{equation*}
for all $\xi,\eta\in\Om^0(S^1,x^*TM)$.
\end{lemma}

\begin{proof}
The operator
$$
A^\eps(x,y)(\xi,\eta)
:=(-\Nabla{t}\eta-R(\xi,\p_tx)y-\Hh_\Vv(x)\xi,\Nabla{t}\xi-\eta)
$$
on
$
     L^2(S^1,x^*TM\oplus x^*TM)
$
with dense domain
$
    W^{1,2}(S^1,x^*TM\oplus x^*TM)
$
is self-adjoint if $y=\p_tx$.
In the case $(x,y)=(x_0,\p_tx_0)$ it is bijective,
because $\Ss_\Vv$ is Morse.
Hence the result is a consequence
of the open mapping theorem.
Since bijectivity is preserved
under small perturbations
(with respect to the operator norm),
the result for general pairs
$(x,y)$ follows from continuous dependence
of the operator family on the pair $(x,y)$
with respect to the $W^{1,2}$-topology on 
$x$ and the $L^\infty$-topology on $y$.
The set $\Pp^a(\Vv)$ is finite, because
$\Ss_\Vv$ is Morse (see \cite{JOA1}). Hence
we may choose the same constants $\delta_0$ and $c$
for all $x_0\in\Pp^a(\Vv)$.
\end{proof}

\begin{lemma}\label{le:nearby-Crit}
Fix a perturbation $\Vv:\Ll M\to\R$ that satisfies~$(V0)$.  
Suppose $\Ss_\Vv$ is Morse and let $a\in\R$ be
a regular value of $\Ss_\Vv$. 
Then, for every $\delta_0>0$,
there is a constant $\delta_1>0$ 
such that the following is true.
If $(x,y):S^1\to TM$ is a smooth loop such that
$$
     \Aa_\Vv(x,gy)\le a,\qquad
     \Norm{\Nabla{t}y+\grad\Vv(x)}_\infty
     +\Norm{\p_tx-y}_\infty
     <\delta_1,
$$
then there is a periodic orbit $x_0\in\Pp^a(\Vv)$
and a pair of vector fields
$\xi_0,\eta_0\in\Omega^0(S^1,{x_0}^*TM)$
such that
$$
    x=\exp_{x_0}(\xi_0),\qquad
    y=\Phi(x_0,\xi_0)(\p_tx_0+\eta_0),
$$
and
$$
    \Norm{\xi_0}_\infty + \Norm{\Nabla{t}\xi_0}_\infty 
    + \Norm{\eta_0}_\infty + \Norm{\Nabla{t}\eta_0}_\infty
    \le\delta_0.
$$
\end{lemma}

\begin{proof}
First note that 
\begin{equation*}
\begin{split}
\frac12\int_0^1\left|y(t)\right|^2
&= \Aa_\Vv(x,gy) + \Vv(x)
- \int_0^t\inner{y(t)}{\dot x(t)-y(t)}\,dt \\
&\le a + C
+ \int_0^1\left(\frac14\Abs{y(t)}^2
+\Abs{\dot x(t)-y(t)}^2\right)\,dt,
\end{split}
\end{equation*}
where $C$ is the constant in~$(V0)$.
Hence, assuming $\delta_1\le 1$, we have 
$$
      \Norm{y}_2^2
      \le 4\left(a+C+1\right).
$$
Now
\begin{equation*}
\begin{split}
\left|\frac{d}{dt}\left|y\right|^2\right|
&=\bigl|2\inner{y}{\Nabla{t}y+\grad\Vv(x)}
-2\inner{y}{\grad\Vv(x)}\bigr| \\
&\le 2\left(\delta_1+C\right)\Abs{y} 
\le \left(C+1\right)^2 
     + \Abs{y}^2.
\end{split}
\end{equation*}
Integrate this inequality to obtain 
$$
\Abs{y(t_1)}^2-\Abs{y(t_0)}^2
\le \left(C+1\right)^2 + \Norm{y}_2^2
$$
for $t_0,t_1\in[0,1]$.  Integrating again over the interval 
$0\le t_0\le 1$ gives
\begin{equation}\label{eq:yc}
\Norm{y}_\infty
\le \sqrt{\left(C+1\right)^2 + 2\Norm{y}_2^2}
\le c
\end{equation}
where $c^2:=\left(C+1\right)^2+8\left(a+C+1\right)$. 

Now suppose that the assertion is wrong. 
Then there is a $\delta_0>0$ and a sequence of smooth 
loops $(x_\nu,y_\nu):S^1\to TM$ satisfying 
$$
     \Aa_\Vv(x_\nu,gy_\nu)\le a,\qquad
     \lim_{\nu\to\infty}\bigl(
     \Norm{\Nabla{t}y_\nu+\grad\Vv(x_\nu)}_\infty
     +\Norm{\p_tx_\nu-y_\nu}_\infty
     \bigr) = 0,
$$
but not the conclusion of the lemma for the given constant 
$\delta_0$.  By~(\ref{eq:yc}), we have 
$\sup_\nu\Norm{y_\nu}_\infty<\infty$.
Hence $\sup_\nu\Norm{\p_tx_\nu}_\infty<\infty$
and also $\sup_\nu\Norm{\Nabla{t}y_\nu}_\infty<\infty$. 
Hence, by the Arzela--Ascoli theorem, there exists a 
subsequence, still denoted by $(x_\nu,y_\nu)$, that 
converges in the $C^0$-topology.  Our assumptions guarantee
that this subsequence actually converges in the 
$C^1$-topology.  Let $(x_0,y_0):S^1\to TM$ be the limit.
Then $\p_tx_0=y_0$ and $\Nabla{t}y_0+\grad\Vv(x_0)=0$.
Hence $x_0\in\Pp^a(\Vv)$ and $(x_\nu,y_\nu)$ converges
to $(x_0,\p_tx_0)$ in the $C^1$-topology.  
This contradicts our assumption on the sequence 
$(x_\nu,y_\nu)$ and hence proves the lemma. 
\end{proof}

\begin{proof}[Proof of Theorem~\ref{thm:exp-decay}]
To begin with note that 
$
     \Ss_\Vv(x)\ge-C_0
$
for every $x\in\Pp(\Vv)$, where $C_0$ is the constant in~$(V0)$. 
Hence, with 
$
     c_0:=a+C_0, 
$
we have 
$$
x^\pm\in\Pp^a(\Vv)\qquad\implies\qquad
\Ss_\Vv(x^-)\le c_0,\quad\Ss_\Vv(x^-)-\Ss_\Vv(x^+)\le c_0.
$$
Let $C>0$ be the constant of Theorem~\ref{thm:gradient}
with this choice of $c_0$. 
Let $\delta_0$ and $c$ be the constants of Lemma~\ref{le:hessian}
and $\delta_1>0$ the constant of Lemma~\ref{le:nearby-Crit}
associated to $a$ and $\delta_0$.  Then choose $\delta>0$ 
such that $\sqrt{C\delta}\le\delta_1$. Below we will shrink 
the constants $\delta_1$ and $\delta$ further if necessary. 

\smallbreak

In the remainder of the proof we will sometimes 
use the notation $u_s(t):=u(s,t)$ and $v_s(t):=v(s,t)$.
Moreover, $\Norm{\cdot}$ will always denote
the $L^2$ norm on $S^1$ and $\Norm{\cdot}_\infty$
the $L^\infty$ norm on $S^1$. 

Now let $x^\pm\in\Pp^a(\Vv)$, $0<\eps\le1$, and $T_0>0$,
and suppose $(u,v)\in\Mm^\eps(x^-,x^+;\Vv)$ 
satisfies~(\ref{eq:small-energy}).
Then, by Theorem~\ref{thm:gradient}, we have 
\begin{equation}\label{eq:delta1}
     \Norm{\p_su_s}_\infty 
     + \Norm{\Nabla{s}v_s}_\infty
     \le \sqrt{CE_{[s-1,s+1]}^\eps(u,v)}
     \le \sqrt{C\delta}\le\delta_1
\end{equation}
for $|s|\ge T_0+1$. Hence, by Lemma~\ref{le:nearby-Crit},
we know that, for every $s\in\R$ with $|s|\ge T_0+1$,
there is a periodic orbit $x_s\in\Pp^a(\Vv)$ such that
the $C^1$-distance between $(u_s,v_s)$ and $(x_s,\p_tx_s)$
is bounded by $\delta_0$.  Hence we can apply Lemma~\ref{le:hessian}
to the pair $(u_s,v_s)$ and the vector fields 
$(\p_su_s,\Nabla{s}v_s)$ for $|s|\ge T_0+1$. Since
$$
\Nabla{t}\Nabla{s}v+R(\p_su,\p_tu)v+\Hh_\Vv(u)\p_su
= \Nabla{s}\p_su,\qquad
\Nabla{t}\p_su-\Nabla{s}v 
= -\eps^2\Nabla{s}\Nabla{s}v,
$$
we obtain from Lemma~\ref{le:hessian} that
\begin{equation}\label{eq:exp-4}
\begin{split}
    &\Norm{\p_su_s}^2
     +\Norm{\Nabla{t} \p_su_s}^2
     +\Norm{\Nabla{s} v_s}^2
     +\Norm{\Nabla{t}\Nabla{s} v_s}^2 \\
    &\le c\left(\Norm{\Nabla{s}\p_su_s}^2
     +\eps^4\Norm{\Nabla{s}\Nabla{s} v_s}^2\right).
\end{split}
\end{equation}
for $|s|\ge T_0+1$. 

Define the function $F:\R\to[0,\infty)$ by 
$$
     F(s):=\frac12\int_0^1\left(
     \Abs{\p_su(s,t)}^2
     +\eps^2\Abs{\Nabla{s}v(s,t)}^2\right)\,dt.
$$
We shall prove that 
\begin{equation}\label{eq:f''}
F''(s)\ge \frac{1}{c}F(s)
\end{equation}
for $|s|\ge T_0+1$.  The proof of~(\ref{eq:f''}) is based on the identity
\begin{equation}\label{eq:2Df}
\begin{split}
F''(s)    
&=2\Norm{\Nabla{s}\p_su}^2
+2\eps^2\Norm{\Nabla{s}\Nabla{s}v}^2 \\
&\quad
+\inner{\p_su}{\Nabla{s}\Nabla{s}\grad\Vv(u)-\Hh_\Vv(u)\Nabla{s}\p_su} \\
&\quad
     +\left\langle\p_su,
     3R(\p_su,\p_tu)\Nabla{s}v\right\rangle
     +\left\langle\p_su,
     (\Nabla{\p_su}R)(\p_su,\p_tu)v
     \right\rangle \\
    &\quad
     +\left\langle\p_su,R(\p_su,\Nabla{s}v)v\right\rangle
     -\eps^2\left\langle\p_su,R(\p_su,\Nabla{s}\Nabla{s}v)v
     \right\rangle \\
&\quad
+\eps^2\inner{\p_su}{R(\Nabla{s}v,v)\Nabla{s}\p_su}.
\end{split}
\end{equation}
Here all norms and inner products are understood 
in $L^2(S^1,u_s^*TM)$ and we have dropped 
the subscript $s$ for $u_s$ and $v_s$. 
The $L^\infty$ norms of $v$ and $\p_tu=v-\eps^2\Nabla{s}v$
are uniformly bounded, by Theorems~\ref{thm:apriori}
and~\ref{thm:gradient}. Hence there is a constant $c'>0$ 
such that 
\begin{equation*}
\begin{split}
    F''(s)
    &\ge2\Norm{\Nabla{s}\p_su_s}^2
     +2\eps^2\Norm{\Nabla{s}\Nabla{s} v_s}^2 \\
    &\quad
     - c'\Norm{\p_su_s}_\infty\biggl(
     \Norm{\p_su_s}^2
     + \Norm{\p_su_s}\Norm{\Nabla{s}v_s}\biggr) \\
    &\quad
     - c'\eps^2\Norm{\p_su_s}_\infty\biggl(
     \Norm{\p_su_s}\Norm{\Nabla{s}\Nabla{s}v}
     +\Norm{\Nabla{s}v_s}\Norm{\Nabla{s}\p_su_s}
     \biggr) \\
    &\ge\Norm{\Nabla{s}\p_su_s}^2
     +\eps^2\Norm{\Nabla{s}\Nabla{s}v_s}^2.
\end{split}
\end{equation*}
Here the first inequality uses~$(V2)$. 
To understand the last step note that, 
by~(\ref{eq:delta1}), we have 
$\Norm{\p_su_s}_\infty\le\sqrt{C\delta}$
and so the inequality follows from~(\ref{eq:exp-4}), 
provided that $\delta>0$ is sufficiently small.
Now use~(\ref{eq:exp-4}) again to obtain~(\ref{eq:f''}). 

Thus we have proved that $F''(s)\ge\rho^2F(s)$ 
for $|s|\ge T_0+1$, where $\rho:=c^{-1/2}$.  
Since $F(s)$ does not diverge to infinity as $|s|\to\infty$
it follows by standard arguments (see for example~\cite{DOSA,Sa97})
that $F(s)\le e^{-\rho(s-T_0-1)}F(T_0+1)$ for $s\ge T_0+1$
and similarly for $s\le-T_0-1$. 

It remains to prove~(\ref{eq:2Df}). By direct computation, 
$$
F''(s) = \Norm{\Nabla{s}\p_su}^2
+\eps^2\Norm{\Nabla{s}\Nabla{s} v}^2+G(s)+H(s),
$$
where
\begin{equation*}
\begin{split}
G(s)
    &:=
     \left\langle\p_su,
     \Nabla{s}\Nabla{s}\p_su\right\rangle \\
    &=
     \left\langle\p_su,
     \Nabla{s}\Nabla{s}(\Nabla{t}v+\grad\Vv(u))\right\rangle \\
    &=
     \left\langle\p_su,
     [\Nabla{s}\Nabla{s},\Nabla{t}] v
     +\Nabla{s}\Nabla{s}\grad\Vv(u)
     +\Nabla{t}\Nabla{s}\Nabla{s} v
     \right\rangle \\
    &=
     \left\langle\p_su,
     [\Nabla{s}\Nabla{s},\Nabla{t}] v
     +\Nabla{s}\Nabla{s}\grad\Vv(u)
     \right\rangle 
     -\left\langle\Nabla{s}\p_tu,
     \Nabla{s}\Nabla{s} v\right\rangle \\
    &=
     \left\langle\p_su,
     \Nabla{s}[\Nabla{s},\Nabla{t}] v
     +[\Nabla{s},\Nabla{t}]\Nabla{s} v
     +\Nabla{s}\Nabla{s}\grad\Vv(u)
     \right\rangle 
     \\ &\quad
     -\left\langle\Nabla{s}(v-\eps^2\Nabla{s}v),
     \Nabla{s}\Nabla{s}v\right\rangle, \\
H(s)
    &:=
     \eps^2\left\langle\Nabla{s} v,
     \Nabla{s}\Nabla{s}\Nabla{s} v
     \right\rangle \\
    &=
     \left\langle\Nabla{s} v,
     \Nabla{s}\Nabla{s}(v-\p_tu)
     \right\rangle \\
    &=
     \left\langle\Nabla{s} v,
     \Nabla{s}\Nabla{s} v
     -[\Nabla{s},\Nabla{t}]\p_su
     -\Nabla{t}\Nabla{s}\p_su\right\rangle \\
    &=
     \left\langle\Nabla{s} v,
     \Nabla{s}\Nabla{s} v
     -[\Nabla{s},\Nabla{t}]\p_su
     \right\rangle 
     +\left\langle\Nabla{t}\Nabla{s}v,
     \Nabla{s}\p_su\right\rangle \\
    &=\left\langle\Nabla{s}v,
     \Nabla{s}\Nabla{s}v
     -[\Nabla{s},\Nabla{t}]\p_su
     \right\rangle 
     +\left\langle[\Nabla{t},\Nabla{s}]v
     +\Nabla{s}(\p_su-\grad\Vv(u)),
     \Nabla{s}\p_su\right\rangle.
\end{split}
\end{equation*}
Here all inner products are in $L^2(S^1,u_s^*TM)$;
in each formula the fourth step uses integration by parts. 
The sum is 
\begin{equation*}
\begin{split}
G(s) + H(s)
    &=\Norm{\Nabla{s}\p_su}^2
     +\eps^2\Norm{\Nabla{s}\Nabla{s}v}^2  \\
&\quad
     +\inner{\p_su}{\Nabla{s}\Nabla{s}\grad\Vv(u)}
     - \inner{\Nabla{s}\grad\Vv(u)}{\Nabla{s}\p_su} \\
    &\quad
     +\left\langle\p_su,
     3R(\p_su,\p_tu)\Nabla{s}v\right\rangle
     +\left\langle\p_su,
     (\Nabla{\p_su}R)(\p_su,\p_tu)v
     \right\rangle \\
    &\quad
     +\left\langle\p_su,
     R(\p_su,\Nabla{s}\p_tu)v\right\rangle \\
    &\quad
     +\left\langle\p_su,
     R(\Nabla{s}\p_su,\p_tu)v\right\rangle
     -\inner{R(\p_su,\p_tu)v}{\Nabla{s}\p_su}.
\end{split}
\end{equation*}
To obtain~(\ref{eq:2Df}) replace $\Nabla{s}\p_tu$ by 
$\Nabla{s}v-\eps^2\Nabla{s}\Nabla{s}v$.
Moreover, by the first Bianchi identity, 
the last two terms can be expressed in the form
\begin{equation*}
\begin{split}
&\inner{\p_su}{R(\Nabla{s}\p_su,\p_tu)v}
-\inner{R(\p_su,\p_tu)v}{\Nabla{s}\p_su} \\
&= \inner{\p_su}{R(\Nabla{s}\p_su,\p_tu)v} 
+ \inner{\p_su}{R(v,\Nabla{s}\p_su)\p_tu} \\
&= - \inner{\p_su}{R(\p_tu,v)\Nabla{s}\p_su} \\
&= \inner{\p_su}{R(v-\p_tu,v)\Nabla{s}\p_su} \\
&= \eps^2\inner{\p_su}{R(\Nabla{s}v,v)\Nabla{s}\p_su}
\end{split}
\end{equation*}
This proves~(\ref{eq:2Df}) and the theorem.
\end{proof}


\section{Time shift}\label{sec:time-shift}

The next theorem establishes local surjectivity for 
the map $\Tt^\eps$ constructed in Definition~\ref{def:T}.
The idea is to prove that, after
a suitable time shift, the pair $\zeta=(\xi,\eta)$ 
with $u^\eps=\exp_u(\xi)$ and $v^\eps=\Phi(u,\xi)(\p_tu+\eta)$
satisfies the hypothesis $\zeta\in\im\,(\Dd^\eps_u)^*$
of Theorem~\ref{thm:unique}.  The neighbourhood, in which 
the next theorem establishes surjectivity, depends on $\eps$.

\begin{theorem}\label{thm:time-shift}
Assume $\Ss_\Vv$ is Morse--Smale and fix a 
regular value $a\in\R$ of $\Ss_\Vv$. 
Fix two constants $C>0$ and $p>1$.
Then there are positive constants $\delta$, $\eps_0$, 
and $c$ such that $\eps_0\le1$ and the following holds. 
If $x^\pm\in\Pp^a(\Vv)$ is a pair of index difference one, 
$$
u\in\Mm^0(x^-,x^+;\Vv),\qquad
(u^\eps,v^\eps)\in\Mm^\eps(x^-,x^+;\Vv)
$$
with $0<\eps\le\eps_0$, and
$$
      u^\eps = \exp_u(\xi^\eps),
$$
where $\xi^\eps\in\Om^0(\R\times S^1,u^*TM)$
satisfies
\begin{equation}\label{eq:assu-time-1}
     \left\|\xi^\eps\right\|_\infty\le\delta\eps^{1/2},\qquad
     \left\|\xi^\eps\right\|_p\le \delta\eps^{1/2},\qquad
     \left\|\Nabla{t}\xi^\eps\right\|_p\le C,
\end{equation}
then there is a real number $\sigma$ such that
$$
     (u^\eps,v^\eps)=\Tt^\eps
     (u(\sigma+\cdot,\cdot)),
     \qquad
     \Abs{\sigma}< c(\left\|\xi^\eps\right\|_p+\eps^2).
$$
\end{theorem}

\begin{proof}
It suffices to prove the result for a fixed pair 
$x^\pm\in\Pp^a(\Vv)$ of index difference one
and a fixed parabolic cylinder $u\in\Mm^0(x^-,x^+;\Vv)$.
(The assumptions and conclusions of the theorem are 
invariant under simultaneous time shift of 
$u$ and $(u^\eps,v^\eps)$; up to time shift there
are only finitely many index one parabolic cylinders
with $\Ss_\Vv\le a$.) Define 
$$
     c^*:=\Ss_\Vv(x^-)-\Ss_\Vv(x^+) >0.
$$
Let $(u^\eps,v^\eps)\in\Mm^\eps(x^-,x^+;\Vv)$ with $\eps\in(0,1]$.
Denote the time shift of $u$ by 
$$
     u_\sigma(s,t):=u(s+\sigma,t)
$$
for $\sigma\in\R$ and define
$
     \zeta = \zeta(\sigma)
     =\left(\xi,\eta\right)
$
by
\begin{equation}\label{eq:zeta-sigma}
     u^\eps=exp_{u_\sigma}(\xi),\qquad
     v^\eps=\Phi(u_\sigma,\xi)\left(\p_tu_\sigma+\eta\right).
\end{equation}
The pair $(\xi,\eta)$ is well defined whenever
$\sigma\left\|\p_su\right\|_{L^\infty}+\left\|\xi^\eps\right\|_{L^\infty}$
is smaller than the injectivity radius $\rho_M$ of $M$
(i.e. when $\sigma$ and $\delta\eps^{1/2}$ are sufficiently small).
We assume throughout that 
$$
\delta\eps^{1/2}\le\frac{\rho_M}{2}
$$
and choose $\sigma_0>0$ so that 
$\sigma_0\left\|\p_su\right\|_{L^\infty}<\rho_M/2$. 

By Theorem~\ref{thm:par-apriori}
and Theorem~\ref{thm:apriori}, there is a constant
$c_0>0$ such that, for every $\eps\in(0,1]$
and every $(u^\eps,v^\eps)\in\Mm^\eps(x^-,x^+;\Vv)$,
we have
\begin{equation}\label{eq:L-infty-1}
     \Norm{\p_su}_\infty
     +\Norm{\p_tu}_\infty
     +\Norm{v^\eps}_\infty
     \le c_0.
\end{equation}
It follows from~(\ref{eq:zeta-sigma}) and~(\ref{eq:L-infty-1})
that $\left\|\eta(\sigma)\right\|_\infty\le c_0$
for every $\sigma\in[-\sigma_0,\sigma_0]$.
Choose $\delta_0>0$ so small that the assertion of
the Uniqueness Theorem~\ref{thm:unique} holds with 
$C=c_0$ and $\delta=\delta_0$.
We shall prove that for every sufficiently small 
$\eps>0$ there is a $\sigma\in[-\sigma_0,\sigma_0]$
such that
\begin{equation}\label{eq:main}
     \zeta(\sigma)\in\im
     (\Dd^\eps_{u_\sigma})^*,\qquad
     \Norm{\xi(\sigma)}_\infty\le\delta_0\eps^{1/2},\qquad
     \Norm{\eta(\sigma)}_\infty\le c_0.
\end{equation}
It then follows from Theorem~\ref{thm:unique}
that
$
     (u^\eps,v^\eps)=\Tt^\eps(u_{\sigma}).
$
The proof of~(\ref{eq:main}) will take five steps
and uses the following estimate.  Choose $q>1$
such that $1/p+1/q=1$.  Then, 
by parabolic exponential decay
(see Theorem~\ref{thm:par-exp-decay}),
there is a constant $c_1>0$
such that, for  $r=p,q,\infty$,
\begin{equation}\label{eq:Lp-1}
     \Norm{\p_su}_r
     +\Norm{\Nabla{t}\p_su}_r
     +\Norm{\Nabla{s}\p_su}_r
     +\Norm{\Nabla{s}\Nabla{t}\p_su}_r
     \le c_1.
\end{equation}

\vspace{.1cm}
\noindent
{\sc Step 1.}
{\it For $\sigma\in[-\sigma_0,\sigma_0]$
and $\eps>0$ sufficiently small define
$$
     \theta^\eps(\sigma)
     :=-\left\langle Z^\eps_\sigma,\zeta
     \right\rangle_\eps,
$$
where $\zeta=\zeta(\sigma)$ is given
by~(\ref{eq:zeta-sigma}) and
\begin{gather*}
     Z^\eps:=\begin{pmatrix}X^\eps\\Y^\eps
     \end{pmatrix}
     :=\begin{pmatrix}\p_su\\
     \Nabla{t}\p_su\end{pmatrix}
     -\begin{pmatrix}\xi^*\\
     \eta^*\end{pmatrix},
     \\
     \zeta^*
     :=\begin{pmatrix}\xi^*\\\eta^*\end{pmatrix}
     :=(\Dd^\eps_u)^*
     \left(\Dd^\eps_u(\Dd^\eps_u)^*\right)^{-1}
     \Dd^\eps_u
     \begin{pmatrix}\p_su\\\Nabla{t}\p_su
     \end{pmatrix}.
\end{gather*}
Then
$
     \theta^\eps(\sigma)=0
$
if and only if
$
     \zeta\in\im
     (\Dd^\eps_{u_\sigma})^*
$.
}

\vspace{.1cm}
\noindent
For $\eps>0$ sufficiently small, the operator 
$\Dd^\eps_u$ is onto, by Theorem~\ref{thm:inverse},
and, by assumption, it has index one (see Remark~\ref{rmk:eps2}). 
Hence $Z^\eps$ is well defined and belongs to the kernel
of $\Dd^\eps_u$. It remains to prove that $Z^\eps\ne0$
for $\eps>0$ sufficiently small. To see this note that 
$\p_su\ne0$ and so the $(0,2,\eps)$-norm of the pair
$(\p_su,\Nabla{t}\p_su)$ is bounded below by 
a positive constant (the parabolic energy
identity gives $c^*$ as a lower bound).
On the other hand, 
\begin{equation}\label{eq:zeta*}
     \zeta^*
     =(\Dd^\eps_u)^*
     \left(\Dd^\eps_u
     (\Dd^\eps_u)^*\right)^{-1}
     \begin{pmatrix}0\\
     \Nabla{s}\Nabla{t}\p_su\end{pmatrix}.
\end{equation}
Hence, by Theorem~\ref{thm:inverse},
the $(0,2,\eps)$-norm
of $\zeta^*$ converges to zero as $\eps$ tends
to zero.  It follows that $Z^\eps\ne0$ for $\eps>0$
sufficiently small and this proves Step~1. 

\vspace{.1cm}
\noindent
{\sc Step 2.}
{\it 
There are positive constants $\eps_0$ and $c_2$ 
such that
$$
     \Abs{\theta^\eps(0)}\le c_2\left(
     \left\|\xi^\eps\right\|_p+\eps^2\right)
$$
for $0<\eps\le\eps_0$ and every 
$(u^\eps=\exp_u(\xi^\eps),v^\eps)\in\Mm^\eps(x^-,x^+;\Vv)$
satisfying~(\ref{eq:assu-time-1}).}

\vspace{.1cm}
\noindent
We first prove that that there are positive constants 
$\eps_0$ and $c_3$ such that 
\begin{equation}\label{eq:XY}
     \left\|X^\eps\right\|_q+\left\|Y^\eps\right\|_q\le c_3
\end{equation}
for $0<\eps\le\eps_0$. For the summands $\p_su$ of $X^\eps$
and $\Nabla{s}\p_tu$ of $Y^\eps$ this follows from~(\ref{eq:Lp-1})
with $r=q$. Moreover, by~(\ref{eq:zeta*})
and Theorem~\ref{thm:inverse}, we have 
\begin{equation*}\label{eq:xieta*}
\begin{split}
     \left\|\xi^*\right\|_q
     +\eps^{1/2}\left\|\eta^*\right\|_q
    &\le 
     c_4\left(
     \eps\Norm{(0,\Nabla{s}\Nabla{t}\p_su)}_{0,q,\eps}
     +\Norm{\pi_\eps(0,\Nabla{s}\Nabla{t}\p_su)}_q
     \right) \\
    &\le 
     c_4\eps^{3/2}(\eps^{1/2}+\kappa_q)
     \Norm{\Nabla{s}\Nabla{t}\p_su}_q.
\end{split}
\end{equation*}
The last step uses Lemma~\ref{le:eat-eps}
with constant $\kappa_q>1$.
This proves~(\ref{eq:XY}). 
It follows from~(\ref{eq:XY}) that
\begin{equation}\label{eq:theta0}
     \left|\theta^\eps(0)\right|
     \le c_3\left(
         \left\|\xi^\eps\right\|_p + \eps^2\left\|\eta^\eps\right\|_p
         \right),
\end{equation}
where $\eta^\eps\in\Om^0(\R\times S^1,u^*TM)$ is defined by 
$$
     v^\eps=:\Phi(u,\xi^\eps)(\p_tu^\eps+\eta^\eps).
$$
Define the linear maps
$
     E_i(x,\xi):T_xM\to T_{exp_x(\xi)}M
$
by the formula
\begin{equation}\label{eq:E-i}
     \frac{d}{d\tau} exp_x(\xi)
     =:E_1(x,\xi)\p_\tau x
     +E_2(x,\xi)\Nabla{\tau}\xi
\end{equation}
for every smooth path $x:\R\to M$ and
every vector field
$\xi\in\Omega^0(\R,x^*TM)$ along $x$.
Abbreviate $\Phi:=\Phi(u,\xi^\eps)$ and 
$E_i:=E_i(u,\xi^\eps)$ for $i=1,2$. 
Then
\begin{equation*}
\begin{split}
     \eta^\eps
    &=
     \Phi^{-1}v^\eps-\p_tu \\
    &=
     \Phi^{-1}(v^\eps-\p_tu^\eps)
     +\Phi^{-1}(E_1\p_tu+E_2\Nabla{t}\xi_\eps)
     -\p_tu \\
    &=
     \eps^2\Phi^{-1}\Nabla{s}v^\eps
     +\Phi^{-1}E_2\Nabla{t}\xi^\eps
     +(\Phi^{-1}E_1-\1)\p_tu.
\end{split}
\end{equation*}
By Corollary~\ref{cor:exp-decay}, there is a constant
$c_5$ such that 
$\eps\left\|\Nabla{s}v^\eps\right\|_p\le c_5$.
Moreover, there is a constant $c_6>0$ such that 
$\left\|\Phi^{-1}E_1-\1\right\|_p\le c_6\left\|\xi^\eps\right\|_p$. 
Hence there is another constant $c_7>0$ such that 
$$
     \left\|\eta^\eps\right\|_p
     \le c_7\left(\eps
     +\left\|\Nabla{t}\xi^\eps\right\|_p+\left\|\xi^\eps\right\|_p\right)
     \le c_7\left(\left\|\xi^\eps\right\|_p + C + 1\right).
$$
Combining this with~(\ref{eq:theta0}) proves Step~2. 

\vspace{.1cm}
\noindent
{\sc Step~3.}
{\it 
There is a constant $c_8>0$ such that
$$
     \Norm{\xi(\sigma)}_\infty\le \delta\eps^{1/2}+c_8\Abs{\sigma},\qquad
     \Norm{\eta(\sigma)}_\infty\le c_0,
$$
$$
     \Norm{\Nabla{s}\xi}_p\le c_8,\quad
     \left\|\Nabla{\sigma}\xi+\p_su_\sigma\right\|_p
     \le c_8\left(\left|\sigma\right|+\delta\eps^{1/2}\right),\quad
     \Norm{\xi(\sigma)}_p\le \delta\eps^{1/2}+c_8|\sigma|
$$
for $0<\eps\le\eps_0$ and $|\sigma|\le\sigma_0$.}

\smallbreak

\vspace{.1cm}
\noindent
For every $\sigma\in\R$, we have
$$
     d\left(u(s+\sigma,t),u(s,t)\right)
     \le L(\gamma)\le \Abs{\sigma}\Norm{\p_su}_\infty,
$$
where $\gamma(r):=u(s+r\sigma,t)$, $0\le r\le1$. 
Moreover, by~(\ref{eq:assu-time-1}),
$d(u(s,t),u^\eps(s,t))\le\delta\eps^{1/2}$. Hence the first
estimate of Step~3 follows from the triangle inequality.
The second estimate follows from the identity
$$
     \eta(\sigma)=\Phi(u_\sigma,\xi(\sigma))^{-1}v^\eps-\p_tu_\sigma
$$
and~(\ref{eq:L-infty-1}). To prove the next two estimates 
we differentiate the identity
$$
    \exp_{u_\sigma}(\xi(\sigma))=u^\eps
$$ 
with respect to $\sigma$ and $s$
to obtain
$$
    E_1(u_\sigma,\xi)\p_su_\sigma 
    + E_2(u_\sigma,\xi)\Nabla{\sigma}\xi = 0,\qquad
    E_1(u_\sigma,\xi)\p_su_\sigma 
    + E_2(u_\sigma,\xi)\Nabla{s}\xi = \p_su^\eps.
$$
By the energy identities the $L^2$ norms of $\p_su$ and $\p_su^\eps$
are uniformly bounded and hence, so is the $L^2$ norm of $\Nabla{s}\xi$.
Moreover,
$$
    \left\|\Nabla{\sigma}\xi+\p_su_\sigma\right\|_p
    = \left\|\left(E_2^{-1}E_1-\1\right)\p_su\right\|_p
    \le c_9\left\|\xi(\sigma)\right\|_\infty
    \le c_{10}\left(\left|\sigma\right|+\delta\eps^{1/2}\right).
$$
Hence the $L^p$ norm of $\Nabla{\sigma}\xi$ is uniformly bounded.
Now differentiate the function $\sigma\mapsto\Norm{\xi(\sigma)}_p$
to obtain the inequality
$\Norm{\xi(\sigma)}_p\le\Norm{\xi(0)}_p+c_{11}|\sigma|$.
Then the last inequality in Step~3 follows from~(\ref{eq:assu-time-1}).

\vspace{.1cm}
\noindent
{\sc Step 4.}
{\it Shrinking $\sigma_0$ and $\eps_0$, if necessary, we have
$$
     \frac{d}{d\sigma}\theta^\eps(\sigma)
     \ge \frac{c^*}{2}
$$
for $0<\eps\le\eps_0$ and $|\sigma|\le\sigma_0$.}

\vspace{.1cm}
\noindent
We will investigate the two terms in the sum
\begin{equation}\label{eq:term-main}
     \frac{d}{d\sigma}\theta^\eps(\sigma)
     =- \frac{d}{d\sigma}\left\langle
     X^\eps_\sigma,\xi(\sigma)\right\rangle
     -\eps^2 \frac{d}{d\sigma}\left\langle
     Y^\eps_\sigma,\eta(\sigma)\right\rangle
\end{equation}
separately. The key term is
$\inner{X_\sigma^\eps}{\Nabla{\sigma}\xi}$.
We have seen that $X_\sigma^\eps$ is
$L^q$-close to $\p_su_\sigma$
and $\Nabla{\sigma}\xi$ is $L^p$-close
to $-\p_su_\sigma$.  
We shall prove that all the other terms are small 
and hence $\p_\sigma\theta^\eps$ is approximately
equal to $\left\|\p_su\right\|_2^2$. More precisely,
for the first term in~(\ref{eq:term-main}) we obtain
\begin{equation*}
\begin{split}
     -\frac{d}{d\sigma}\left\langle
     X^\eps_\sigma,\xi\right\rangle 
    &=
     -\inner{X^\eps_\sigma}{\Nabla{\sigma}\xi}
     -\inner{\Nabla{s}X^\eps_\sigma}{\xi} \\
    &=
     \left\|\p_su\right\|_2^2
     -\inner{X^\eps_\sigma}{\p_su_\sigma+\Nabla{\sigma}\xi}
     -\inner{\xi^*}{\p_su_\sigma} \\
    &
     -\inner{\Nabla{s}\p_su_\sigma}{\xi}-\inner{\xi^*_\sigma}{\Nabla{s}\xi} \\
    &\ge
     \left\|\p_su\right\|_2^2
     -c_{12}\left(\Norm{\p_su_\sigma+\Nabla{\sigma}\xi}_p
     +\Norm{\xi^*}_q+\Norm{\xi}_p\right) \\
    &\ge
     \left\|\p_su\right\|_2^2
     -c_{13}\left(|\sigma|+\delta\eps^{1/2}+\eps^{3/2}\right).
\end{split}
\end{equation*}
Here the second step follows from integration by parts.
The third step uses the inequalities 
$\Norm{X^\eps}_q\le c$ (see~(\ref{eq:XY})), 
$\Norm{\p_su}_p+\Norm{\Nabla{s}\p_su}_q\le c$
(see~(\ref{eq:Lp-1})),
and $\Norm{\Nabla{s}\xi}_p\le c$ (see Step~3). 
The last step uses Step~3 and~(\ref{eq:xieta*}).

To estimate the second term in~(\ref{eq:term-main}) we 
differentiate the identity
$$
     \Phi(u_\sigma,\xi(\sigma))
     (\p_tu_\sigma+\eta(\sigma))=v^\eps
$$
with respect to $\sigma$ to obtain
$$
     \Norm{\Nabla{s}\p_tu_\sigma+\Nabla{\sigma}\eta}_p
     \le c_{14}\left(
     \Norm{\p_su}_p+\Norm{\Nabla{\sigma}\xi}_p
     \right)
     \le c_{15}.
$$
In the first inequality we have used the
fact that the $L^\infty$ norms
of $\eta(\sigma)$ and $\p_tu_\sigma$
are uniformly bounded.
In the second inequality we have used Step~3.
Combining this estimate with~(\ref{eq:Lp-1}) we 
find that the $L^p$ norm of $\Nabla{\sigma}\eta$
is uniformly bounded.
Differentiating the same identity
with respect to $s$ we obtain 
$$
     \Norm{\Nabla{s}\p_tu_\sigma+\Nabla{s}\eta}_p
     \le c_{16}\left(
     \Norm{\Nabla{s}v^\eps}_p+\Norm{\p_su}_p
     +\Norm{\Nabla{s}\xi}_p
     \right)
     \le c_{17}\eps^{-1}.
$$
Here the last inequality follows from Step~3 and 
Corollary~\ref{cor:exp-decay}. Using~(\ref{eq:Lp-1})
again, we obtain that the $L^p$ norm of 
$\eps\Nabla{s}\eta$
is uniformly bounded. Now
\begin{equation*}
\begin{split}
     \eps^2\frac{d}{d\sigma}\left\langle
     Y^\eps_\sigma,\eta\right\rangle 
    &=
     \eps^2\inner{\Nabla{s}Y^\eps_\sigma}{\eta}  
     +\eps^2\inner{Y^\eps_\sigma}{\Nabla{\sigma}\eta} \\
    &=
     -\eps^2\inner{Y^\eps_\sigma}{\Nabla{s}\eta} 
     +\eps^2\inner{Y^\eps_\sigma}{\Nabla{\sigma}\eta} \\
    &\le
     c_{18}\eps.
\end{split}
\end{equation*}
In the last estimate we have used~(\ref{eq:XY})
and the uniform estimates on the $L^p$ norms
of $\Nabla{\sigma}\eta$ and $\eps\Nabla{s}\eta$.
Putting things together we obtain
$$
     \frac{d}{d\sigma}
     \theta^\eps(\sigma)
     \ge \left\|\p_su\right\|_2^2 
     - c_{19}\left(|\sigma|+\eps^{1/2}\right).
$$
Since $\left\|\p_su\right\|_2^2=c^*$, the assertion of Step~4
holds whenever $0<\eps\le\eps_0$, $|\sigma|\le\sigma_0$,
and $c_{19}(\sigma_0+\eps_0^{1/2})\le c^*/2$.

\vspace{.1cm}
\noindent
{\sc Step 5.}
{\it We prove Theorem~\ref{thm:time-shift}.
}

\vspace{.1cm}
\noindent
Suppose the pair $(u^\eps,v^\eps)$ satisfies the 
requirements of the theorem with $\eps$ and $\delta$
sufficiently small. Then, by Steps~2 and~4, there is 
a $\sigma\in[-\sigma_0,\sigma_0]$ such that
$$
     \theta^\eps(\sigma)=0,\qquad
     |\sigma|\le c_{20}(\Norm{\xi^\eps}_p+\eps^2),\qquad
     c_{20}:=\frac{2c_2}{c^*}.
$$
Let $\xi:=\xi(\sigma)$ and $\eta:=\eta(\sigma)$.
Then, by Step~3,
$$
     \Norm{\xi}_\infty
     \le(\delta+c_8c_{20}(\delta+\eps^{3/2}))\eps^{1/2},\qquad
     \Norm{\eta}_\infty\le c_0.
$$
If $\delta+c_8c_{20}(\delta+\eps^{3/2})\le\delta_0$ then,
by Step~1, $\zeta:=(\xi,\eta)\in\im\,(\Dd^\eps_{u_\sigma})^*$.
Hence, by Theorem~\ref{thm:unique},
$(u^\eps,v^\eps)=\Tt^\eps(u_\sigma)$.
\end{proof}


\section{Surjectivity}\label{sec:onto}

\begin{theorem}\label{thm:onto}
Assume $\Ss_\Vv$ is Morse--Smale and fix 
a constant $a\in\R$.  Then there is a constant $\eps_0>0$ 
such that, for every $\eps\in(0,\eps_0)$ and
every pair $x^\pm\in\Pp^a(\Vv)$ of index difference one,
the map $\Tt^\eps:\Mm^0(x^-,x^+;\Vv)\to\Mm^\eps(x^-,x^+;\Vv)$,
constructed in Definition~\ref{def:T}, is bijective. 
\end{theorem}

\begin{lemma}\label{le:compact}
Assume $\Ss_\Vv$ is Morse.  Let 
$x^\pm\in\Pp(\Vv)$ and $u_i\in\Mm^{\eps_i}(x^-,x^+;\Vv)$
where $\eps_i$ is a sequence of positive real numbers 
converging to zero.  Then there is a pair $x_0,x_1\in\Pp(\Vv)$,
a parabolic cylinder $u\in\Mm^0(x_0,x_1;\Vv)$,
and a subsequence, still denoted by $(u_i,v_i)$,
such that the following holds.
\begin{enumerate}
\item[\rm\bfseries(i)]
$(u_i,v_i)$ converges to $(u,v)$ strongly in $C^1$ and weakly 
in $W^{2,p}$ on every compact subset
of $\R\times S^1$ and for every $p>1$.
Moreover, $v_i-\p_tu_i$ converges to zero in the 
$C^1$ norm on every compact subset of $\R\times S^1$.
\item[\rm\bfseries(ii)]
For all $s\in\R$ and $T>0$,
\begin{align*}
\Ss_\Vv(u(s,\cdot))
&=\lim_{i\to\infty}\Aa_\Vv(u_i(s,\cdot),v_i(s,\cdot)), \\
E_{[-T,T]}(u)
&= \lim_{i\to\infty}E^\eps_{[-T,T]}(u_i,v_i).
\end{align*}
\end{enumerate}
\end{lemma}

\begin{proof}
By Theorems~\ref{thm:apriori}, \ref{thm:gradient},
and~\ref{thm:2nd-derivs} there is a constant $c>0$ such that
\begin{equation}\label{eq:i1}
\Norm{v_i}_\infty + \Norm{\p_tu_i}_\infty
+ \Norm{\p_su_i}_\infty
+ \Norm{\Nabla{t}v_i}_\infty
+ \Norm{\Nabla{s}v_i}_\infty
\le c,
\end{equation}
\begin{equation}\label{eq:i2}
\Norm{\Nabla{s}\p_tu_i}_p
+ \Norm{\Nabla{s}\p_su_i}_p
+ \Norm{\Nabla{t}\Nabla{s}v_i}_p
+ \Norm{\Nabla{s}\Nabla{s}v_i}_p
\le c,
\end{equation}
\begin{equation}\label{eq:i3}
\Norm{\Nabla{t}\p_tu_i}_\infty
+ \Norm{\Nabla{t}\Nabla{t}v_i}_\infty
\le c
\end{equation}
for every $i\in\N$ and every $p\in[2,\infty]$.
In~(\ref{eq:i1}) the estimate for $\Nabla{t}v_i$ 
follows from the one for $\p_su_i$ and the identity 
$\Nabla{t}v_i=\p_su_i-\grad\Vv(u_i)$. 
The estimate for $\p_tu_i$ follows from the ones 
for $v_i$ and $\Nabla{s}v_i$ and the identity 
$\p_tu_i=v_i-\eps_i^2\Nabla{s}v_i$. In~(\ref{eq:i3})
the estimate for $\Nabla{t}\p_tu_i$ follows from 
the ones for $\Nabla{t}v_i$ and $\Nabla{t}\Nabla{s}v_i$ 
and the identity $\Nabla{t}\p_tu_i=\Nabla{t}v_i
-\eps_i^2\Nabla{t}\Nabla{s}v_i$.
The estimate for $\Nabla{t}\Nabla{t}v_i$ follows from the 
ones for $\Nabla{t}\p_su_i$ and $\p_tu_i$ and the identity 
$\Nabla{t}\Nabla{t}v_i=\Nabla{t}\p_su_i-\Nabla{t}\grad\Vv(u_i)$.

By~(\ref{eq:i1}), (\ref{eq:i2}), and~(\ref{eq:i3}) the 
sequence $(u_i,v_i)$ is bounded in 
$C^2$ and hence in $W^{2,p}([-T,T]\times S^1)$
for every $T>0$ and every $p>1$. Hence, by the Arzela-Ascoli
theorem and the Banach--Alaoglu theorem, a suitable 
subsequence, still denoted by $(u_i,v_i)$, 
converges strongly in $C^1$ and weakly in $W^{2,p}$ 
on every compact subset of $\R\times S^1$ to some 
$C^2$-funtion $(u,v):\R\times S^1\to TM$. 
By~(\ref{eq:i1}) and~(\ref{eq:i2}), the sequence 
$$
v_i-\p_tu_i=\eps_i^2\Nabla{s}v_i
$$
converges to zero in the $C^1$ norm. Hence $v=\p_tu$. 
Moreover, the sequence 
$$
     \p_su_i-\Nabla{t}\p_tu_i-\grad\Vv(u_i)
     =\eps_i^2\Nabla{t}\Nabla{s}v_i
$$
converges to zero in the sup-norm, by~(\ref{eq:i2}),
so the limit $u:\R\times S^1\to M$ satisfies 
the parabolic equation~(\ref{eq:heat-V}). By the 
parabolic regularity theorem~\ref{thm:par-regularity}, 
$u$ is smooth and so is $v=\p_tu$.  
This proves~(i). 

To prove~(ii) note that
\begin{equation*}
\begin{split}
E_{[-T,T]}(u)
&=\int_{-T}^{T}\int_0^1\left|\p_su\right|^2\,dsdt \\
&=\lim_{i\to\infty}\int_{-T}^{T}\int_0^1
  \left|\p_su_i\right|^2\,dsdt \\
&=\lim_{i\to\infty}\int_{-T}^{T}\int_0^1
  \bigl(\left|\p_su_i\right|^2
  +\eps_i^2\left|\Nabla{s}v_i\right|^2\bigr)\,dsdt \\
&=\lim_{i\to\infty}E_{[-T,T]}(u_i,v_i)
\end{split}
\end{equation*}
for every $T$; here the third identity followws 
from~(\ref{eq:i1}). Hence the limit $u$ has finite energy 
and so belongs to the moduli space $\Mm^0(x_0,x_1;\Vv)$ for 
some pair $x_0,x_1\in\Pp(\Vv)$. To prove convergence 
of the symplectic action at time $s$ note that
$$
\Vv(u(s,\cdot)) = \lim_{i\to\infty}\Vv(u_i(s,\cdot)),
$$
because $\Vv$ is continuous with respect to the 
$C^0$ topology on $\Ll M$.  Moreover
\begin{align*}
\Ss_0(u(s,\cdot))
&= \int_0^1\Abs{\p_tu(s,t)}^2\,dt  \\
&= \lim_{i\to\infty}
   \int_0^1\left(
   \inner{\p_tu_i(s,t)}{v_i(s,t)}
   - \frac12\Abs{v_i(s,t)}\right)\,dt \\
&= \lim_{i\to\infty}\Aa_0(u_i(s,\cdot),v_i(s,\cdot)).
\end{align*}
Here the second equality follows from the fact that
$\p_tu_i(s,\cdot)$ and $v_i(s,\cdot)$ both converge
to $\p_tu(s,\cdot)$ in the sup-norm.
This proves the lemma. 
\end{proof}

\begin{lemma}\label{le:catenation}
Assume $\Ss_\Vv$ is Morse.  Let 
$x^\pm\in\Pp(\Vv)$ and $u_i\in\Mm^{\eps_i}(x^-,x^+;\Vv)$
where $\eps_i$ is a sequence of positive real numbers 
converging to zero.  Then there exist periodic orbits 
$
     x^-=x^0,x^1,\dots,x^\ell=x^+\in\Pp(\Vv),
$
parabolic cylinders $u^k\in\Mm^0(x^{k-1},x^k;\Vv)$
for $k\in\{1,\dots,\ell\}$,
a subsequence, still denoted by $(u_i,v_i)$, 
and sequences $s^k_i\in\R$, $k\in\{1,\dots,\ell\}$,
such that the following holds.
\begin{enumerate}
\item[\rm\bfseries(i)]
For every $k\in\{1,\dots,\ell\}$ the sequence
$(s,t)\mapsto(u_i(s^k_i+s,t),v_i(s^k_i+s,t))$
converges to $(u^k,\p_tu^k)$
as in Lemma~\ref{le:compact}.
\item[\rm\bfseries(ii)]
$s^k_i-s^{k-1}_i$ diverges to infinity
for $k=2,\dots,\ell$
and $\p_su^k\not\equiv0$
for $k=1,\dots,\ell$.
\item[\rm\bfseries(iii)]
For every $k\in\{0,\dots,\ell\}$ and every 
$\rho>0$ there is a constant 
$T>0$ such that, for every $i$ and 
every $(s,t)\in\R\times S^1$,
$$
      s^k_i+T\le s\le s^{k+1}_i-T\qquad\IMP\qquad
      d(u_i(s,t),x^k(t))<\rho.
$$
(Here we abbreviate $s^0_i:=-\infty$ 
and $s^{\ell+1}_i:=\infty$.)
\end{enumerate}
\end{lemma}

\begin{proof}
Denote $a:=\Ss_\Vv(x^-)$ and choose $\rho>0$ so small 
that $d(x(t),x'(t))>2\rho$ for every $t\in\R$ and any 
two distinct periodic orbits $x,x'\in\Pp^a(\Vv)$.  
Choose $s^1_i$ such that 
\begin{equation}\label{eq:cat1}
       \sup_{s\le s^1_i}\sup_td(x^-(t),u_i(s,t))
       \le\rho,\qquad
       \sup_td(x^-(t),u_i(s^1_i,t))=\rho.
\end{equation}
Passing to a subsequence we may assume, 
by Lemma~\ref{le:compact}, that the sequence
$(u_i(s^1_i+\cdot,\cdot),v_i(s^1_i+\cdot,\cdot))$
converges in the required sense to a parabolic cylinder
$u^1\in\Mm^0(x^0,x^1;\Vv)$, where $x^0,x^1\in\Pp^a(\Vv)$.
By~(\ref{eq:cat1}), we have $x^0=x^-$ and $x^1\ne x^0$.
Hence $\p_su^1\not\equiv0$ and so 
$\Ss_\Vv(x^1)<\Ss_\Vv(x^0)$. 
If $x^1=x^+$ the lemma is proved.
If $x^1\ne x^+$ choose $T>0$ such 
that $d(u^1(s,t),x^1(t))<\rho$ 
for every~$t$ and every $s\ge T$.  
Passing to a subsequence,
we may assume that $d(u_i(s^1_i+T,t),x^1(t))<\rho$ 
for every $t$.  Since $x^1\ne x^+$
there exists a sequence $s^2_i>s^1_i+T$ such that 
$$
       \sup_{s_i^1+T\le s\le s^2_i}
       \sup_td(x^1(t),u_i(s,t))
       \le\rho,\qquad
       \sup_td(x^1(t),u_i(s^2_i,t))=\rho.
$$
The difference $s^2_i-s^1_i$ diverges to infinity
and, by Lemma~\ref{le:compact}, there is a further 
subsequence such that 
$(u_i(s^2_i+\cdot,\cdot),v_i(s^2_i+\cdot,\cdot))$
converges to a parabolic cylinder 
$u^2\in\Mm^0(x^1,x^2;\Vv)$,
where $\Ss_\Vv(x^2)<\Ss_\Vv(x^1)$. Continue by induction.  
The induction can only terminate if $x^\ell=x^+$.  
It must terminate because $\Pp^a(\Vv)$ is a finite set. 
This proves the lemma. 
\end{proof}

\begin{proof}[Proof of Theorem~\ref{thm:onto}.]
By Theorem~\ref{thm:unique} the map $\Tt^\eps$ is injective 
for $\eps>0$ sufficiently small. We will prove surjectivity 
by contradiction.

Assume the result is false.
Then there exist periodic orbits $x^\pm\in\Pp^a(\Vv)$ of 
Morse index difference one and sequences $\eps_i>0$ and 
$(u_i,v_i)\in\Mm^{\eps_i}(x^-,x^+;\Vv)$ such that
\begin{equation} \label{eq:surject-statement}
     \lim_{i\to\infty}\eps_i=0,\qquad
     (u_i,v_i)\notin
     \Tt^{\eps_i}(\Mm^0(x^-,x^+;\Vv)).
\end{equation}
Applying a time shift, if necessary, we assume 
without loss of generality that 
\begin{equation}\label{eq:normalization-action}
     \Aa_\Vv\bigl(u_i(0,\cdot),v_i(0,\cdot)\bigr)
     =\frac12\bigl(\Ss_\Vv(x^-)+\Ss_\Vv(x^+)\bigr).
\end{equation}
Fix a constant $p>2$.
We shall prove in two steps that, 
after passing to a subsequence
if necessary, there is a sequence 
$u^0_i\in\Mm^0(x^-,x^+;\Vv)$
and a constant $C>0$ such that 
$$
    u_i=exp_{u^0_i}(\xi_i),
$$
where the sequence
$\xi_i\in\Om^0(\R\times S^1,(u^0_i)^*TM)$ satisfies
\begin{equation}\label{eq:estimate-zeta-i-0-new}
\begin{gathered}
    \lim_{i\to\infty}\eps_i^{-1/2}\bigl(\|\xi_i\|_\infty
    +\|\xi_i\|_p\bigr) = 0,\qquad
    \|\Nabla{t}\xi_i\|_p\le C.
\end{gathered}
\end{equation}
Hence it follows from 
Theorem~\ref{thm:time-shift} that, 
for $i$ sufficiently large, there is a real number 
$\sigma_i$ such that
$(u_i,v_i)=\Tt^{\eps_i}(u^0_i(\sigma_i+\cdot,\cdot))$.
This contradicts~(\ref{eq:surject-statement})
and hence proves Theorem~\ref{thm:onto}.

\medskip
\noindent
{\sc Step 1.}
{\it 
For every $\delta>0$ there is a constant
$T_0>0$ such that
\begin{equation}\label{SU-step2}
     E_{\R\setminus[-T_0,T_0]}^{\eps_i}(u_i,v_i)
     <\delta
\end{equation}
for every $i\in\N$.}

\medskip
\noindent
Assume, by contradiction, that the statement is false.
Then there is a constant $\delta>0$,
a sequence of positive real numbers $T_i\to\infty$, 
and a subsequence, still denoted by 
$(\eps_i,u_i,v_i)$, such that, for every $i\in\N$,
\begin{equation}\label{eq:xpmdelta}
     E_{[-T_i,T_i]}^{\eps_i}(u_i,v_i)
     \le\Ss_\Vv(x^-)-\Ss_\Vv(x^+)-\delta.
\end{equation}
Choose a further subsequence, still denoted by $(u_i,v_i)$, 
that converges as in Lemma~\ref{le:catenation} to a
finite collection of parabolic cylinders
$u^k\in\Mm^0(x^{k-1},x^k;\Vv)$, $k=1,\dots,\ell$,
with $x^-=x^0,x^1,\dots,x^{\ell-1},x^\ell=x^+\in\Pp(\Vv)$.
We claim that $\ell\ge 2$. Otherwise, $u_i(s_i+\cdot,\cdot)$
converges to $u:=u^1\in\Mm^0(x^-,x^+;\Vv)$ as in Lemma~\ref{le:compact}
(for some sequence $s_i\in\R$).  
By~(\ref{eq:normalization-action}) and Lemma~\ref{le:compact}~(iv), 
the sequence $s_i$ must be bounded.  By~(\ref{eq:xpmdelta}),
this implies that 
$$
   E_{[-T,T]}(u)
   = \int_{-T}^T\int_0^1\left|\p_su\right|^2\,dtds
   \le\Ss_\Vv(x^-)-\Ss_\Vv(x^+)-\delta
$$ 
for every $T>0$.
This contradicts the fact that $u$ connects $x^-$
with $x^+$. Thus we have proved that $\ell\ge 2$
as claimed. Since $\Ss_\Vv$ is Morse--Smale 
it follows that the Morse index difference of $x^-$
and $x^+$ is at least two.  This contradicts our
assumption and proves Step~1. 

\medskip
\noindent
{\sc Step 2.}
{\it For $i$ sufficiently large 
there is a parabolic cylinder $u_i^0\in\Mm^0(x^-,x^+;\Vv)$ 
and a vector field 
$\xi_i\in\Om^0(\R\times S^1,(u^0_i)^*TM)$
such that
$
      u_i=exp_{u_i^0}(\xi_i)
$
and $\xi_i$ satisfies~(\ref{eq:estimate-zeta-i-0-new}).}

\medskip
\noindent
Let $\delta$, $c$ and $\rho$ denote the constants
in Theorem~\ref{thm:exp-decay} and choose $T_0>0$,
according to Step~1, such that~(\ref{SU-step2}) 
holds with this constant $\delta$.
Then, by Corollary~\ref{cor:exp-decay},
\begin{equation}\label{eq:SU-1}
     \Abs{\p_su_i(s,t)}^2
     +\Abs{\Nabla{s}v_i(s,t)}^2
     \le c_3e^{-\rho|s|}
     E_{\R\setminus[-T_0,T_0]}^{\eps_i}(u_i,v_i)
\end{equation}
for $|s|\ge T_0+2$ and a suitable constant $c_3>0$. 
By Theorem~\ref{thm:apriori} and Theorem~\ref{thm:gradient},
there is a constant $c_4>0$ such that
\begin{equation}\label{eq:a-priori-proof-new}
     \Norm{v_i}_\infty
     +\Norm{\p_su_i}_\infty
     +\Norm{\p_tu_i}_\infty
     +\Norm{\Nabla{s} v_i}_\infty
     \le c_4
\end{equation}
for every $i$.  Here we have also used the identity
$\p_tu_i=v_i-\eps_i^2\Nabla{s}v_i$.  It follows 
from~(\ref{eq:SU-1}) and~(\ref{eq:a-priori-proof-new})
that there is a constant $c_5\ge c_4$ such that 
$$
     \Abs{\p_su_i(s,t)}+\Abs{\Nabla{s}v_i(s,t)}
     \le\frac{c_5}{1+s^2}
$$
for every $(s,t)\in\R\times S^1$ and every $i\in\N$. 
Moreover, it follows from Theorem~\ref{thm:2nd-derivs} that 
$$
     \Norm{\p_su_i-\Nabla{t}\p_tu_i-\grad\Vv(u_i)}_p
     = \eps_i^2 \Norm{\Nabla{t}\Nabla{s} v_i}_p
     \le c_6\eps_i^2.
$$
for a suitable constant $c_6>0$. 
Now let $\delta_0=\delta_0(p,c_5)$ and $c=c(p,c_5)$ 
be the constants in the parabolic implicit 
function theorem~\ref{thm:par-mod-IFT}.
Then the function $u_i$ satisfies the hypotheses
of Theorem~\ref{thm:par-mod-IFT}, whenever
$c_6\eps_i^2<\delta_0$. Hence, for $i$ 
sufficiently large, there is a parabolic cylinder 
$u_i^0\in\Mm^0(x^-,x^+;\Vv)$ and a vector field 
$\xi_i\in\Om^0(\R\times S^1,(u^0_i)^*TM)$ such that
$$
     u_i=exp_{u_i^0}(\xi_i),
$$
$$
     \Norm{\xi_i}_{\Ww_{u_i^0}}
     \le c_7\Norm{\p_su_i-\Nabla{t}\p_tu_i-\grad\Vv(u_i)}_p
     \le c_6c_7\eps_i^2.
$$
By the Sobolev embedding theorem, we have
$$
     \Norm{\xi_i}_\infty
     \le c_8\Norm{\xi_i}_{\Ww_{u_i^0}}
     \le c_6c_7c_8\eps_i^2
$$
for large $i$. Moreover, by definition of 
the $\Ww_{u_i^0}$-norm we have 
$$
     \Norm{\xi_i}_p +  \Norm{\Nabla{t}\xi_i}_p
     \le 2\Norm{\xi_i}_{\Ww_{u_i^0}}
     \le 2c_6c_7\eps_i^2.
$$
Hence $\xi_i$ satisfies~(\ref{eq:estimate-zeta-i-0-new}).
This proves Step~2 and the theorem.
\end{proof}

\begin{corollary}\label{cor:one-to-one}
Assume $\Ss_\Vv$ is Morse--Smale and fix a regular
value $a$ of $\Ss_\Vv$.  Then there is a constant 
$\eps_0>0$ such that, for every $\eps\in(0,\eps_0]$, 
the following holds.
\begin{enumerate}
\item[\rm\bfseries(i)]
If $x^\pm\in\Pp^a(\Vv)$ have index difference less than or equal 
to zero and $x^+\ne x^-$ then $\Mm^\eps(x^-,x^+;\Vv)=\emptyset$. 
\item[\rm\bfseries(ii)]
If $x^\pm\in\Pp^a(\Vv)$ have index difference one then
$$
    \#\Mm^0(x^-,x^+;\Vv)/\R = \#\Mm^\eps(x^-,x^+;\Vv)/\R.
$$
\item[\rm\bfseries(iii)]
If $x^\pm\in\Pp^a(\Vv)$ have index difference one and 
$(u,v)\in\Mm^\eps(x^-,x^+;\Vv)$ then $\Dd^\eps_{u,v}$ 
is surjective. 
\end{enumerate}
\end{corollary}

\begin{proof}
Assertion~(i) follows from Lemma~\ref{le:catenation}.
Assertion~(ii) follows from Theorems~\ref{thm:existence}
and~\ref{thm:onto}.  Assertion~(iii) 
follows from Theorems~\ref{thm:existence}, 
\ref{thm:inverse}, and~\ref{thm:onto}. 
\end{proof}


\section{Proof of the main result}\label{sec:proof}

\begin{theorem}\label{thm:mainZ2}
The assertion of Theorem~\ref{thm:main} 
holds with $\Z_2$-coefficients.
\end{theorem}

\begin{proof}
Let $V_t$ be a potential such that
$\Ss_V$ is a Morse function on the loop space
and denote 
$$
\Vv(x):=\int_0^1V_t(x(t))\,dt.
$$ 
Fix a regular value $a$ of $\Ss_V$. 
Choose a sequence of perturbations
${\Vv_i:\Ll M\to\R}$, converging to $\Vv$ in the 
$\Cinf$ topology, such that $\Ss_{\Vv_i}:\Ll M\to\R$
is Morse--Smale for every $i$. We may assume without 
loss of generality that the perturbations agree with 
$\Vv$ near the critical points and that $\Pp(\Vv_i)=\Pp(V)$
for all $i$. Let $\eps_i>0$ be the constant 
of Corollary~\ref{cor:one-to-one} for $\Vv=\Vv_i$. 
Then, by Corollary~\ref{cor:one-to-one},
$$
    \#\Mm^0(x^-,x^+;\Vv_i)/\R = \#\Mm^{\eps_i}(x^-,x^+;\Vv_i)/\R
$$
for every pair $x^\pm\in\Pp^a(V)$ with index difference one. 
Hence the Floer boundary operator on the chain complex 
$$
C^a(V;\Z_2) := \bigoplus_{x\in\Pp^a(V)}\Z_2 x,
$$
defined by counting modulo $2$ the solutions 
of~(\ref{eq:floer-V}) with $\Vv=\Vv_i$ and $\eps=\eps_i$
agrees with the Morse boundary operator defined
by counting the solutions of~(\ref{eq:heat-V})
with $\Vv=\Vv_i$.  Let us denote the resulting 
Floer homology groups by $\HF^a_*(T^*M,\Vv_i,\eps_i;\Z_2)$.  
Then, by what we have just observed, there is a natural
isomorphism
$$
\HF^a_*(T^*M,\Vv_i,\eps_i;\Z_2)
\cong \HM^a_*(\Ll M,\Ss_{\Vv_i};\Z_2)
\cong\mathrm{H}_*(\{\Ss_{\Vv_i}\le a\};\Z_2).
$$
Here the last isomorphism follows from 
Theorem~\ref{thm:morse}. The assertion of 
Theorem~\ref{thm:main} with $\Z_2$ coefficients 
now follows from the isomorphisms
$$
\HF^a_*(T^*M,H_V;\Z_2)
\cong\HF^a_*(T^*M;\Vv_i,\eps_i;\Z_2)
$$
and
$$
\mathrm{H}_*(\{\Ss_{\Vv_i}\le a\};\Z_2)
\cong\mathrm{H}_*(\{\Ss_V\le a\};\Z_2)
$$
for $i$ sufficiently large. 
Here the second isomorphism follows by varying 
the level $a$ and noting that the inclusions
$\{\Ss_V\le a\}\hookrightarrow\{\Ss_{\Vv_i}\le b\}
\hookrightarrow\{\Ss_V\le c\}$
are homotopy equivalences for $a<b<c$, 
$c$ sufficiently close to $a$, 
and $i$ sufficiently large. 
To understand the isomorphism on Floer 
homology, we first recall that the Floer homology
groups $\HF^a_*(T^*M,H_V;\Z_2)$ (for a nonregular Hamiltonian
$H_V$ and a regular value $a$ of the symplectic action $\Aa_V$)
are defined in terms of almost complex structures $J$
and nearby Hamiltonian functions $H$, such that $(J,H)$
is a regular pair in the sense of Floer; one then defines
$\HF^a_*(T^*M,H_V;\Z_2):=\HF^a_*(T^*M,H,J;\Z_2)$
and observes that the resulting Floer homology groups
are independent of $J$ and of the nearby Hamiltonian $H$.
Now let $J=J_{\eps_i}$ be the almost complex structure
of Remark~\ref{rmk:eps} and choose a $J_{\eps_i}$-regular
Hamiltonian $H=H_V+W$ with $W$ sufficiently close to zero.
Then the Floer equation for the pair
$(J_{\eps_i},H)$ can be written in the form
\begin{equation}\label{eq:floer-W}
\p_su + \Nabla{t}v = \nabla V_t(u) + \Nabla{1}W_t(u,v),\qquad
\eps_i^2\Nabla{s}v + \p_tu = v + \Nabla{2}W_t(u,v).
\end{equation}
Now the standard Floer homotopy argument can be used 
to relate the Floer complex associated to~(\ref{eq:floer-W})
to that of 
\begin{equation}\label{eq:floer-Vi}
\p_su + \Nabla{t}v = \grad\Vv_i(u),\qquad
\eps_i^2\Nabla{s}v + \p_tu = v.
\end{equation}
This shows that $\HF^a_*(T^*M,H_V;\Z_2)$
is isomorphic to $\HF^a_*(T^*M,\Vv_i,\eps_i;\Z_2)$
for $i$ sufficiently large. 
This proves Theorem~\ref{thm:main} 
with $\Z_2$ coefficients.
\end{proof}

To prove the result with integer coefficients it remains 
to examine the orientations of the moduli spaces.
The first step is a result about abstract Fredholm operators
on Hilbert spaces.

Let $W\subset H$ be an inclusion of Hilbert spaces 
that is compact and has a dense image. Let 
$\R\mapsto\Ll(W,H):s\mapsto A(s)$ be a family of 
bounded linear operators satisfying the following 
conditions.
\begin{enumerate}
\item[\rm\bfseries(A1)]
The map $s\mapsto A(s)$ is continuously 
differentiable in the norm topology.  
Moreover, there is a constant $c>0$ such that 
$$
      \left\|A(s)\xi\right\|_H
      + \|\dot A(s)\xi\|_H
      \le c\left\|\xi\right\|_W
$$
for every $s\in\R$ and every $\xi\in W$.
\item[\rm\bfseries(A2)]
The operators $A(s)$ are uniformly self-adjoint.
This means that, for each $s$, the operator $A(s)$,
when considered as an unbounded operator
on $H$, is self adjoint, and there is a constant 
$c$ such that 
$$
     \left\|\xi\right\|_W\le c\left(
     \left\|A(s)\xi\right\|_H+\left\|\xi\right\|_H
     \right)
$$
for every $s\in\R$ and every $\xi\in W$.
\item[\rm\bfseries(A3)]
There are invertible operators $A^\pm:W\to H$
such that 
$$
     \lim_{s\to\pm\infty}\left\|A(s)-A^\pm\right\|_{\Ll(W,H)}
     = 0.
$$
\item[\rm\bfseries(A4)]
The operator $A(s)$ has finitely many negative 
eigenvalues for every $s\in\R$. 
\end{enumerate}
Denote by $\Ss(W,H)$ the set of invertible 
self-adjoint operators $A:W\to H$ 
with finitely many negative eigenvalues.
For $A\in\Ss(W,H)$ denote by $E(A)$ the direct sum of the 
eigenspaces of $A$ with negative eigenvalues.
Given $A^\pm\in\Ss(W,H)$ denote by $\Pp(A^-,A^+)$ the set of
functions $A:\R\to\Ll(W,H)$ that satisfy~(A1-4) and 
by $\Pp$ the union of the spaces $\Pp(A^-,A^+)$ 
over all pairs $A^\pm\in\Ss(W,H)$.  This is an open subset
of a Banach space. 

Denote 
$$
    \Ww := L^2(\R,W)\cap W^{1,2}(\R,H),\qquad
    \Hh := L^2(\R,H)
$$
and, for every pair $A^\pm\in\Ss(W,H)$ and every
$A\in\Pp(A^-,A^+)$, consider the operator 
$\Dd_A:\Ww\to\Hh$ defined by
$$
    (\Dd_A\xi)(s) := \dot\xi(s)+A(s)\xi(s)
$$
for $\xi\in\Ww$. 
This operator is Fredholm and its index 
is the spectral flow, i.e.
$$
    \INDEX(\Dd_A) = \dim E(A^-) - \dim E(A^+)
$$
(see Robbin--Salamon~\cite{ROSA}). 
The formal adjoint operator $\Dd_A^*:\Ww\to\Hh$ is given 
by $\Dd_A^*\eta=-\dot\eta+A\eta$. Denote by 
$$
    \det(\Dd_A) := \Lambda^{\rm max} \left(\ker\,\Dd_A\right)
       \otimes \Lambda^{\rm max} \left(\ker\,(\Dd_A)^*\right)
$$
the determinant line of $\Dd_A$ and by $\Or(\Dd_A)$
the set of orientations of $\det(\Dd_A)$.
For $A\in\Ss(W,H)$ denote by $\Or(A)$ the set of 
orientations of $E(A)$. 

\begin{remark}[The finite dimensional case]\label{rmk:fdor}\rm
Assume $W=H=\R^n$. Let $A^\pm$ be nonsingular symmetric 
$(n\times n)$-matrices and $A\in\Pp(A^-,A^+)$.
Suppose that $A(s)=A^\pm$ for $\pm s\ge T$.
Define $\Phi(s,s_0)\in\R^{n\times n}$ by 
$$
    \p_s\Phi(s,s_0)+A(s)\Phi(s,s_0)=0,\qquad
    \Phi(s_0,s_0)=\1.
$$
Define 
$$
    E^\pm(s) := \left\{\xi\in\R^n\,|\,
    \lim_{r\to\pm\infty}\Phi(r,s)\xi=0\right\}.
$$
Then $E^-(s)=E(A^-)$ for $s\le -T$ and
$E^+(s)=E(A^+)^\perp$ for $s\ge T$. Moreover, 
$$
    \ker\Dd_A\cong E^-(s)\cap E^+(s),\qquad
    (\im\Dd_A)^\perp\cong (E^-(s)+E^+(s))^\perp.
$$
Hence there is a natural map 
$$
    \tau_A:\Or(A^-)\times\Or(A^+)\to\Or(\Dd_A)
$$
defined as follows.  Given orientations of $E(A^-)\cong E^-(s)$
and $E(A^+)\cong E^+(s)^\perp$, pick any basis 
$u_1,\dots,u_\ell$ of $E^-(s)\cap E^+(s)\cong\ker\,\Dd_A$.
Extend it to a positive basis of $E^-(s)$ by picking a 
suitable basis $v_1,\dots,v_m$ of $E^-(s)\cap E^+(s)^\perp$.
Now extend the vectors $v_j$ to a positive basis 
of $E^+(s)^\perp$ by picking a suitable basis
$w_1,\dots,w_n$ of 
$(E^-(s)+E^+(s))^\perp\cong(\im\,\Dd_A)^\perp$.
Then the bases $u_1,\dots,u_\ell$ of $ker\,\Dd_A$
and $w_1,\dots,w_n$ of $(\im\,\Dd_A)^\perp$
determine the induced orientation of $\det(\Dd_A)$.
Note that this is well defined (a sign change in the 
$u_i$ leads to a sign change in the $w_k$).
\end{remark}

\begin{remark}[Catenation]\label{rmk:catenation}\rm
Let $A_0,A_1,A_2\in\Ss(W,H)$ and suppose that 
$A_{01}\in\Pp(A_0,A_1)$ and $A_{12}\in\Pp(A_1,A_2)$ 
satisfy
\begin{equation}\label{eq:012}
    A_{01}(s) = \begin{cases}
      A_0&\mbox{if }s\le -T, \\
      A_1&\mbox{if }s\ge T,
     \end{cases}\qquad
    A_{12}(s) = \begin{cases}
      A_1&\mbox{if }s\le -T, \\
      A_2&\mbox{if }s\ge T.
     \end{cases}
\end{equation}
For $R>T$ define $A^R_{02}\in\Pp(A_0,A_2)$ by
\begin{equation}\label{eq:catenation}
    A^R_{02}(s) = \begin{cases}
      A_{01}(s+R)&\mbox{if }s\le0, \\
      A_{12}(s-R)&\mbox{if }s\ge0.
     \end{cases}
\end{equation}
If $\Dd_{A_{01}}$ and $\Dd_{A_{12}}$ are onto then,
for $R$ sufficiently large, the operator $\Dd_{A^R_{02}}$ 
is onto and there is a natural isomorphism
$$
    S^R:\ker\,\Dd_{A_{01}}\oplus\ker\,\Dd_{A_{12}}
    \to\ker\,\Dd_{A^R_{02}}
$$
The isomorphism $S^R$ is defined by composing a pre-gluing operator
with the orthogonal projection onto the kernel. That this gives an 
isomorphism follows from exponential decay estimates for the elements
in the kernel and a uniform estimate for suitable right inverses
of the operators $\Dd_{A^R_{02}}$ (see for example~\cite{Sa97}). 
\end{remark}

\begin{theorem}\label{thm:or}
There is a family of maps
$$
    \tau_A:\Or(A^-)\times\Or(A^+)\to\Or(\Dd_A),
$$
one for each pair of Hilbert spaces $W\subset H$
with a compact dense inclusion, each pair
$A^\pm\in\Ss(W,H)$, and each $A\in\Pp(A^-,A^+)$,
satisfying the following axioms.
\begin{description}
\item[(Equivariant)]
$\tau_A$ is equivariant with respect
to the $\Z_2$-action on each factor. 
\item[(Homotopy)]
The map 
$(A,o^-,o^+)$$\mapsto$$(A,\tau_A(o^-,o^+))$
from the topological space
$
\{(A,o^-,o^+)\,|\,A\in\Pp,\,o^\pm\in\Or(A^\pm)\}
$
to
$
\{(A,o)\,|\,A\in\Pp,\,o\in\Or(\Dd_A)\}
$
is continuous.
\item[(Naturality)]
Let $\Phi(s):(W,H)\to(W',H')$ be a family of 
(pairs of) Hilbert space isomorphisms that 
is continuously differentiable in the operator 
norm on $H$ and continuous in the operator norm 
on $W$.  Suppose that there exist
Hilbert space isomorphisms $\Phi^\pm:(W,H)\to(W',H')$
such that $\Phi(s)$ converges to $\Phi^\pm$ in the operator 
norm on both spaces and $\dot\Phi(s)$ converges to zero in 
$\Ll(H)$ as $s\to\pm\infty$.  Then
$$
      \tau_{\Phi_*A}(\Phi^-_*o^-,\Phi^+_*o^+)
      = \Phi_*\tau_A(o^-,o^+)
$$
for all $A^\pm\in\Ss(W,H)$, $A\in\Pp(A^-,A^+)$, 
and $o^\pm\in\Or(A^\pm)$.
\item[(Direct Sum)]
If $A^\pm_j\in\Ss(W_j,H_j)$ and $A_j\in\Pp(A_j^-,A_j^+)$
for $j=0,1$ then 
$$
    \tau_{A_0\oplus A_1}(o_0^-\otimes o_1^-,o_0^+\otimes o_1^+)
    = \tau_{A_0}(o_0^-,o_0^+)\otimes\tau_{A_1}(o_1^-,o_1^+).
$$
for all $o_j^\pm\in\Or(A_j^\pm)$. 
\item[(Catenation)]
Let $A_0,A_1,A_2\in\Ss(W,H)$, suppose that 
$A_{01}\in\Pp(A_0,A_1)$ and $A_{12}\in\Pp(A_1,A_2)$ 
satisfy~(\ref{eq:012}) and, for $R>T$, define 
$A^R_{02}$ by~(\ref{eq:catenation}). 
Assume $\Dd_{A_{01}}$ and $\Dd_{A_{12}}$ are onto.
Then $\Dd_{A^R_{02}}$ is onto for large $R$ and
$$
    \tau_{A_{02}}(o_0,o_2)
    = \sigma^R\left(\tau_{A_{01}}(o_0,o_1),
      \tau_{A_{12}}(o_1,o_2)\right).
$$
for $o_0\in\Or(A_0)$, $o_1\in\Or(A_1)$, and $o_2\in\Or(A_2)$. 
Here the map 
$$
    \sigma^R:\det(\Dd_{A_{01}})\times\det(\Dd_{A_{12}})
    \to\det(\Dd_{A^R_{02}})
$$
is induced by the isomorphism $S^R$ of 
Remark~\ref{rmk:catenation}.
\item[(Constant)]
If $A(s)\equiv A^+=A^-$ and $o^+=o^-\in\Or(A^\pm)$
then $\tau_A(o^-,o^+)$ is the standard orientation of
$\det(\Dd_A)\cong\R$.
\item[(Normalization)]
If $W=H=\R^n$ then $\tau_A$ is the 
map defined in Remark~\ref{rmk:fdor}. 
\end{description}
The maps $\tau_A$ are uniquely determined
by the {\it (Homotopy)}, {\it (Direct Sum)}, 
{\it (Constant)}, and {\it (Normalization)} axioms.
\end{theorem}

Theorem~\ref{thm:or} is standard with the techniques of~\cite{FH} 
(although the assumptions are not quite the same as in the 
work of Floer and Hofer).

\begin{proof}[Proof of Theorem~\ref{thm:main}.]
Assume $\Ss_\Vv$ is Morse--Smale.
For $x\in\Pp(\Vv)$ denote by $W^u(x)$ the unstable 
manifold of $x$ with respect to the negative gradient 
flow of $\Ss_\Vv$. Thus $W^u(x)$ is the space of all smooth 
loops $y:S^1\to M$ such that there exists a solution
$u:(-\infty,0]\times S^1\to M$ of the nonlinear heat 
equation~(\ref{eq:heat}) that converges to $x$ as $s\to-\infty$
and satisfies $u(0,t)=y(t)$.  Then $W^u(x)$ is a finite dimensional
manifold (see for example~\cite{DAVIES}). It is diffeomorphic
to $\R^k$ where $k=\IND_\Vv(x)$ is the Morse index of $x$
as a critical point of $\Ss_\Vv$. Fix an orientation of $W^u(x)$
for every periodic orbit $x\in\Pp(\Vv)$. 
These orientations determine a system of coherent 
orientations for the heat flow as follows.

Fix a pair $x^\pm\in\Pp(\Vv)$ of periodic orbits 
that represent the same component of $\Ll M$.
Denote by $\Pp^0(x^-,x^+)$ the set of smooth maps
$u:\R\times S^1\to M$ such that $u(s,\cdot)$ 
converges to $x^\pm$ in the $C^2$ norm and
$\p_su(s,\cdot)$ converges to zero in the $C^1$ norm 
as $s$ tends to $\pm\infty$. Then, in a suitable 
trivialization of the tangent bundle $u^*TM$, 
the linearized operator $\Dd^0_u$ has the form of 
an operator $\Dd_A$ as in Theorem~\ref{thm:or}
where the spaces $E(A^\pm)$ correspond to the tangent
spaces $T_{x^\pm}W^u(x^\pm)$ of the unstable manifolds.
Hence, by Theorem~\ref{thm:or}, the given orientations
of the unstable manifolds determine orientations 
$$
     \nu^0(u)\in\Or(\det(\Dd^0_u))
$$
of the determinant lines for all $u\in\Pp^0(x^-,x^+)$
and all $x^\pm\in\Pp(\Vv)$. By the {\it (Naturality)} axiom,
these orientations are independent 
of the choice of the trivializations used to define them. 
By the {\it (Catenation)} axiom, they form a system 
of {\it coherent orientations} in the sense of 
Floer--Hofer~\cite{FH}. 

Next we show how the coherent orientations for the 
heat flow induce a system of coherent orientations 
$$
     \nu^\eps(u,v)\in\Or(\det(\Dd^\eps_{u,v}))
$$
for the Floer equations~(\ref{eq:floer-V}). 
Let us denote by $\Pp(x^-,x^+)$ the set of smooth maps
$(u,v):\R\times S^1\to TM$ such that 
$(u(s,\cdot),v(s,\cdot))$ converges 
to $(x^\pm,\dot x^\pm)$ in the $C^1$ norm 
and $(\p_su,\Nabla{s}v)$ converges to zero,
uniformly in $t$, as $s$ tends to $\pm\infty$.
By the obvious homotopy arguments it suffices to 
assume $u\in\Pp^0(x^-,x^+)$ and $v=\p_tu$. 
We abbreviate 
$$
\Dd^\eps_u:=\Dd^\eps_{u,\p_tu}.
$$
It follows from the definition of the operators
in~(\ref{eq:D-eps}) that 
$$
    \Dd^0_u\xi=0\qquad\IMP\qquad
    \Dd^\eps_u\begin{pmatrix}
     \xi \\ \Nabla{t}\xi \end{pmatrix}
     = \begin{pmatrix}
       0 \\
       \Nabla{s}\Nabla{t}\xi + R(\xi,\p_su)\p_tu
       \end{pmatrix}.
$$
Hence $\Dd^\eps_u(\xi,\Nabla{t}\xi)$
is small in the $(0,2,\eps)$-norm.  If the operator
$\Dd^0_u$ is onto then the estimate of 
Theorem~\ref{thm:inverse} shows that the map
$$
    \ker\,\Dd^0_u\to\ker\,\Dd^\eps_u:
    \xi\mapsto
    \begin{pmatrix}
    \xi \\ \Nabla{t}\xi \end{pmatrix}
    - {\Dd^\eps_u}^*
    \left(\Dd^\eps_u{\Dd^\eps_u}^*\right)^{-1}
    \Dd^\eps_u\begin{pmatrix}
    \xi \\ \Nabla{t}\xi \end{pmatrix}
$$
is an isomorphism between the kernels and we define 
$\nu^\eps(u,\p_tu)$ to be the image of $\nu^0(u)$
under the induced isomorphism of the top exterior powers.
If $\Dd^0_u$ is not onto we obtain a similar isomorphism
between the determinant lines of $\Dd^0_u$ and $\Dd^\eps_u$
by augmenting the operators first to make them surjective. 
It follows again from the {\it (Catenation)} axiom that
the $\nu^\eps(u,v)$ form a system of coherent orientations
for the Floer equations. 

Now assume that $x^\pm\in\Pp(\Vv)$ have Morse index 
difference one. Consider the map
$$
    \Tt^\eps:\Mm^0(x^-,x^+;\Vv)\to\Mm^\eps(x^-,x^+;\Vv)
$$
of Definition~\ref{def:T} and recall that,
by Theorem~\ref{thm:onto}, it is bijective. 
It follows from the proof of Theorem~\ref{thm:time-shift} 
that the map $\Tt^\eps$ satisfies the following.  
Let $u\in\Mm^0(x^-,x^+;\Vv)$ and 
$$
(u^\eps,v^\eps):=\Tt^\eps(u)\in\Mm^\eps(x^-,x^+;\Vv).
$$
Then the vector $\p_su\in\ker\,\Dd^0_u$ is positively 
oriented with respect to $\nu^0(u)$ if and only if the 
vector $(\p_su^\eps,\Nabla{s}v^\eps)\in\ker\Dd^\eps_{u^\eps,v^\eps}$
is positively oriented with respect to $\nu^\eps(u^\eps,v^\eps)$.
Hence the bijection $\Tt^\eps$ preserves the signs 
for the definitions of the two boundary operators.
This shows that the Morse complex of the heat flow 
has the same boundary operator as the Floer 
complex for $\eps$ sufficiently small. Hence 
the resulting homologies are naturally isomorphic, 
i.e. for every regular value $a$ of $\Ss_\Vv$
there is a constant $\eps_0>0$ such that
$$
     \HF^a_*(T^*M,\Vv,\eps;\Z)\cong\HM^a_*(\Ll M,\Ss_\Vv;\Z)
$$
for $0<\eps\le\eps_0$. In fact, we have established this 
isomorphism on the chain level and with integer 
coefficients. To complete the proof of 
Theorem~\ref{thm:main} one can now argue as in the proof 
of Theorem~\ref{thm:mainZ2} to show that, given a
potential $V$ such that $\Ss_V$ is Morse
and a regular value $a$ of $\Ss_V$, we have two 
isomorphisms 
$$
\HF^a_*(T^*M,H_V;\Z)\cong\HF^a_*(T^*M,\Vv,\eps;\Z)
$$
and 
$$
\HM^a_*(\Ll M,\Ss_\Vv;\Z)\cong\mathrm{H}_*(\{\Ss_V\le a\};\Z)
$$
for a suitable perturbation $\Vv$ and $\eps>0$ sufficiently
small.  This proves the result for integer 
coefficients and $a<\infty$.  The argument for general 
coefficient rings is exactly the same.  The result for 
$a=\infty$ follows by taking the direct limit $a\to\infty$
and noting that there are natural isomorphisms
$$
\HF_*(T^*M,H_V)\cong\underset{a\in\R}\varinjlim\;\HF^a_*(T^*M,H_V)
$$
and 
$$
\mathrm{H}_*(\Ll M)\cong
\underset{a\in\R}\varinjlim\;\mathrm{H}_*(\{\Ss_V\le a\}).
$$
This proves Theorem~\ref{thm:main}.
\end{proof}


\appendix

\section{The heat flow}\label{app:MS}

In this appendix we summarize results 
from~\cite{WEBER} that are used in this paper.
We assume throughout this appendix that $M$ is a 
closed Riemannian manifold. Let $\Vv:\Ll M\to\R$ 
be a smooth function that satisfies the axioms~$(V0-V4)$. 
Consider the action functional
$$
     \Ss_\Vv(x) = \frac12\int_0^1 
     \abs{\dot x(t)}^2\,dt - \Vv(x)
$$
and the corresponding heat equation
\begin{equation}\label{eq:HEAT}
\p_su-\Nabla{t}\p_tu-\grad\Vv(u)=0
\end{equation}
for smooth functions 
$\R\times S^1\to M:(s,t)\mapsto u(s,t)$. 
In the following we denote by $\Pp(\Vv)\subset\Cinf(S^1,M)$
the set of critical points $x$ of $\Ss_\Vv$ (i.e. of solutions
of the equation $\Nabla{t}\dot x+\grad\Vv(x)=0$), and by 
$\Pp^a(\Vv)$ the set of all $x\in\Pp(\Vv)$ with action 
$\Ss_\Vv(x)\le a$. For two nondegenerate critical points 
$x^\pm\in\Pp(\Vv)$ we denote by $\Mm^0(x^-,x^+;\Vv)$ the set of 
all solutions $u$ of~(\ref{eq:HEAT}) that converge to $x^\pm(t)$ as 
$s\to\pm\infty$.  The energy of such a solution
is given by 
$$
     E(u) := \int_{-\infty}^\infty\int_0^1\Abs{\p_su}^2\,dtds
     = \Ss_\Vv(x^-)-\Ss_\Vv(x^+).
$$

\begin{theorem}[Apriori estimates]
\label{thm:par-apriori}
Fix a perturbation $\Vv:\Ll M\to\R$ that 
satisfies~$(V0-V1)$ and a constant $c_0>0$. 
Then there is a constant
$C=C(c_0,\Vv)>0$ such that the following holds.
If $u:\R\times S^1\to M$ is a solution
of~(\ref{eq:HEAT}) such that
$\Ss_\Vv(u(s,\cdot))\le c_0$
for every $s\in\R$ then
$$
     \left\|\p_su\right\|_\infty
     +\left\|\p_tu\right\|_\infty
     +\left\|\Nabla{t}\p_tu\right\|_\infty
     \le C.
$$
\end{theorem}

\begin{theorem}[Exponential decay]
\label{thm:par-exp-decay}
Fix a perturbation $\Vv:\Ll M\to\R$ that 
satisfies~$(V0-V4)$ and assume $\Ss_\Vv$ is Morse. 
\begin{enumerate}
\item[\rm\bfseries(F)]
Let $u:[0,\infty)\times S^1\to M$ be a solution 
of~(\ref{eq:HEAT}).  Then there are positive 
constants $\rho$ and $c_1,c_2,c_3,\dots$
such that 
$$
     \Norm{\p_su}_{C^k([T,\infty)\times S^1)}
     \le c_ke^{-\rho T}
$$
for every $T\ge1$. Moreover, there is 
a periodic orbit $x\in\Pp(\Vv)$ such that
$u(s,t)$ converges to $x(t)$ as $s\to\infty$. 
\item[\rm\bfseries(B)]
Let $u:(-\infty,0]\times S^1\to M$ be a solution 
of~(\ref{eq:HEAT}) with finite energy. 
Then there are positive 
constants $\rho$ and $c_1,c_2,c_3,\dots$
such that 
$$
     \Norm{\p_su}_{C^k((-\infty,-T]\times S^1)}
     \le c_ke^{-\rho T}
$$
for every $T\ge1$. Moreover, there is 
a periodic orbit $x\in\Pp(\Vv)$ such that
$u(s,t)$ converges to $x(t)$ as $s\to-\infty$. 
\end{enumerate}
\end{theorem}

\begin{theorem}[Regularity]
\label{thm:par-regularity}
Fix a constant $p>2$ and a perturbation 
$\Vv:\Ll M\to\R$ that satisfies~$(V0-V4)$.
Let $u:\R\times S^1\to M$ be a continuous function
which is locally of class $W^{1,p}$. Assume further
that $u$ is a weak solution of~(\ref{eq:HEAT}).
Then $u$ is smooth. 
\end{theorem}

The covariant Hessian of $\Ss_\Vv$ at a loop $x:S^1\to M$
is the operator $A(x):W^{2,2}(S^1,x^*TM)\to L^2(S^1,x^*TM)$,
given by
$$
A(x)\xi
:= - \Nabla{t}\Nabla{t}\xi 
- R(\xi,\dot x)\dot x - \Hh_\Vv(x)\xi.
$$
This operator is self-adjoint with respect to 
the standard $L^2$ inner product on $\Om^0(S^1,x^*TM)$. 
In this notation the linearized operator 
$\Dd_u^0:\Ww_u^p\to\Ll_u^p$ is given by 
$$
\Dd_u^0\xi := \Nabla{s}\xi + A(u_s)\xi
$$
where $u_s(t):=u(s,t)$.
(See Section~\ref{sec:linear} for the definition of the spaces 
$\Ww_u=\Ww_u^p$ and $\Ll_u=\Ll_u^p$.)

\begin{theorem}[Fredholm]
\label{thm:par-Fredholm}
Fix a perturbation $\Vv:\Ll M\to\R$ that 
satisfies~$(V0-V4)$ and assume $\Ss_\Vv$ is Morse. 
Let $x^\pm\in\Pp(\Vv)$ and $u:\R\times S^1\to M$
be a smooth map such that $u(s,\cdot)$ converges 
to $x^\pm$ in the $C^2$ norm and $\p_su$
converges uniformly to zero as
$s\to\pm\infty$. Then, for every $p>1$,
the operator $\Dd_u^0:\Ww_u^p\to\Ll_u^p$ is Fredholm
and its Fredholm index is given by
$$
\INDEX\,\Dd_u^0 = \IND_\Vv(x^-)-\IND_\Vv(x^+).
$$
Here $\IND_\Vv(x^\pm)$ denotes the Morse
index of $x^\pm$, i.e. the number of negative 
eigenvalues of $A(x^\pm)$. 
\end{theorem}

\begin{theorem}[Implicit function theorem]
\label{thm:par-mod-IFT}
Fix a perturbation $\Vv:\Ll M\to\R$ that satisfies~$(V0-V4)$.
Assume $\Ss_\Vv$ is Morse and that $\Dd_u^0$ is onto
for every $u\in\Mm^0(x^-,x^+;\Vv)$ and every pair
$x^\pm\in\Pp^a(\Vv)$. Fix two critical points 
$x^\pm\in\Pp^a(\Vv)$ with Morse index difference one.
Then, for all $c_0>0$ and $p>2$,
there exist positive constants $\delta_0$ and $c$
such that the following holds.
If $u:\R\times S^1\to M$ is a smooth map
such that $\lim_{s\to\pm\infty}u(s,\cdot)=x^\pm(\cdot)$
exists, uniformly in $t$, and such that
$$
\Abs{\p_su(s,t)}\le\frac{c_0}{1+s^2},\qquad
\Abs{\p_tu(s,t)}\le c_0
$$
for all $(s,t)\in\R\times S^1$ and 
$$
\Norm{\p_su-\Nabla{t}\p_tu-\grad\Vv(u)}_p\le\delta_0.
$$
Then there exist elements $u_0\in\Mm^0(x^-,x^+;\Vv)$
and $\xi\in\im(\Dd_{u_0}^0)^*\cap\Ww_{u_0}$
satisfying
$$
     u=\exp_{u_0}(\xi),\qquad
     \Norm{\xi}_{\Ww_{u_0}}
     \le c\Norm{\p_su-\Nabla{t}\p_tu-\grad\Vv(u)}_p.
$$
\end{theorem}

\begin{theorem}[Transversality]\label{thm:par-transversality}
For a generic perturbation $\Vv:\Ll M\to\R$ satisfying~$(V0-V4)$
the function $\Ss_\Vv:\Ll M\to\R$ is {\bf Morse--Smale} in the sense
that every critical point $x$ of $\Ss_\Vv$ is nondegenerate
(i.e. the Hessian $A(x)$ is bijective) and every finite energy
solution $u:\R\times S^1\to M$ of~(\ref{eq:HEAT}) is regular
(i.e. the Fredholm operator $\Dd^0_u$ is surjective).
\end{theorem}

\begin{theorem}\label{thm:morse}
Let $\Vv:\Ll M\to\R$ be a perturbation that satisfies~$(V0-V4)$
and assume that $\Ss_\Vv$ is Morse--Smale.
Then, for every regular value $a$ of $\Ss_\Vv$ and every
principal ideal domain $R$, there is a natural isomorphism
$$
\HM^a_*(\Ll M,\Ss_\Vv;R) \cong\mathrm{H}_*(\Ll^aM;R),\qquad
\Ll^aM:=\{x\in\Ll M\,|\,\Ss_\Vv(x)\le a\}.
$$
If $M$ is not simply connected then there is a separate
isomorphism for each component of the loop space.
The isomorphism commutes with the homomorphisms
$\HM^a_*(\Ll M,\Ss_\Vv)\to\HM^b_*(\Ll M,\Ss_\Vv)$
and $\mathrm{H}_*(\Ll^aM)\to\mathrm{H}_*(\Ll^bM)$ for $a<b$. 
\end{theorem}

The proof of Theorem~\ref{thm:morse}
is similar to the finite dimensional case
(see~\cite{SAL,JOA0}) since the gradient 
flow of $\Ss_\Vv$ defines a wellposed
initial value problem. 


\section{Mean value inequalities}\label{app:mean}

\begin{figure}
  \centering
  \epsfig{figure=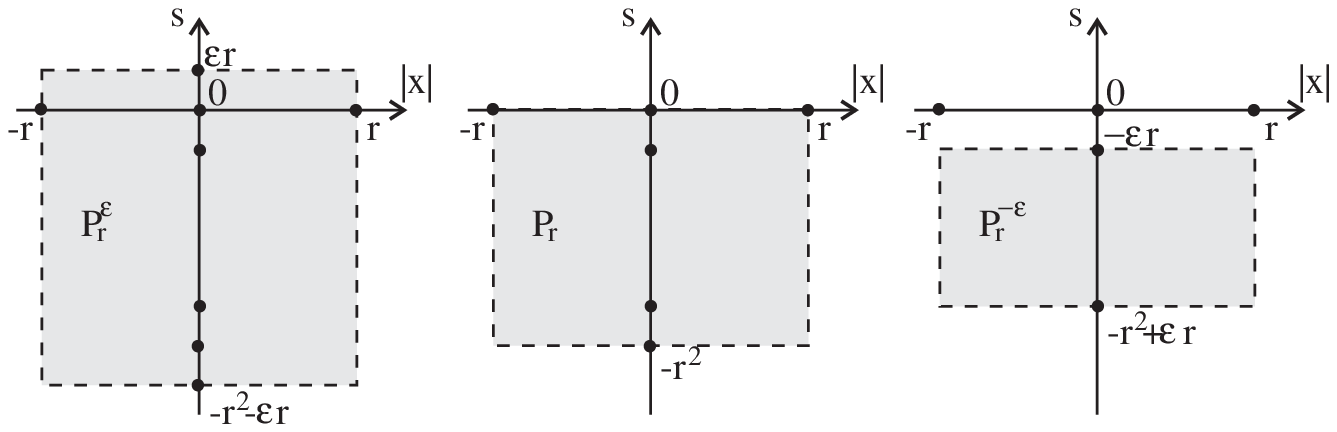,width=12cm}
  \caption{Parabolic cylinders} 
  \label{fig:fig-parabdom}
\end{figure}

Let $n$ be a positive integer and denote by
$$
     \Delta
     :={\p_1}^2+\cdots+{\p_n}^2
$$
the standard Laplacian on $\R^n$.
Given positive real numbers $r$ and $\eps$ 
let $B_r=B_r(0)$ be the open ball
of radius $r$ in $\R^n$ 
and define the \emph{parabolic cylinders}
$P_r,P_r^\eps,P_r^{-\eps}\subset\R^{n+1}$
by 
\begin{equation*}
\begin{split}
     P_r^\eps
     &:=(-r^2-\eps r,\eps r)\times B_r,\\
     P_r
     &:=(-r^2,0)\times B_r,\\
     P_r^{-\eps}
     &:=(-r^2+\eps r,-\eps r)\times B_r.
\end{split}
\end{equation*}
(see Figure~\ref{fig:fig-parabdom}).
The elements of $P_r$ are denoted
by $(s,x)=(s,x_1,\dots,x_n)$.

\begin{lemma}\label{le:apriori-basic}
For every $n\in\N$ there is a constant $c_n>0$
such that the following holds for every $r\in(0,1]$.
If $a\ge 0$ and $w:\R\times\R^n\supset P_r\to\R$
is $C^1$ in the $s$-variable and $C^2$ in the $x$-variable
such that
$$
     (\Delta-\p_s)w\ge -aw, \qquad w\ge0,
$$
then
$$
     w(0)\le\frac{c_ne^{ar^2}}{r^{n+2}}
     \int_{P_r} w.
$$
\end{lemma}

\begin{proof}
For $a=0$ this is a special case of a theorem by 
Gruber for parabolic differential operators
with variable coefficients.
(See Gruber~\cite[Theorem~2.1]{Gr84} with 
$p=1$, $\theta=1$, $\lambda=1$, $\sigma=1/2$, $R=r$ and $f=0$;
for an another proof see 
Lieberman~\cite[Theorem~7.21]{Li96} 
with $R=r$, $p=1$, $\rho=1/2$, $f=0$.)

To prove the result in general assume that $w$ satisfies
the hypotheses of the lemma and define 
$
     f(s,x):=e^{-as}w(s,x).
$
Then
$$
     (\Delta-\p_s)f=e^{-as}(\Delta-\p_s+a)w
     \ge0.
$$
Hence, by Gruber's theorem, 
$$
     w(0)=f(0) 
     \le\frac{c_n}{r^{n+2}}\int_{P_r}f
     \le\frac{c_ne^{ar^2}}{r^{n+2}}\int_{P_r}w.
$$
This proves the lemma.
\end{proof}

\begin{lemma} \label{le:apriori-basic-eps}
Let $c_2$ be the constant in Lemma~\ref{le:apriori-basic} 
with $n=2$. Let $\eps>0$, $r\in(0,1]$, and $a\ge 0$.
If $w:\R\times\R\supset P_r^\eps\to\R$ is $C^1$ 
in the $s$-variable and $C^2$ in the $t$-variable 
and satisfies
\begin{equation}\label{eq:assu-apriori-basic-eps}
     L_\eps w
     :=\left(\eps^2{\p_s}^2+{\p_t}^2-\p_s\right)
     w\ge-aw, \qquad w\ge0,
\end{equation}
then
$$
     w(0)
     \le\frac{2c_2e^{ar^2}}{r^3}\int_{P_r^\eps}w.
$$
\end{lemma}

\begin{proof}
The idea of proof was suggested to us by Tom Ilmanen.
Define a function $W$ on the domain
$P_r\subset \R \times \R^2$ by
$$
     W(s,t,q):=w(s+\eps q,t).
$$
(Note that $(s+\eps q,t)\in P_r^\eps\subset\R\times\R$
for every $(s,t,q)\in P_r\subset\R\times\R^2$.)
Then, by assumption, we have
\begin{equation*}
      \left(\Delta-\p_s\right)W(s,t,q)
      =\left(L_\eps w\right)(s+\eps q,t) 
      \ge-aw(s+\eps q,t) 
      =-aW(s,t,q),
\end{equation*}
where $\Delta:=\p_t^2+\p_q^2$. 
Hence it follows from Lemma~\ref{le:apriori-basic} with $n=2$
that
$$
     w(0) = W(0)
     \le\frac{c_2e^{ar^2}}{r^4}
     \int_{P_r}W.
$$
It remains to
estimate the integral on the right hand side:
\begin{equation}\label{eq:helper-1}
\begin{split}
     \int_{P_r}W
    &\le\int_{-r}^r\int_{-r}^r\int_{-r^2}^0W(s,t,q)
     \:dsdqdt \\
    &=\int_{-r}^r\int_{-r}^r\int_{-r^2+\eps q}
     ^{\eps q}
     w(z,t)\:dzdqdt \\
    &\le\int_{-r}^r\int_{-r}^r\int_{-r^2-\eps r}
     ^{\eps r}
     w(z,t)\:dzdqdt \\
    &=2r\int_{P_r^\eps}w.
\end{split}
\end{equation}
The first step uses the fact that $W\ge 0$ and
$B_r \subset [-r,r]\times[-r,r]$. 
The third step uses the fact that $w\ge 0$ 
and $(-r^2+\eps q,\eps q)\subset(-r^2-\eps r,\eps r)$,
since $0\le q\le r$. 
\end{proof}

\begin{lemma}\label{le:apriori-1eps}
Fix three constants $r>0$, $\eps\ge 0$, and $\mu\ge 0$. 
Let $c_2$ be the constant of Lemma~\ref{le:apriori-basic}.
If $f:[-r^2-\eps r,\eps r]\to\R$ is a $C^2$ function
satisfying
$$
\eps^2f''-f'+\mu f\ge 0,\qquad f\ge 0,
$$
then
$$
f(0) \le \frac{2c_2e^{\mu r^2}}{r^3}
\int_{-r^2-\eps r}^{\eps r}f(s)\,ds.
$$
\end{lemma}

\begin{proof}
This follows immediately from Lemma~\ref{le:apriori-basic-eps}
with $w(s,t):=f(s)$. 
\end{proof}

\begin{lemma} \label{le:apriori-L2}
Let $u:\R\times\R^n\supset P_{R+r}\to\R$ be $C^1$
in the $s$-variable 
and $C^2$ in the $x$-variable and $f,g:P_{R+r}\to\R$ 
be continous functions such that
$$
     \left(\Delta -\p_s\right)u\ge g-f, \qquad
     u\ge 0, \qquad f\ge 0, \qquad g\ge 0.
$$
Then
$$
     \int_{P_R}g
     \le\int_{P_{R+r}}f+\left(\frac{4}{r^2}
     +\frac{1}{Rr}\right)\int_{P_{R+r}\setminus P_R}u.
$$
\end{lemma}

\begin{proof}
The proof rests on the following two inequalities.
Let $B_r\subset\R^n$ be the open ball of radius
$r$ centered at zero.  Then, for every smooth function
$u:\R^n\to[0,\infty)$, we have
\begin{equation}\label{eq:ID-I}
     \int_{\p B_r}\frac{\p u}{\p\nu}
     =-\frac{n-1}{r}\int_{\p B_r}u
     +\frac{d}{dr}\int_{\p B_r}u 
     \le\frac{d}{dr}\int_{\p B_r}u
\end{equation}
(see~\cite[Theorem 2.1]{GT77}).
Secondly, every smooth function $u:\R\times\R^n\to[0,\infty)$
satisfies
\begin{equation} \label{eq:ID-II}
\begin{split}
    \frac{d}{d\sigma}\int_{-(R+\sigma)^2}^0
     \biggl(\int_{\p B_{R+\sigma}}u(s,\cdot)\biggr)ds
    &=\int_{-(R+\sigma)^2}^0
     \biggl(\frac{d}{d\sigma}\int_{\p B_{R+\sigma}}
     u(s,\cdot)\biggr)ds \\
    &\quad+2(R+\sigma)\int_{\p B_{R+\sigma}}
     u(-(R+\sigma)^2,\cdot) \\
    &\ge\int_{-(R+\sigma)^2}^0\biggl(\frac{d}{d\sigma}
     \int_{\p B_{R+\sigma}}u(s,\cdot)\biggr)ds.
\end{split}
\end{equation}
Now suppose $u$, $f$, $g$ satisfy the assumptions
of the lemma. Then, for $0\le\sigma\le r$,
\begin{equation*}
\begin{split}
    &\int_{P_R}g-\int_{P_{R+r}}f \\
    &\le\int_{P_{R+\sigma}}(\Delta u-\p_s u) \\
    &=\int_{-(R+\sigma)^2}^0\left(
     \int_{\p B_{R+\sigma}}\frac{\p u}{\p \nu}
     (s,\cdot)\right)ds
     -\int_{B_{R+\sigma}}\Bigl(u(0,\cdot)
     -u(-(R+\sigma)^2,\cdot)\Bigr)dx \\
    &\le\int_{-(R+\sigma)^2}^0\left(\frac{d}{d\sigma}
     \int_{\p B_{R+\sigma}}u(s,\cdot)\right)ds
     +\int_{B_{R+\sigma}}u(-(R+\sigma)^2,x)\:dx\\
    &\le\frac{d}{d\sigma}\int_{-(R+\sigma)^2}^0
     \left(\int_{\p B_{R+\sigma}}u(s,\cdot)\right)ds
     +\int_{B_{R+\sigma}}u(-(R+\sigma)^2,x)\:dx.
\end{split}
\end{equation*}
Here the first step uses the inclusions
$P_R\subset P_{R+\sigma}\subset P_{R+r}$.
The third step follows from~(\ref{eq:ID-I}) and 
the last from~(\ref{eq:ID-II}).  Now integrate this
inequality over the interval $0\le\sigma\le t$, 
with $r/2\le t\le r$, to obtain
\begin{equation*}
\begin{split}
    &\frac{r}{2}\left(\int_{P_R}g-\int_{P_{R+r}}f\right) \\
    &\le\int_{-(R+t)^2}^0\left(
     \int_{\p B_{R+t}}u(s,\cdot)\right)ds
     +\int_0^r\left(
     \int_{B_{R+\sigma}}u(-(R+\sigma)^2,\cdot)
     \right)d\sigma \\
%
%
    &\le\int_{-(R+r)^2}^0\left(
     \int_{\p B_{R+t}}u(s,\cdot)\right)ds
     +\frac{1}{2R}\int_{-(R+r)^2}^{-R^2}
     \int_{B_{R+r}}u(s,x)\:dxds.
\end{split}
\end{equation*}
Here the last step follows by substituting $s=-(R+\sigma)^2$
and using $\sqrt{-s}\ge R$. 
Integrate this inequality again over 
the interval $r/2\le t\le r$ to obtain
$$
     \frac{r}{2}\left(\int_{P_R}g-\int_{P_{R+r}}f\right)
     \le\frac{2}{r}\int_{P_{R+r}\setminus P_R}u
     +\frac{1}{2R}\int_{P_{R+r}\setminus P_R}u.
$$
This proves Lemma~\ref{le:apriori-L2}.
\end{proof}

\begin{lemma}\label{le:apriori-L2-eps}
Let $\eps,R,r$ be positive real numbers.
Let $u:\R^2\supset P_{R+r}^\eps\to\R$ be
a $C^2$ function and $f,g:P_{R+r}^\eps\to\R$
be continous functions such that
$$
     \left(\eps^2{\p_s}^2+{\p_t}^2-\p_s\right)u
     \ge g-f, \qquad
     u\ge0, \qquad f\ge0, \qquad g\ge0.
$$
Then
$$
     \int_{P_{R/2}^{-\eps}}g
     \le \frac{2(R+r)}{R}\int_{P_{R+r}^\eps}f
     +\frac{2(R+r)}{R}
     \left(\frac{4}{r^2}+\frac{1}{Rr}\right)
     \int_{P_{R+r}^\eps}u.
$$
\end{lemma}

\begin{proof}
The idea of proof is as in
Lemma~\ref{le:apriori-basic-eps}.
Increase the dimension of the domain
from two to three and
apply Lemma~\ref{le:apriori-L2} with $n=2$.
Define functions $U,F,G$ on
$P_{R+r} \subset \R \times \R^2$ by
$$
     U(s,t,q)
     :=u(s+\eps q,t),\quad
     F(s,t,q)
     :=f(s+\eps q,t),\quad
     G(s,t,q)
     :=g(s+\eps q,t).
$$
The new variable $\sigma:=s+\eps q$ satisfies
$(\sigma,t)\in P_{R+r}^\eps\subset\R\times\R$
whenever $(s,t,q)\in P_{R+r}\subset\R\times\R^2$.
Use the differential inequality
in the assumption of the lemma to conclude
$
     \left(\Delta-\p_s\right)U
     \ge G-F,
$
where $\Delta:=\p_t^2+\p_q^2$. 
Thus Lemma~\ref{le:apriori-L2} with $n=2$
yields
\begin{equation*}
\begin{split}
     \int_{P_R}G
    &\le\int_{P_{R+r}}F+\left(\frac{4}{r^2}
     +\frac{1}{Rr}\right)\int_{P_{R+r}}U \\
    &\le2(R+r)\int_{P_{R+r}^\eps}f
     +2(R+r)\left(\frac{4}{r^2}
     +\frac{1}{Rr}\right)\int_{P_{R+r}^\eps}u.
\end{split}
\end{equation*}
The last step uses (\ref{eq:helper-1}).
By definition of $G$
\begin{equation*}
\begin{split}
     \int_{P_R}G
    &=\int_{B_R\subset\R^2}\left(
     \int_{-R^2}^0g(s+\eps q,t)\:ds\right)dqdt \\
    &\ge\int_{-R/2}^{R/2}\int_{-R/2}^{R/2}
     \int_{-R^2+\eps q}^{\eps q}
     g(\sigma,t)\:d\sigma dqdt \\
    &\ge\int_{-R/2}^{R/2}\int_{-R/2}^{R/2}
     \int_{-R^2+\eps R/2}^{-\eps R/2}
     g(\sigma,t)\:d\sigma dqdt \\
    &=R\int_{-R/2}^{R/2}\int_{-R^2+\eps R/2}^{-\eps R/2}
     g(\sigma,t)\:d\sigma dt \\
    &\ge R\int_{P_{R/2}^{-\eps}}g.
\end{split}
\end{equation*}
This proves the lemma.
\end{proof}

\begin{lemma}\label{le:apriori-1-L2}
Fix three positive constants $r,R,\eps$ and three functions
$u,f,g:[-(R+r)^2-\eps(R+r),\eps(R+r)]\to\R$ such that $u$ 
is $C^2$ and $f,g$ are continuous. If
$$
\eps^2u''-u'\ge g-f,\qquad u\ge 0,\qquad
f\ge 0,\qquad g\ge 0,
$$
then
\begin{align*}
\int_{-R^2/4+R\eps/2}^{-R\eps/2}g(s)\,ds 
&\le \frac{2(R+r)}{R}\int_{-(R+r)^2-\eps(R+r)}^{\eps(R+r)}f(s)\,ds\\
&\quad
+\frac{2(R+r)}{R}\left(\frac{4}{r^2}+\frac{1}{Rr}\right)
\int_{-(R+r)^2-\eps(R+r)}^{\eps(R+r)}u(s)\,ds.
\end{align*}
\end{lemma}

\begin{proof}
This follows immediately from Lemma~\ref{le:apriori-L2-eps}
with $u$, $f$, and $g$ independent of the $t$-variable.
\end{proof}


\section{Two fundamental $L^p$ estimates}\label{app:Lp}

\begin{theorem}\label{thm:CZ}
For every $p>1$ there is a constant $c=c(p)>0$ such that
\begin{equation}\label{eq:CZ}
     \left\|\p_su\right\|_{L^p} + \left\|\p_sv\right\|_{L^p}
     \le c\left(\left\|\p_su-\p_tv\right\|_{L^p}
         + \left\|\p_sv+\p_tu-v\right\|_{L^p}\right)
\end{equation}
for all $u,v\in\Cinf_0(\R^2)$.
\end{theorem}

\begin{theorem}\label{thm:parabolic}
For every $p>1$ there is a constant $c=c(p)>0$ such that
\begin{equation}\label{eq:parabolic}
     \left\|\p_su\right\|_{L^p} + \left\|\p_t\p_tu\right\|_{L^p}
     \le c\left\|\p_su-\p_t\p_tu\right\|_{L^p}
\end{equation}
for every $u\in\Cinf_0(\R^2)$.
\end{theorem}

If we assume $v=\p_tu$ then~(\ref{eq:CZ}) follows
from~(\ref{eq:parabolic}) (but not conversely).  
On the other hand if the term
$\p_sv+\p_tu-v$ on the right is replaced by 
$\p_sv+\p_tu$, then~(\ref{eq:CZ}) becomes the 
Calderon-Zygmund inequality.  However, it seems that
the estimate~(\ref{eq:CZ}) in its full strength
cannot be deduced directly from the Calderon--Zygmund
inequality and the parabolic estimate~(\ref{eq:parabolic}). 
Theorems~\ref{thm:CZ} and~\ref{thm:parabolic}
will be proved below. 

\begin{corollary}\label{cor:CZ-eps}
Let $p>1$ and denote by $c=c(p)$ the constant 
of Theorem~\ref{thm:CZ}.  Then
\begin{equation}\label{eq:CZ-eps}
     \left\|\p_su\right\|_{L^p} + \eps\left\|\p_sv\right\|_{L^p}
     \le c\left(\left\|\p_su-\p_tv\right\|_{L^p}
         + \eps\left\|\p_sv+\eps^{-2}(\p_tu-v)\right\|_{L^p}\right)
\end{equation}
for every $\eps>0$ and every pair $u,v\in\Cinf_0(\R^2)$.
\end{corollary}

\begin{proof}  Denote 
$$
     f:=\p_su-\p_tv,\qquad g := \p_sv+\eps^{-2}(\p_tu-v).
$$
Now consider the rescaled functions
$$
     \tilde{u}(s,t) := u(\eps^2s,\eps t),\qquad
     \tilde{v}(s,t) := \eps v(\eps^2 s,\eps t)
$$
and
$$
     \tilde{f}(s,t) := \eps^2f(\eps^2 s,\eps t),\qquad
     \tilde{g}(s,t) := \eps^3g(\eps^2 s,\eps t).
$$
Then 
$$
     \p_s\tilde{u}-\p_t\tilde{v} = \tilde{f},\qquad
     \p_s\tilde{v}+\p_t\tilde{u}-\tilde{v} = \tilde{g}.
$$
Hence, by Theorem~\ref{thm:CZ},
$$
     \bigl\|\p_s\tilde{u}\bigr\|_{L^p}
     +\bigl\|\p_s\tilde{v}\bigr\|_{L^p}
     \le c\left(\bigl\|\tilde{f}\bigr\|_{L^p}
     +\bigl\|\tilde{g}\bigr\|_{L^p}\right).
$$
Now the result follows from the fact that 
$$
     \bigl\|\p_s\tilde{u}\bigr\|_{L^p}
     =\eps^{2-3/p}\bigl\|\p_su\bigr\|_{L^p},\qquad
     \bigl\|\tilde{f}\bigr\|_{L^p}
     =\eps^{2-3/p}\bigl\|f\bigr\|_{L^p},
$$
and similarly for the other terms.
\end{proof}

We give a proof of~(\ref{eq:CZ}) and~(\ref{eq:parabolic})
that is based on the Marcinkiewicz--Mihlin multiplier method. 
To formulate the result, we consider the Fourier transform
$$
      \Ff:L^2(\R^2,\C)\to L^2(\R^2,\C),
$$
given by 
$$
      (\Ff f)(\sigma,\tau)
      := \frac{1}{2\pi}\int_{-\infty}^\infty\int_{-\infty}^\infty
         e^{-i(\sigma s+\tau t)}f(s,t)\,dsdt
$$
for $f\in L^2(\R^2,\C)\cap L^1(\R^2,\C)$.  Given a bounded 
measurable complex valued function $m:\R^2\to\C$ define the 
bounded linear operator 
$$
      \Tt_m:L^2(\R^2,\C)\to L^2(\R^2,\C)
$$
by 
$$
      \Tt_mf := \Ff^{-1}(m\Ff f).
$$
The following theorem is proved in~\cite{LSU}.

\begin{theorem}[Marcinkiewicz--Mihlin]\label{thm:LSU}
For every $c>0$ and every $p>1$ there is 
a constant $c_p=c_p(c)>0$ such that the following holds.
If $m:\R^2\to\C$ is a measurable function
such that the restriction of $m$ to each of the four 
open quadrants in $\R^2$ is twice continuously differentiable
and
\begin{equation}\label{eq:LSU}
     |m(\sigma,\tau)| + \left|\sigma\p_\sigma m(\sigma,\tau)\right|
     + \left|\tau\p_\tau m(\sigma,\tau)\right|
     + \left|\sigma\tau\p_\sigma\p_\tau m(\sigma,\tau)\right|
     \le c
\end{equation}
for $\sigma,\tau\in\R\setminus\{0\}$ then
$$
     f\in L^p(\R^2,\C)\cap L^2(\R^2,\C)\qquad\IMP\qquad
     \Tt_mf\in L^p(\R^2,\C)
$$
and 
$$
     \left\|\Tt_mf\right\|_{L^p}\le c_p\left\|f\right\|_{L^p}
$$
for every $f\in L^p(\R^2,\C)\cap L^2(\R^2,\C)$. 
\end{theorem}

\begin{remark}\label{rmk:mm}\rm
The theorem of Marcinkiewicz--Mihlin in its original form
is slightly stronger than Theorem~\ref{thm:LSU}, namely
condition~(\ref{eq:LSU}) is replaced by the weaker conditions
\begin{equation}\label{eq:mm0}
     \sup_{\sigma,\tau}|m(\sigma,\tau)| \le c,
\end{equation}
\begin{equation}\label{eq:mm1}
    \sup_{\sigma\ne0}\int_{2^\ell}^{2^{\ell+1}}\left|\p_\tau 
    m(\sigma,\pm\tau)\right|\,d\tau\le c,\qquad
    \sup_{\tau\ne0}\int_{2^k}^{2^{k+1}}\left|\p_\sigma
    m(\pm\sigma,\tau)\right|\,d\sigma\le c
\end{equation}
and
\begin{equation}\label{eq:mm2}
    \int_{2^k}^{2^{k+1}} \int_{2^\ell}^{2^{\ell+1}}
    \left|\p_\sigma\p_\tau m(\pm\sigma,\pm\tau)\right|\,d\tau\le c
\end{equation}
for all integers $k$ and $\ell$ (and all choices of signs).
In this form the result is proved in Stein~\cite[Theorem~6']{STEIN}.
It is easy to see that~(\ref{eq:LSU}) implies
(\ref{eq:mm1}) with $c$ replaced by $c\log 2$ 
and~(\ref{eq:mm2}) with $c$ replaced by $c(\log 2)^2$.
\end{remark}

\begin{proof}[Proof of Theorem~\ref{thm:parabolic}.]
Let $u\in\Cinf_0(\R^2,\C)$ and define 
$f\in\Cinf_0(\R^2)$ by 
$$
     f:=\p_su-\p_t\p_tu.
$$
Denote the Fourier transforms of $f$ and $u$ by 
$$
     \Hat f:=\Ff f,\qquad \Hat u:=\Ff u.
$$
Then 
$$
     \Hat f = i\sigma\Hat u + \tau^2\Hat u
$$
and hence 
$$
     \Hat{\p_su} = i\sigma\Hat u 
     = \frac{i\sigma}{\tau^2+i\sigma}\Hat{f}.
$$
Denote the multiplier in this equation by
$$
     m(\sigma,\tau):=\frac{i\sigma}{\tau^2+i\sigma}.
$$ 
The formulae
$$
     \p_\sigma m = 
     \frac{i\tau^2}{\left(\tau^2+i\sigma\right)^2},\qquad
     \p_\tau m = 
     \frac{-2i\sigma\tau}{\left(\tau^2+i\sigma\right)^2},\qquad
     \p_\sigma\p_\tau m = 
     \frac{-2i\tau(\tau^2-i\sigma)}{\left(\tau^2+i\sigma\right)^3}
$$
show that the functions $m$, $\sigma\p_\sigma m$,
$\tau\p_\tau m$, and $\sigma\tau\p_\sigma\p_\tau m$ 
are bounded.  Hence the result follows from 
Theorem~\ref{thm:LSU}. 
\end{proof}

\begin{proof}[Proof of Theorem~\ref{thm:CZ}.]
Let $u,v\in\Cinf_0(\R^2,\C)$ and define 
$f,g\in\Cinf_0(\R^2)$ by 
$$
     f:=\p_su-\p_tv,\qquad
     g:=\p_sv+\p_tu-v.
$$
Then 
$$
     \Hat f = i\sigma\Hat u - i\tau\Hat v,\qquad
     \Hat g = i\sigma\Hat v + i\tau\Hat u - \Hat v.
$$
Solving this equation for $\Hat u$ and $\Hat v$ we find
$$
     \Hat u = \frac{1-i\sigma}{\sigma^2+\tau^2+i\sigma}\Hat f
              -  \frac{i\tau}{\sigma^2+\tau^2+i\sigma}\Hat g,
$$
$$
     \Hat v = \frac{i\tau}{\sigma^2+\tau^2+i\sigma}\Hat f
              -  \frac{i\sigma}{\sigma^2+\tau^2+i\sigma}\Hat g,
$$
and hence 
$$
     \Hat{\p_su} 
     = \frac{\sigma^2+i\sigma}{\sigma^2+\tau^2+i\sigma}\Hat f
       + \frac{\sigma\tau}{\sigma^2+\tau^2+i\sigma}\Hat g,
$$
$$
     \Hat{\p_sv} 
     = \frac{-\sigma\tau}{\sigma^2+\tau^2+i\sigma}\Hat f
       + \frac{\sigma^2}{\sigma^2+\tau^2+i\sigma}\Hat g.
$$
The four multipliers in the last two equations 
satisfy~(\ref{eq:LSU}).  Hence the result follows from 
Theorem~\ref{thm:LSU}. 
\end{proof}


\section{The estimate for the inverse}\label{app:inverse}

We begin by proving a weaker version of the estimate in 
Theorem~\ref{thm:elliptic-eps}.

\begin{proposition}\label{prop:elliptic-linear}
Let $u\in\Cinf(\R \times S^1,M)$ 
and $v\in\Omega^0(\R\times S^1,u^*TM)$ such that
$\|\p_s u\|_\infty$, $\|\p_t u\|_\infty$ and $\|v\|_\infty$
are finite and
$
     \lim_{s\to\pm\infty}u(s,t)
$
exists, uniformly in $t$. Then, for every $p>1$,
there is a constant $c>0$ such that
\begin{equation}\label{eq:elliptic-linear}
\begin{split}
    &\eps^{-1}\|\Nabla{t}\xi-\eta\|_p+\|\Nabla{t}\eta\|_p
     +\|\Nabla{s} \xi \|_p+\eps\|\Nabla{s}\eta\|_p \\
    &\le c\left( \| \Dd_{u,v}^\eps\zeta \|_{0,p,\eps} 
     +\eps^{-1}\|\zeta\|_{0,p,\eps}\right)
\end{split}
\end{equation}
for every $\eps\in(0,1]$ and every pair of 
compactly supported vector fields
$\zeta=(\xi,\eta)\in\Omega^0(\R\times S^1,u^*TM\oplus u^*TM)$.
The formal adjoint operator $(\Dd_{u,v}^\eps)^*$ satisfies
the same estimate.
\end{proposition}

\begin{proof}
Choose a finite open cover $\{U_\alpha\}_\alpha$ 
of the cylinder $\R \times S^1$
with the following properties.
\begin{enumerate}
\item[(i)]
For each $\alpha$ the set
$U_\alpha \subset \R \times S^1$
is contractible.
\item[(ii)]
For each $\alpha$ the closure of the 
image of $U_\alpha$ under $u$
is contained in a coordinate chart on $M$.
\item[(iii)]
There is a constant $T>0$
and an open cover $\{I_\alpha\}_\alpha$
of $S^1$ such that
$U_\alpha \cap [T,\infty) \times S^1
=[T,\infty) \times I_\alpha$
for every $\alpha$.
Similarly for the interval $(-\infty,-T]$.
\end{enumerate}
We prove~(\ref{eq:elliptic-linear}) for $\Dd^\eps_{u,v}$.
The estimate for $(\Dd^\eps_{u,v})^*$ is analoguous.
Assume first that $\xi$ and $\eta$ are compactly supported
in $U_\alpha$ for some $\alpha$ and denote by 
$\xi_\alpha,\eta_\alpha:U_\alpha\to\R^n$ the vector 
fields in local coordinates. By Corollary~\ref{cor:CZ-eps},
there is a constant $c_\alpha$, depending only on $p$ and 
the metric, such that 
$$
     \left\|\p_s\xi_\alpha\right\|_p
     +\eps\left\|\p_s\eta_\alpha\right\|_p 
     \le c_\alpha\Bigl(\left\|\p_s\xi_\alpha-\p_t\eta_\alpha\right\|_p
     +\eps\left\|\p_s\eta_\alpha
     +\eps^{-2}(\p_t\xi_\alpha-\eta_\alpha)\right\|_p
     \Bigr).
$$
Here we denote by $\left\|\cdot\right\|_p$
the $L^p$ norm
with respect to the Riemannian metric in
the coordinate charts
on $M$.  Replacing the partial derivatives
$\p_s$ and $\p_t$
by the covariant derivatives $\Nabla{s}$
and $\Nabla{t}$ we obtain 
\begin{equation}\label{eq:p1}
\begin{split}
     \left\|\Nabla{s}\xi\right\|_p
     +\eps\left\|\Nabla{s}\eta\right\|_p 
    &\le c\Bigl(\left\|\Nabla{s}\xi-\Nabla{t}\eta\right\|_p
     +\eps\left\|\Nabla{s}\eta
     +\eps^{-2}(\Nabla{t}\xi-\eta)\right\|_p \\
    &\qquad +\eps^{-1}\left\|\xi\right\|_p+\left\|\eta\right\|_p
     \Bigr).
\end{split}
\end{equation}
for every $\xi$ with support in one
of the sets $U_\alpha$. 
Here we have used the $L^\infty$ bounds
on $\p_su$ and $\p_tu$.
Observe that the constant $c$ depends on the Christoffel
symbols determined by our coordinate chart on $M$. 
Now let $\{\beta_\alpha\}_\alpha$ be a partition
of unity subordinate to the cover $\{U_\alpha\}_\alpha$
such that $\left\|\p_s\beta_\alpha\right\|_\infty
+\left\|\p_t\beta_\alpha\right\|_\infty<\infty$ for 
every~$\alpha$. (Note that $\beta_\alpha$ need not have 
compact support when $U_\alpha$ is unbounded.)
Given any two compactly supported vector fields
$\xi,\eta\in\Om^0(\R\times S^1,u^*TM)$
apply~(\ref{eq:p1})
to the (compactly supported) pair 
$(\beta_\alpha\xi,\beta_\alpha\eta)$ and take
the sum to deduce that~(\ref{eq:p1})
continues to hold for 
the pair $(\xi,\eta)$ with an
appropriate larger constant $c$. 
Using the $L^\infty$ bounds on
$\p_su$, $\p_tu$, $v$, and the curvature 
(as well as the axioms~$(V0-V1)$ for $\Vv$)
we obtain
\begin{equation*}
\begin{split}
     \left\|\Nabla{s}\xi\right\|_p
     +\eps\left\|\Nabla{s}\eta\right\|_p 
     &\le c'\Bigl(\left\|\Nabla{s}\xi-\Nabla{t}\eta
      -R(\xi,\p_tu)v-\Hh_\Vv(u)\xi\right\|_p \\
     &\qquad+\eps\left\|\Nabla{s}\eta+R(\xi,\p_su)v
     +\eps^{-2}(\Nabla{t}\xi-\eta)\right\|_p  \\
     &\qquad+\eps^{-1}\left\|\xi\right\|_p+\left\|\eta\right\|_p
     \Bigr).
\end{split}
\end{equation*}
This implies~(\ref{eq:elliptic-linear}).
\end{proof}

Under the assumptions of Proposition~\ref{prop:elliptic-linear}
it follows immediately that
\begin{equation}\label{eq:elliptic-standard}
     \|\zeta\|_{1,p,\eps} 
     \le c\left( \eps^2 \| \Dd_{u,v}^\eps\zeta \|_{0,p,\eps} 
     +\|\zeta\|_{0,p,\eps}\right)
\end{equation}
and similarly for $(\Dd_{u,v}^\eps)^*$.  
Moreover, note that the difference between 
Proposition~\ref{prop:elliptic-linear} and 
Theorem~\ref{thm:elliptic-eps} lies in the $\eps$-factors
in front of $\left\|\xi\right\|_p$ and $\left\|\eta\right\|_p$
on the right hand sides of the estimates. 
To prove Theorem~\ref{thm:elliptic-eps} we must improve 
these these factors by $\eps$ for $\xi$ and by $\eps^2$ 
for $\eta$.  This requires the following parabolic estimate.
Let $1/p+1/q=1$.
The \emph{formal adjoint operator} 
$$
     (\Dd^0_u)^*:\Ww_u^q\to\Ll_u^q
$$
of $\Dd_u^0:\Ww_u^p\to\Ll_u^p$
is given by
\begin{equation}\label{eq:D0^*}
     (\Dd^0_u)^*\xi 
     = -\Nabla{s}\xi
       - \Nabla{t}\Nabla{t}\xi 
       - R(\xi,\p_tu)\p_tu
       - \Hh_\Vv(u)\xi.
\end{equation}

\begin{proposition}\label{prop:parabolic-linear}
Let $u\in\Cinf(\R \times S^1,M)$ such that
$\|\p_s u\|_\infty$, $\|\p_t u\|_\infty$
and $\|\Nabla{t}\p_t u\|_\infty$ are finite and
$
     \lim_{s\to\pm\infty}u(s,t)
$
exists, uniformly in $t$. Then, for every $p>1$,
there is a constant $c>0$ such that
\begin{equation}\label{eq:parabolic-linear}
     \| \Nabla{s} \xi \|_p + \| \Nabla{t} \Nabla{t}  \xi \|_p
     \le c \left( \| \Dd_u^0 \xi \|_p + \| \xi \|_p \right)
\end{equation}
for every compactly supported vector field
$\xi \in \Omega^0(\R \times S^1, u^*TM)$.
The formal adjoint operator $(\Dd_u^0)^*$ satisfies
the same estimate.
\end{proposition}

\begin{lemma}\label{le:eat-eps}
Let $x:S^1\to M$ be a smooth map, $p>1$ and
\begin{equation}\label{eq:kappa-p}
     \kappa_p:=
     \begin{cases}
       p
      &\mbox{if }p\ge2, \\
       p/(p-1)
      &\mbox{if }p\le2.
     \end{cases} 
\end{equation}
Then, for every $\eps>0$ and every $\xi\in\Om^0(S^1,x^*TM)$,
we have
\begin{align*}
     \|(\1-\eps\Nabla{t}\Nabla{t})^{-1}\xi\|_p
    &\le \|\xi\|_p,\\
     \sqrt{\eps}\|(\1-\eps\Nabla{t}\Nabla{t})^{-1}
     \Nabla{t}\xi\|_p
    &\le \kappa_p \|\xi\|_p,\\
     \eps\|(\1-\eps\Nabla{t}\Nabla{t})^{-1}
     \Nabla{t}\Nabla{t}\xi\|_p
    &\le 2 \|\xi\|_p.
\end{align*}
These estimates continue to hold for
$u\in\Cinf(\R \times S^1,M)$ and compactly supported
vector fields $\xi\in\Omega^0(\R \times S^1,u^*TM)$.
\end{lemma}

\begin{proof}
First consider the case $p\ge2$:
Let $\eps>0$ and $\xi\in\Om^0(S^1,x^*TM)$.  Define
$$
     \eta := (\1 -\eps\Nabla{t}\Nabla{t})^{-1} \xi.
$$
(The operator
$
      (\1 -\eps\Nabla{t}\Nabla{t}): 
      W^{2,p}(S^1,x^*TM)\to L^p(S^1,x^*TM)
$
is bijective.)  Then
\begin{equation*}
\begin{split}
     \frac{d^2}{dt^2} \abs{\eta}^p
    &=\frac{d}{dt} \Bigl( p
     \left|\eta\right|^{p-2}\inner{\Nabla{t}\eta}{\eta} \Bigr) \\
    &= p(p-2) \abs{\eta}^{p-4} 
     \langle\Nabla{t} \eta,\eta\rangle^2
     +p \abs{\eta}^{p-2} \Bigl(
     \langle\Nabla{t}\Nabla{t}\eta,\eta\rangle
     +\abs{\Nabla{t}\eta}^2 \Bigr) \\
    &\ge p\eps^{-1} \abs{\eta}^p -p\eps^{-1} \abs{\eta}^{p-2}
      \langle\xi,\eta\rangle \\
    &\ge p\eps^{-1}\abs{\eta}^p 
     -p\eps^{-1}\abs{\eta}^{p-1}\abs{\xi} \\
    &\ge \eps^{-1}\abs{\eta}^p -\eps^{-1}\abs{\xi}^p.
\end{split}
\end{equation*}
The third step uses the identity 
$\Nabla{t}\Nabla{t}\eta=\eps^{-1}\eta-\eps^{-1}\xi$.
The last step uses Young's inequality
\begin{equation} \label{eq:Young}
     ab \le \frac{a^r}{r} + \frac{b^s}{s},\qquad
     \frac{1}{r}+\frac{1}{s}=1,
\end{equation}
with $r=p$, $a=\abs{\xi}$ and $s=p/(p-1)$, $b=\abs{\eta}^{p-1}$.
Moreover,
\begin{equation*}
\begin{split}
     &\frac{d}{d t}
      \Bigl( 
      \left|\Nabla{t}\eta\right|^{p-2}\inner{\Nabla{t}\eta}{\eta} \Bigr) \\
     &=\abs{\Nabla{t} \eta}^p
      +\abs{\Nabla{t} \eta}^{p-2}\inner{\Nabla{t}\Nabla{t}\eta}{\eta} 
      +(p-2)\abs{\Nabla{t} \eta}^{p-4} \inner{\Nabla{t}\eta}{\eta}
      \langle\Nabla{t}\Nabla{t}\eta,\Nabla{t}\eta\rangle \\
     &= \abs{\Nabla{t} \eta}^p
      +\eps^{-1}\abs{\Nabla{t} \eta}^{p-2}\abs{\eta}^2 
      -\eps^{-1}\abs{\Nabla{t} \eta}^{p-2}\inner{\xi}{\eta} \\
     &\quad 
      -\eps^{-1}(p-2) \abs{\Nabla{t} \eta}^{p-4}
      \langle \Nabla{t} \eta,\eta \rangle
      \langle \xi ,\Nabla{t} \eta \rangle 
      +\eps^{-1}(p-2) \abs{\Nabla{t} \eta}^{p-4} 
      \langle \Nabla{t} \eta,\eta  \rangle^2 \\
     &\ge \abs{\Nabla{t} \eta}^p
      +\tfrac{1}{2}\eps^{-1}\abs{\Nabla{t} \eta}^{p-2}\abs{\eta}^2 
      - \tfrac{p-1}{2}\eps^{-1}\abs{\Nabla{t} \eta}^{p-2}\abs{\xi}^2
      +\tfrac{p-2}{2} \eps^{-1}\abs{\Nabla{t} \eta}^{p-4} 
        \langle \Nabla{t} \eta , \eta \rangle^2 \\
     &\ge \abs{\Nabla{t} \eta}^p
      - \tfrac{p-1}{2}\eps^{-1}\abs{\Nabla{t} \eta}^{p-2}\abs{\xi}^2 \\
     &\ge \tfrac{2}{p}
      \abs{\Nabla{t} \eta}^p 
      - \tfrac{2}{p} \bigl( \tfrac{p-1}{2} \bigr)^{p/2} 
      \eps^{-p/2}\abs{\xi}^p.
\end{split}
\end{equation*}
The third step uses~(\ref{eq:Young}) with $r=s=2$.
The last step uses~(\ref{eq:Young}) with
$r=p/2$, $a=\tfrac{p-1}{2}\eps^{-1}\abs{\xi}^2$
and $s=p/(p-2)$, $b=\abs{\Nabla{t}\eta}^{p-2}$.
Now the first two estimates of the lemma follow
by integration over $S^1$, respectively $\R\times S^1$.
The last estimate is an easy consequence of the first:
$$
     \eps \| \Nabla{t}\Nabla{t} \eta \|_p
     = \| \eta - \xi \|_p
     \le \| \eta \|_p + \| \xi \|_p
     \le 2 \| \xi \|_p.
$$
This proves the lemma for $p\ge2$.
Now assume $1<p<2$ and let $q:=p/(p-1)$.
Then $q>2$ and hence
\begin{equation*}
\begin{split}
     \sqrt{\eps}\Norm{
     (\1-\eps\Nabla{t}\Nabla{t})^{-1}\Nabla{t}\xi}_p
    &=\sqrt{\eps}\sup_{0\not=\eta\in L^q}\frac{
     \left\langle(\1-\eps\Nabla{t}\Nabla{t})^{-1}
     \Nabla{t}\xi,\eta\right\rangle}
     {\Norm{\eta}_q} \\
    &\le\sqrt{\eps}\sup_{0\not=\eta\in L^q}\frac{
     \bigl\|\xi\bigr\|_p
     \Norm{(\1-\eps\Nabla{t}\Nabla{t})^{-1}
     \Nabla{t}\eta}_q}{\Norm{\eta}_q} \\
    &\le q\Norm{\xi}_p.
\end{split}
\end{equation*}
This prove the second estimate for $p<2$.
The other estimates follow similarly.
This proves the lemma.
\end{proof} 

\begin{lemma}\label{le:nabla-t-xi}
Let $x\in\Cinf(S^1,M)$ and $p>1$.  Then
$$
     \|\Nabla{t}\xi\|_p
     \le \kappa_p \left( \delta^{-1} \| \xi \|_p
     + \delta \| \Nabla{t} \Nabla{t} \xi \|_p\right)     
$$
for $\delta>0$ and $\xi\in\Omega^0(S^1,x^*TM)$,
where $\kappa_p$ is defined by~(\ref{eq:kappa-p}).
This estimate continues to hold for
$u\in\Cinf(\R \times S^1,M)$ and compactly supported
vector fields $\xi\in\Omega^0(\R \times S^1,u^*TM)$.
\end{lemma}

\begin{proof}
Let $1/p+1/q=1$. Since the operator 
$$
     W^{2,q}(S^1,x^*TM)\to L^q(S^1,x^*TM):
     \eta\mapsto\delta^{-1}\eta+\delta\Nabla{t}\Nabla{t}\eta
$$
is bijective, we have
\begin{equation*}
\begin{split}
     \Norm{\Nabla{t}\xi}_p
    &=\sup_{\eta\in W^{2,q}} 
     \frac{\inner{\Nabla{t}\xi}
     {\delta^{-1}\eta-\delta\Nabla{t}\Nabla{t}\eta}}
     {\Norm{\delta^{-1}\eta-\delta\Nabla{t}\Nabla{t}\eta}_q} \\
    &=\sup_{\eta\in W^{2,q}} 
     \frac{-\inner{\xi}{\delta^{-1}\Nabla{t}\eta}
     +\inner{\Nabla{t}\Nabla{t}\xi}{\delta\Nabla{t}\eta}}
     {\Norm{\delta^{-1}\eta-\delta\Nabla{t}\Nabla{t}\eta}_q} \\
    &\le 
     \left(\delta^{-1}\Norm{\xi}_p
     +\delta\Norm{\Nabla{t}\Nabla{t}\xi}_p\right)
     \sup_{\eta\in W^{2,q}} 
     \frac{\Norm{\Nabla{t}\eta}_q}
     {\Norm{\delta^{-1}\eta-\delta\Nabla{t}\Nabla{t}\eta}_q} \\
    &\le \kappa_p \left(\delta^{-1}\Norm{\xi}_p
     +\delta\Norm{\Nabla{t}\Nabla{t}\xi}_p\right).
\end{split}
\end{equation*}
To prove the last step, denote
$$
     \zeta:=\eta-\delta^2\Nabla{t}\Nabla{t}\eta.
$$
Then
$$
     \Nabla{t}\eta
     =\left(\1-\delta^2\Nabla{t}\Nabla{t}\right)^{-1}\Nabla{t}\zeta
$$
and hence, by Lemma~\ref{le:eat-eps} with $\eps=\delta^2$,
we have
$$
     \Norm{\Nabla{t}\eta}_q
     \le \kappa_q\delta^{-1}\Norm{\zeta}_q
     =\kappa_p\Norm{\delta^{-1}\eta-\delta\Nabla{t}\Nabla{t}\eta}_q.
$$
We have used the fact that $\kappa_p=\kappa_q$. 
This proves the lemma.
\end{proof}

\begin{proof}[Proof of Proposition~\ref{prop:parabolic-linear}.]
The proof follows the same pattern as that of
Proposition~\ref{prop:elliptic-linear}.
Let $\{U_\alpha\}_\alpha$ be as above. 
If $\xi$ is (compactly) supported in $U_\alpha$ then,
by Theorem~\ref{thm:parabolic}, 
$$
     \|\p_s\xi_\alpha\|_p + \|\p_t\p_t\xi_\alpha\|_p
     \le c_\alpha\left\|\p_s\xi_\alpha-\p_t\p_t\xi_\alpha\right\|_p
$$
Replacing $\p_s$ and $\p_t$ by $\Nabla{s}$ and $\Nabla{t}$,
and using the $L^\infty$ bounds on $\p_su$, $\p_tu$, and 
$\Nabla{t}\p_tu$, we find 
$$
     \|\Nabla{s}\xi\|_p + \|\Nabla{t}\Nabla{t}\xi\|_p
     \le c\Bigl(\left\|\Nabla{s}\xi-\Nabla{t}\Nabla{t}\xi\right\|_p
          + \left\|\xi\right\|_p+\left\|\Nabla{t}\xi\right\|_p\Bigr).
$$
Using a partition of unity $\{\beta_\alpha\}_\alpha$,
subordinate to the cover $\{U_\alpha\}_\alpha$, such that 
$$
     \left\|\p_s\beta_\alpha\right\|_\infty
     +\left\|\p_t\beta_\alpha\right\|_\infty
     + \left\|\p_t\p_t\beta_\alpha\right\|_\infty<\infty,
$$
we deduce that the last estimate continues to hold for
every compactly supported vector field 
$\xi \in \Omega^0(\R \times S^1, u^*TM)$.
Now apply Lemma~\ref{le:nabla-t-xi} with 
$\delta cp<1/2$ to obtain 
$$
     \|\Nabla{s}\xi\|_p + \|\Nabla{t}\Nabla{t}\xi\|_p
     \le c'\Bigl(\left\|\Nabla{s}\xi-\Nabla{t}\Nabla{t}\xi\right\|_p
          + \left\|\xi\right\|_p\Bigr).
$$
Hence 
$$
     \|\Nabla{s}\xi\|_p + \|\Nabla{t}\Nabla{t}\xi\|_p
     \le c''\Bigl(\left\|\Nabla{s}\xi-\Nabla{t}\Nabla{t}\xi
         -R(\xi,\p_tu)\p_tu-\Hh_\Vv(u)\xi\right\|_p
          + \left\|\xi\right\|_p\Bigr)
$$
as required.
\end{proof}

\begin{proof}[Proof of Theorem~\ref{thm:elliptic-eps}]
Fix a constant $p>1$ and define
$$
  f(\xi,\eta):=\Nabla{s}\xi-\Nabla{t}\eta, \qquad
  g(\xi,\eta):=\Nabla{s}\eta+\eps^{-2}(\Nabla{t}\xi-\eta),
$$
for compactly supported vector fields
$\zeta=(\xi,\eta)\in\Omega^0(\R\times S^1,u^*TM\oplus u^*TM)$.
It suffices to show that
\begin{equation}\label{eq:main-elliptic}
\begin{split}
     &\eps^{-1}\left\|\Nabla{t}\xi-\eta\right\|_p
      + \left\|\Nabla{t}\eta\right\|_p
      + \left\|\Nabla{s}\xi\right\|_p
      + \eps\left\|\Nabla{s}\eta\right\|_p \\
     &\le 
      c\left(
      \left\|f\right\|_p+\eps\left\|g\right\|_p
      + \left\|\xi\right\|_p
      + \eps^2\left\|\eta\right\|_p
      \right)
\end{split}
\end{equation}
for some constant $c>0$ independent of $\eps$ and $(\xi,\eta)$.
The general case (for $\Dd^\eps_{u,v}$) then follows 
easily:
\begin{equation*}
\begin{split}
     &\eps^{-1}\left\|\Nabla{t}\xi-\eta\right\|_p
      + \left\|\Nabla{t}\eta\right\|_p
      + \left\|\Nabla{s}\xi\right\|_p
      + \eps\left\|\Nabla{s}\eta\right\|_p \\
     &\le 
      c'\Bigl(
      \left\|f-R(\xi,\p_tu)v-\Hh_\Vv (u)\xi\right\|_p
      +\eps\left\|g+R(\xi,\p_su)v\right\|_p 
      + \left\|\xi\right\|_p + \eps^2\left\|\eta\right\|_p
      \Bigr).
\end{split}
\end{equation*}
To prove the estimate for the formal adjoint 
operator $(\Dd^\eps_{u,v})^*$ apply~(\ref{eq:main-elliptic})
to the vector fields $\xi(-s,t)$ and $\eta(-s,t)$
and then proceed as above. 

To prove~(\ref{eq:main-elliptic}) we
split $\zeta$ into two components. Let
$$
     \pi_\eps(\xi,\eta) 
     := (\1-\eps\Nabla{t}\Nabla{t})^{-1}
       (\xi-\eps^2\Nabla{t}\eta),\qquad
     \i(\xi) := (\xi,\Nabla{t}\xi),
$$
and define 
$$
  \zeta_0 
  := \begin{pmatrix}\xi_0\\ \eta_0\end{pmatrix}
  := \iota \pi_\eps \zeta
  =\begin{pmatrix}(\1-\eps\Nabla{t}\Nabla{t})^{-1}
  (\xi-\eps^2\Nabla{t}\eta) \\
  \Nabla{t}(\1-\eps\Nabla{t}\Nabla{t})^{-1}
  (\xi-\eps^2\Nabla{t}\eta)\end{pmatrix},
$$
$$
     \zeta_1
     := \begin{pmatrix} \xi_1 \\ \eta_1 \end{pmatrix}
     :=\zeta-\zeta_0
     =\begin{pmatrix}(\1-\eps\Nabla{t}\Nabla{t})^{-1}
     (\eps^2\Nabla{t}\eta-\eps\Nabla{t}\Nabla{t}\xi) \\
     (\1-\eps\Nabla{t}\Nabla{t})^{-1}
     (\eta-\Nabla{t}\xi+(\eps^2-\eps)\Nabla{t}\Nabla{t}\eta)
     \end{pmatrix}.
$$
Note that $\eta_0=\Nabla{t}\xi_0$ and 
\begin{equation} \label{eq:difference}
     \xi_1-\eps\Nabla{t}\eta_1=(\eps^2-\eps)\Nabla{t}\eta.
\end{equation}
Since $f$ and $g$ are linear, we obtain
the splitting $f=f_0+f_1$ and $g=g_0+g_1$, where
$f_i:=f(\xi_i,\eta_i)$ and $g_i:=g(\xi_i,\eta_i)$
for $i=0,1$.  Thus 
\begin{equation*}
     f_0 =\Nabla{s}\xi_0-\Nabla{t}\Nabla{t}\xi_0, \qquad 
     g_0 =\Nabla{s}\Nabla{t}\xi_0.
\end{equation*}
Now apply the parabolic estimate of
Proposition~\ref{prop:parabolic-linear}, with a constant $c_0>0$,
to $\xi_0$ and the elliptic estimate of
Proposition~\ref{prop:elliptic-linear}, with a constant $c_1>0$,
to $(\xi_1,\eta_1)$.  This gives
\begin{equation} \label{eq:ref-est}
\begin{split}
    &\eps^{-1}\left\|\Nabla{t}\xi-\eta\right\|_p 
     +\left\|\Nabla{t}\eta\right\|_p+\left\|\Nabla{s}\xi\right\|_p
     +\eps\left\|\Nabla{s}\eta\right\|_p \\
    &\le\left\|\Nabla{t}\Nabla{t}\xi_0\right\|_p
     +\left\|\Nabla{s}\xi_0\right\|_p
     +\eps\left\|\Nabla{s}\Nabla{t}\xi_0\right\|_p \\
    &\quad+\eps^{-1}\left\|\Nabla{t}\xi_1-\eta_1\right\|_p
     +\left\|\Nabla{t}\eta_1\right\|_p
     +\left\|\Nabla{s}\xi_1\right\|_p
     +\eps\left\|\Nabla{s}\eta_1\right\|_p \\
    &\le c_0\bigl(\left\|f_0\right\|_p
     +\left\|\xi_0\right\|_p\bigr)+\eps\left\|g_0\right\|_p \\
    &\quad +c_1\bigl(\left\|f_1\right\|_p+\eps\left\|g_1\right\|_p
     +\eps^{-1}\left\|\xi_1\right\|_p
     +\left\|\eta_1\right\|_p\bigr) \\
    &\le c_1\bigl(\left\|f\right\|_p
     +\eps\left\|g\right\|_p+\eps^{-1}\left\|\xi_1\right\|_p
     +\left\|\eta_1\right\|_p \bigr) \\
    &\quad +(c_0+c_1)\left\|f_0\right\|_p
     +(1+c_1)\eps\left\|g_0\right\|_p 
     +c_0\left\|\xi_0\right\|_p.
\end{split}
\end{equation}
We examine the last five terms on the right individually. 
For this we shall need the commutator identities
\begin{align}
\begin{split}\label{eq:comm1} 
     [\Nabla{s},\Nabla{t}]
    &=R(\p_su,\p_tu),
\end{split}\\
\begin{split}\label{eq:comm2}
     {[\Nabla{s},\Nabla{t}\Nabla{t}]}
    &=2\Nabla{t}[\Nabla{s},\Nabla{t}]
      - (\Nabla{\p_tu}R)(\p_su,\p_tu) \\
    &\quad-\,R(\Nabla{t}\p_su,\p_tu)+R(\p_su,\Nabla{t}\p_tu),
\end{split}\\
\begin{split}\label{eq:comm3}
    {[\Nabla{s},(\1-\eps\Nabla{t}\Nabla{t})^{-1}]}
    &=(\1-\eps\Nabla{t}\Nabla{t})^{-1} 
     {[\1-\eps\Nabla{t}\Nabla{t},\Nabla{s}]}
     (\1-\eps\Nabla{t}\Nabla{t})^{-1} \\
    &=\eps(\1-\eps\Nabla{t}\Nabla{t})^{-1} 
     {[\Nabla{s},\Nabla{t}\Nabla{t}]}
     (\1-\eps\Nabla{t}\Nabla{t})^{-1}. 
\end{split} 
\end{align}
By Lemma~\ref{le:eat-eps} and~(\ref{eq:comm2}), we have
\begin{equation}\label{eq:comm4}
\begin{split}
    &\eps^{1/2}\left\|(\1-\eps\Nabla{t}\Nabla{t})^{-1}
     [\Nabla{s},\Nabla{t}\Nabla{t}]\xi\right\|_p \\
    &\le
     2\eps^{1/2}\left\|(\1-\eps\Nabla{t}\Nabla{t})^{-1}
     \Nabla{t}[\Nabla{s},\Nabla{t}]\xi\right\|_p
     + c_1\eps^{1/2}\left\|\xi\right\|_p \\
    &\le
     2\kappa_p\left\|[\Nabla{s},\Nabla{t}]\xi\right\|_p
     + c_1\eps^{1/2}\left\|\xi\right\|_p \\
     &\le c_2\left\|\xi\right\|_p.
\end{split}
\end{equation}
Here we have used the $L^\infty$ bounds on 
$\p_su$, $\p_tu$, $\Nabla{t}\p_tu$, and $\Nabla{t}\p_su$. 
Now the five relevant terms are estimated as follows.

\smallskip
\noindent
{\bf The term \boldmath$\left\| \xi_0 \right\|_p$:} 
By definition, 
$$
     \xi_0=(\1-\eps\Nabla{t}\Nabla{t})^{-1}(\xi-\eps^2\Nabla{t}\eta).
$$ 
Hence, by Lemma~\ref{le:eat-eps},
\begin{equation} \label{eq:xi_0}
     \left\|\xi_0\right\|_p
     \le\left\|\xi\right\|_p+\eps^2\left\|\Nabla{t}\eta\right\|_p.
\end{equation}

\smallskip
\noindent
{\bf The term \boldmath$\| f_0 \|_p$:} 
Consider the identity
\begin{equation*}
\begin{split}
    &(\1-\eps\Nabla{t}\Nabla{t})f_0 
     -f+\eps^2\Nabla{t}g \\
    &=\Nabla{s}\xi_0-\eps\Nabla{t}\Nabla{t}\Nabla{s}\xi_0
     -\Nabla{t}\Nabla{t}\xi_0
     +\eps\Nabla{t}\Nabla{t}\Nabla{t}\Nabla{t}\xi_0
     -\Nabla{s}\xi
     +\eps^2\Nabla{t}\Nabla{s}\eta
     +\Nabla{t}\Nabla{t}\xi \\
    &=\eps^2\Nabla{t}\Nabla{t}\Nabla{t}\eta
     +\eps^2R(\p_tu,\p_su)\eta
     + \eps [\Nabla{s},\Nabla{t}\Nabla{t}] \xi_0. \\
\end{split}
\end{equation*}
Apply the operator $(\1-\eps\Nabla{t}\Nabla{t})^{-1}$
to this equation and use Lemma~\ref{le:eat-eps} 
and~(\ref{eq:comm4}) to obtain
\begin{equation}\label{eq:f-0}
     \left\|f_0\right\|_p
     \le\left\|f\right\|_p+\kappa_p\eps^{3/2}\left\|g\right\|_p
     +2\eps\left\|\Nabla{t}\eta\right\|_p 
     +\eps^2c_3\left\|\eta\right\|_p
     +\eps^{1/2}c_2\left\|\xi_0\right\|_p,
\end{equation}
where $c_3:=\left\|R\right\|_\infty\left\|\p_su\right\|_\infty
\left\|\p_tu\right\|_\infty$. 

\smallskip
\noindent
{\bf The term \boldmath$\eps \| g_0 \|_p$:}
By~(\ref{eq:comm3}), we have
\begin{equation*}
\begin{split}
     g_0
    &=\Nabla{s}\Nabla{t}\xi_0 \\
    &=(\1-\eps\Nabla{t}\Nabla{t})^{-1}
     \left(\Nabla{t}\Nabla{s}\xi+[\Nabla{s},\Nabla{t}]\xi
     -\eps^2\Nabla{t}\Nabla{t}\Nabla{s}\eta
     -\eps^2[\Nabla{s},\Nabla{t}\Nabla{t}]\eta\right) \\
   &\quad +\eps(\1-\eps\Nabla{t}\Nabla{t})^{-1}
     [\Nabla{s},\Nabla{t}\Nabla{t}]
     (\1-\eps\Nabla{t}\Nabla{t})^{-1}
     \left(\Nabla{t}\xi-\eps^2\Nabla{t}\Nabla{t}\eta\right).
\end{split}
\end{equation*}
Hence, by Lemma~\ref{le:eat-eps}, (\ref{eq:comm1}), 
and~(\ref{eq:comm4}),
\begin{equation}\label{eq:g-0}
\begin{split}
     \eps\left\|g_0\right\|_p
    &\le \kappa_p\eps^{1/2}\left\|\Nabla{s}\xi\right\|_p
     +c_3\eps\left\|\xi\right\|_p
     +2\eps^2\left\|\Nabla{s}\eta\right\|_p 
     +c_2\eps^{5/2}\left\|\eta\right\|_p \\
    &\quad 
     +c_2\eps^{3/2}\left\|(\1-\eps\Nabla{t}\Nabla{t})^{-1}
     \left(\Nabla{t}\xi-\eps^2\Nabla{t}\Nabla{t}\eta\right)\right\|_p \\
    &\le \kappa_p\eps^{1/2}\left\|\Nabla{s}\xi\right\|_p
     +2\eps^2\left\|\Nabla{s}\eta\right\|_p \\
    &\quad
     +\eps(\kappa_pc_2+c_3)\left\|\xi\right\|_p
     +3c_2\eps^{5/2}\left\|\eta\right\|_p.
\end{split}
\end{equation}

\smallskip
\noindent
{\bf The term \boldmath$\eps^{-1} \| \xi_1 \|_p$:}
By~(\ref{eq:difference}), we have 
$$
     \eps^{-1}\xi_1
     = \Nabla{t}\eta_1+\eps\Nabla{t}\eta-\Nabla{t}\eta
     = \eps\Nabla{t}\eta-\Nabla{t}\eta_0.  
$$
Hence
\begin{equation}\label{eq:eps-xi_1}
\begin{split}
     \eps^{-1}\|\xi_1\|_p
    &\le\eps\left\|\Nabla{t}\eta\right\|_p
     +\left\|\Nabla{t}\Nabla{t}\xi_0\right\|_p \\
    &\le\eps\left\|\Nabla{t}\eta\right\|_p
     +c_0\left(\left\|f_0\right\|_p
     +\left\|\xi_0\right\|_p\right).
\end{split}
\end{equation}
In the last step we have used the parabolic estimate
of Proposition~\ref{prop:parabolic-linear}.

\smallskip
\noindent
{\bf The term \boldmath$\| \eta_1 \|_p$:} 
By definition, 
$$
     \eta_1=\left(\1-\eps\Nabla{t}\Nabla{t}\right)^{-1}
     \bigl(\eta-\Nabla{t}\xi+(\eps^2-\eps)\Nabla{t}\Nabla{t}\eta
     \bigr).
$$
Hence, by the triangle inequality and 
Lemma~\ref{le:eat-eps}, we have 
\begin{equation}\label{eq:eta-1}
\begin{split}
\Norm{\eta_1}_p
&\le \Norm{\left(\1-\eps\Nabla{t}\Nabla{t}\right)^{-1}
(\eta-\Nabla{t}\xi)}_p 
+ \eps\Norm{\left(\1-\eps\Nabla{t}\Nabla{t}\right)^{-1}
\Nabla{t}\Nabla{t}\eta}_p \\
&\le
\Norm{\eta-\Nabla{t}\xi}_p 
+ \kappa_p\sqrt{\eps}\Norm{\Nabla{t}\eta}_p.
\end{split}
\end{equation}
Insert the five estimates~(\ref{eq:xi_0}-\ref{eq:eta-1})
into~(\ref{eq:ref-est}) to obtain~(\ref{eq:elliptic-eps}),
provided that $\eps$ is sufficiently small. 
This proves Theorem~\ref{thm:elliptic-eps}.
\end{proof} 

\subsection*{The estimate for the inverse}

Geometrically, the difference between the operators $\Dd_u^0$ 
and $\Dd_{u,v}^\eps$ is the difference between
configuration space and phase space, or between loops in 
$M$ and loops in $T^*M\cong TM$.  Consider the embedding 
$$
     \Ll M\to \Ll TM:x\mapsto(x,\dot x).
$$
The differential of this embedding is given by 
$$
     \Omega^0(S^1,x^*TM)
     \to \Omega^0(S^1,x^*TM\oplus x^*TM):
     \xi\mapsto(\xi,\Nabla{t}\xi).
$$
To compare the operators $\Dd_u^0$ and
$\Dd_u^\eps:=\Dd_{u,\p_tu}^\eps$
we must choose a projection onto the image of 
this embedding (along $u$). 
At first glance it might seem natural 
to choose the orthogonal projection with respect to 
the inner product determined by the 
$(0,2,\eps)$-Hilbert space structure.  
This is given by 
$$
     (\xi,\eta)
     \mapsto
     (\1-\eps^\alpha\Nabla{t}\Nabla{t})^{-1}
     (\xi-\eps^\beta\Nabla{t}\eta)
$$
with $\alpha=\beta=2$. 
Instead we introduce the projection operator
$$
     \pi_\eps:L^p(S^1,u^*TM)\times L^p(S^1,u^*TM)\to 
              W^{1,p}(S^1,u^*TM)
$$
given by 
\begin{equation}\label{eq:projection}
     \pi_\eps(\xi,\eta) 
     := (\1-\eps\Nabla{t}\Nabla{t})^{-1}
       (\xi-\eps^2\Nabla{t}\eta).
\end{equation}
The reason for this choice
becomes visible in the proof of 
Proposition~\ref{prop:4.2p>2} below,
which requires $\beta=2$. Moreover, 
the estimates in Step~1
of the proof of Theorem~\ref{thm:inverse}
are optimized for $\alpha=1$.
We denote by 
$
     \i:W^{1,p}(\R\times S^1,u^*TM)\to
     L^p(S^1,u^*TM)\times L^p(S^1,u^*TM)
$
the inclusion
\begin{equation}\label{eq:inclusion}
     \i\xi_0 := (\xi_0,\Nabla{t}\xi_0).
\end{equation}
The significance of these definitions lies
in the next proposition and lemma.
The proofs rely on Lemma~\ref{le:eat-eps}.

\begin{proposition}\label{prop:4.2p>2}
Let $u\in\Cinf(\R \times S^1,M)$
be a smooth map such that the derivatives
$\p_su,\p_tu,\Nabla{t}\p_su,\Nabla{t}\p_tu,
\Nabla{t}\Nabla{t}\p_tu$ are bounded
and define $v:=\p_t u$.
Then, for every $p>1$,
there exists a constant $c>0$ such that
$$
     \left\|\Dd_u^0\pi_\eps\zeta
     -\pi_\eps\Dd_u^\eps\zeta\right\|_p
     \le c\eps^{1/2}\left\|\xi\right\|_p
     +c\eps^2\left\|\eta\right\|_p
     +c\eps\left\|\Nabla{t}\eta\right\|_p
$$
for $\eps\in(0,1]$ and
compactly supported $\zeta=(\xi,\eta)\in 
\Omega^0(\R \times S^1,u^*TM\oplus u^*TM)$.
The same estimate holds for
$(\Dd_u^0)^*\pi_\eps-\pi_\eps(\Dd_u^\eps)^*$.
Moreover, the constant $c$
is invariant under $s$-shifts of $u$.
\end{proposition}

\begin{lemma} \label{le:4.3p>2}
For $u\in\Cinf(\R \times S^1,M)$, $p>1$,
$\kappa_p$ as in~(\ref{eq:kappa-p}), 
and $0<\eps\le1$,
\begin{equation*}
\begin{split}
     \left\|\xi-\pi_\eps\zeta\right\|_p
     &\le \kappa_p\eps^{1/2}
      \left\|\Nabla{t}\xi-\eta\right\|_p
      +\eps\left\|\Nabla{t}\eta\right\|_p \\
      \left\|\eta-\Nabla{t}\pi_\eps\zeta\right\|_p
     &\le\left\|\Nabla{t}\xi-\eta\right\|_p
      +\kappa_p\eps^{1/2}
      \left\|\Nabla{t}\eta\right\|_p \\
      \left\|\zeta-\iota\pi_\eps\zeta\right\|_{0,p,\eps}
     &\le 2\kappa_p\eps^{1/2}
      \left\|\Nabla{t}\xi-\eta\right\|_p
      +2\kappa_p\eps\left\|\Nabla{t}\eta\right\|_p \\
      \left\|\pi_\eps\zeta\right\|_p
     &\le\left\|\iota\pi_\eps\zeta\right\|_{0,p,\eps}
      \le 2\kappa_p\left\|\zeta\right\|_{0,p,\eps} \\
\end{split}
\end{equation*}
for every compactly supported
$\zeta =(\xi,\eta)\in \Omega^0(\R \times S^1,u^*TM\oplus u^*TM)$.
\end{lemma}

\begin{proof}
Denote
$$
     \xi_0 
     := \pi_\eps\zeta 
     =(\1-\eps\Nabla{t}\Nabla{t})^{-1}
      (\xi-\eps^2\Nabla{t}\eta).
$$
Then
$$
     \xi-\xi_0
     = \eps(\1-\eps\Nabla{t}\Nabla{t})^{-1}
       \Nabla{t}(\eta-\Nabla{t}\xi)
       +(\eps^2-\eps)(\1-\eps\Nabla{t}\Nabla{t})^{-1}
       \Nabla{t}\eta
$$
and hence, by Lemma~\ref{le:eat-eps}, 
$$
     \left\|\xi-\xi_0\right\|_p
     \le \kappa_p\eps^{1/2}\left\|\Nabla{t}\xi-\eta\right\|_p
        +\eps\left\|\Nabla{t}\eta\right\|_p
$$
Similarly, 
$$
     \eta-\Nabla{t}\xi_0
     =(\1-\eps\Nabla{t}\Nabla{t})^{-1}(\eta-\Nabla{t}\xi)
      +(\eps^2-\eps)(\1-\eps\Nabla{t}\Nabla{t})^{-1}
      \Nabla{t}\Nabla{t}\eta
$$
and hence, again by Lemma~\ref{le:eat-eps},
$$
     \eps\left\|\eta-\Nabla{t}\xi_0\right\|_p
     \le \eps\left\|\Nabla{t}\xi-\eta\right\|_p
        +\kappa_p\eps^{3/2}\left\|\Nabla{t}\eta\right\|_p.
$$
Take the sum of these two inequalities to obtain
\begin{equation*}
\begin{split}
     \left\|\zeta-\i\pi_\eps\zeta\right\|_{0,p,\eps}
&\le
     \left\|\xi-\xi_0\right\|_p
        +\eps\left\|\eta-\Nabla{t}\xi_0\right\|_p  \\
&\le
      2\kappa_p\eps^{1/2}\left\|\Nabla{t}\xi-\eta\right\|_p
        +2\kappa_p\eps\left\|\Nabla{t}\eta\right\|_p
\end{split}
\end{equation*}
for $0<\eps\le 1$.  Moreover, using Lemma~\ref{le:eat-eps} 
the formula for $\xi_0$ gives
$$
     \left\|\xi_0\right\|_p
     \le\left\|\xi\right\|_p
        +\kappa_p\eps^{3/2}\left\|\eta\right\|_p,\qquad
        \eps\left\|\Nabla{t}\xi_0\right\|_p
     \le \kappa_p\eps^{1/2}\left\|\xi\right\|_p
         +2\eps^2\left\|\eta\right\|_p.
$$
Take these two inequalities to the power $p$
and take the sum to obtain
\begin{equation*}
\begin{split}
     \left\|\i\pi_\eps\zeta\right\|_{0,p,\eps}^p
    &=\left\|\xi_0\right\|_p^p
     +\eps^p\left\|\Nabla{t}\xi_0\right\|_p^p \\
    &\le (1+\kappa_p^p\eps^{p/2})
     \left\|\xi\right\|_p^p
     +(\kappa_p^p\eps^{p/2}+2^p\eps^p)
     \eps^p\left\|\eta\right\|_p^p \\
    &\le (2\kappa_p)^p\left\|\zeta\right\|_{0,p,\eps}^p
\end{split}
\end{equation*}
for $0<\eps\le 1$. This proves Lemma~\ref{le:4.3p>2}.
\end{proof} 

\begin{proof}[Proof of Proposition~\ref{prop:4.2p>2}.]
As above, denote
$$
     \xi_0:=\pi_\eps\zeta 
     = (\1-\eps\Nabla{t}\Nabla{t})^{-1}
       (\xi-\eps^2\Nabla{t}\eta).
$$
Then
\begin{equation*}
\begin{split}
     \Dd_u^0\pi_\eps\zeta
    &=\Nabla{s}\xi_0-\Nabla{t}\Nabla{t}\xi_0
     -R(\xi_0,\p_tu)\p_tu-\Hh_\Vv(u)\xi_0 \\
    &=(\1-\eps\Nabla{t}\Nabla{t})^{-1}\left(
     \Nabla{s}\xi-\eps^2\Nabla{s}\Nabla{t}\eta
     -\Nabla{t}\Nabla{t}\xi
     +\eps^2\Nabla{t}\Nabla{t}\Nabla{t}\eta\right) \\
    &\quad+\eps(\1-\eps\Nabla{t}\Nabla{t})^{-1}
     [\Nabla{s},\Nabla{t}\Nabla{t}] \xi_0 \\
    &\quad+R\bigl((\1-\eps\Nabla{t}\Nabla{t})^{-1}
     \eps^2\Nabla{t}\eta\:,\:\p_tu\bigr)\p_tu
     +\Hh_\Vv(u)(\1-\eps\Nabla{t}\Nabla{t})^{-1}
     \eps^2\Nabla{t}\eta \\
    &\quad-R\bigl((\1-\eps\Nabla{t}\Nabla{t})^{-1}\xi\:,\:
     \p_tu\bigr)\p_tu
     -\Hh_\Vv(u)(\1-\eps\Nabla{t}\Nabla{t})^{-1}\xi.
\end{split}
\end{equation*}
Denote $\zeta':=(\xi',\eta')
:=\Dd_u^\eps\zeta$, then
\begin{equation*}
\begin{split}
     \pi_\eps\Dd_u^\eps\zeta
    &=(\1-\eps\Nabla{t}\Nabla{t})^{-1}(\xi'
     -\eps^2\Nabla{t}\eta') \\
    &=(\1-\eps\Nabla{t}\Nabla{t})^{-1}\bigl(
     \Nabla{s}\xi-R(\xi,\p_tu)\p_tu-\Hh_\Vv(u)\xi \\
    &\quad-\eps^2\Nabla{t}\Nabla{s}\eta
     -\eps^2\Nabla{t}\bigl(R(\xi,\p_su)\p_tu\bigr)
     -\Nabla{t}\Nabla{t}\xi\bigr).
\end{split}
\end{equation*}
Taking the difference we find
\begin{equation}\label{eq:op-diff}
\begin{split}
    &\Dd_u^0\pi_\eps\zeta
     -\pi_\eps\Dd_u^\eps\zeta \\
    &=(\1-\eps\Nabla{t}\Nabla{t})^{-1}\left(
     -\eps^2[\Nabla{s},\Nabla{t}]\eta
     +\eps^2\Nabla{t}\Nabla{t}\Nabla{t}\eta
     +\eps^2\Nabla{t}\bigl(R(\xi,\p_su)\p_tu\bigr)
     \right) \\
    &\quad+\eps(\1-\eps\Nabla{t}\Nabla{t})^{-1}
     [\Nabla{s},\Nabla{t}\Nabla{t}] \xi_0 \\
    &\quad+R\bigl((\1-\eps\Nabla{t}\Nabla{t})^{-1}
     \eps^2\Nabla{t}\eta\:,\:\p_tu\bigr)\p_tu
     +\Hh_\Vv(u)(\1-\eps\Nabla{t}\Nabla{t})^{-1}
     \eps^2\Nabla{t}\eta \\
    &\quad+(\1-\eps\Nabla{t}\Nabla{t})^{-1}R(\xi,\p_tu)\p_tu
     -R\bigl((\1-\eps\Nabla{t}\Nabla{t})^{-1}\xi\:,\:
     \p_tu\bigr)\p_tu \\
    &\quad+(\1-\eps\Nabla{t}\Nabla{t})^{-1}\Hh_\Vv(u)\xi
     -\Hh_\Vv(u)(\1-\eps\Nabla{t}\Nabla{t})^{-1}\xi.
\end{split}
\end{equation}
To finish the proof it remains to inspect 
the $L^p$ norm of this expression line by line.
Using Lemma~\ref{le:eat-eps}, we obtain
for the first line
\begin{equation}\label{eq:line1}
\begin{split}
    &\left\|
     (\1-\eps\Nabla{t}\Nabla{t})^{-1}\left(
     -\eps^2[\Nabla{s},\Nabla{t}]\eta
     +\eps^2\Nabla{t}\Nabla{t}\Nabla{t}\eta
     +\eps^2\Nabla{t}\bigl(R(\xi,\p_su)\p_tu\bigr)
     \right)
     \right\|_p \\
    &\le\eps^2\left\|R\right\|_\infty
     \left\|\p_su\right\|_\infty
     \left\|\p_tu\right\|_\infty\left\|\eta\right\|_p
     +2\eps\left\|\Nabla{t}\eta\right\|_p \\
    &\quad
     +\kappa_p\eps^{3/2}\left\|R\right\|_\infty
     \left\|\p_su\right\|_\infty
     \left\|\p_tu\right\|_\infty\left\|\xi\right\|_p.
\end{split}
\end{equation}
Application of $(\ref{eq:comm4})$ with
constant $C_1:=C$ results in an estimate 
for the second line in $(\ref{eq:op-diff})$, namely
\begin{equation}\label{eq:line2}
     \left\|
     \eps(\1-\eps\Nabla{t}\Nabla{t})^{-1}
     [\Nabla{s},\Nabla{t}\Nabla{t}] \xi_0
     \right\|_p
   \le\eps^{1/2}C_1\left\|\xi\right\|_p
     +\eps^{5/2}C_1\left\|\Nabla{t}\eta\right\|_p.
\end{equation}
Lemma~\ref{le:eat-eps} yields
for line three in~(\ref{eq:op-diff})
\begin{equation}\label{eq:line3}
\begin{split}
    &\left\|
     R\bigl((\1-\eps\Nabla{t}\Nabla{t})^{-1}
     \eps^2\Nabla{t}\eta\:,\:\p_tu\bigr)\p_tu
     +\Hh_\Vv(u)(\1-\eps\Nabla{t}\Nabla{t})^{-1}
     \eps^2\Nabla{t}\eta
     \right\|_p \\
    &\le\bigl(\left\|R\right\|_\infty
     \left\|\p_tu\right\|_\infty^2 + C\bigr)
     \eps^2\left\|\Nabla{t}\eta\right\|_p,
\end{split}
\end{equation}
where $C$ is the constant in~$(V1)$. 
Let us temporarily denote 
$$
    T:=\1-\eps\Nabla{t}\Nabla{t}.
$$
Then the penultimate line in~(\ref{eq:op-diff}) has the form
$
    [T^{-1},\Phi]
    = T^{-1}[\Phi,T]T^{-1}
$
where the endomorphism $\Phi:u^*TM\to u^*TM$
is given by $\Phi\xi=R(\xi,\p_tu)\p_tu$.
This term can be expressed in the form
$$
    [T^{-1},\Phi]\xi
    = \eps T^{-1}
      \left((\Nabla{t}\Nabla{t}\Phi)T^{-1}\xi
      + 2(\Nabla{t}\Phi)T^{-1}\Nabla{t}\xi\right)
$$
and hence
$$
    \Norm{[T^{-1},\Phi]\xi}_p
    \le \eps^{1/2}\kappa_pC\Norm{\xi}_p.
$$
Thus 
\begin{equation}\label{eq:line4}
\begin{split}
     &\left\|
     (\1-\eps\Nabla{t}\Nabla{t})^{-1}R(\xi,\p_tu)\p_tu
     -R\bigl((\1-\eps\Nabla{t}\Nabla{t})^{-1}
     \xi,\p_tu\bigr)\p_tu
     \right\|_p  \\
     &\le \eps^{1/2}\kappa_pC_2\left\|\xi\right\|_p,
\end{split}
\end{equation}
where $C_2$ depends on $\Norm{R}_{C^2}$
$\left\|\p_tu\right\|_\infty$,
$\left\|\Nabla{t}\p_tu\right\|_\infty$, and
$\left\|\Nabla{t}\Nabla{t}\p_tu\right\|_\infty$.
Similarly,
\begin{equation}\label{eq:line5}
    \left\|
     (\1-\eps\Nabla{t}\Nabla{t})^{-1}\Hh_\Vv(u)\xi
     -\Hh_\Vv(u)(\1-\eps\Nabla{t}\Nabla{t})^{-1}\xi
     \right\|_p
     \le \eps^{1/2}\kappa_pC_3\left\|\xi\right\|_p,
\end{equation}
where $C_3$ depends on the constants in~$(V1-V3)$
and on $\left\|\p_tu\right\|_\infty$
and $\left\|\Nabla{t}\p_tu\right\|_\infty$.
The estimates~(\ref{eq:line1}-\ref{eq:line5})
together give the desired $L^p$ bound for~(\ref{eq:op-diff})
and this proves the first claim of Proposition~\ref{prop:4.2p>2}.
The estimate for 
$(\Dd_u^0)^*\pi_\eps\zeta-(\pi_\eps\Dd_u^\eps)^*\zeta$
follows analoguously. Since all constants appearing
in the proof depend on $L^\infty$ norms of derivatives
of $u$, they are invariant under $s$-shifts of $u$.
This completes the proof of Proposition~\ref{prop:4.2p>2}.
\end{proof} 

The next lemma establishes the relevant estimates
for the operator $\Dd^0_u$ and its adjoint in the 
Morse--Smale case, i.e. when $\Dd^0_u$ is onto. 

\begin{lemma} \label{le:step3}
Let $\Vv:\Ll M\to\R$ be a perturbation that 
satisfies~$(V0-V4)$.  Assume $\Ss_\Vv$ is Morse-Smale 
and let $u\in \Mm^0(x^-,x^+;\Vv)$. Then, for every $p>1$, 
there is a constant $c>0$ such that
$$
     \left\|\eta\right\|_p
     +\left\|\Nabla{s}\eta\right\|_p
     + \left\|\Nabla{t}\Nabla{t}\eta\right\|_p
     \le c\left\|(\Dd_u^0)^*\eta\right\|_p
$$
and
$$
     \left\|\xi\right\|_p
     +\left\|\Nabla{s}\xi\right\|_p
     + \left\|\Nabla{t}\Nabla{t}\xi\right\|_p
     \le c\left( 
     \left\|\xi-(\Dd_u^0)^*\eta\right\|_p
     +\left\|\Dd_u^0\xi\right\|_p
     \right)
$$
for all compactly supported vector fields
$\xi,\eta\in\Omega^0(\R \times S^1,u^*TM)$.
\end{lemma}

\begin{proof}
By Theorem~\ref{thm:par-Fredholm},
the operators $\Dd_u^0$ and $(\Dd_u^0)^*$ 
are Fredholm.  Since $\Ss_\Vv$ is Morse--Smale, 
the operator $\Dd^0_u$ is onto and 
$(\Dd_u^0)^*$ is injective. Moreover, the 
operator 
$$
     \Ww^p_u\to\Ll^p_u\oplus\Ll^p_u/\im\,(\Dd^0_u)^*:
     \xi\mapsto(\Dd^0_u\xi,[\xi])
$$
is also an injective Fredholm operator.
Hence the estimates follow from the open mapping theorem. 
\end{proof}

\begin{proof}[Proof of Theorem~\ref{thm:inverse}]
Fix a constant $p>1$. Then the $L^\infty$ norms of
$\p_su$, $\p_tu$ and $\Nabla{t}\p_tu$ are finite
by Theorem~\ref{thm:par-apriori} and
$\|\Nabla{t}\p_su\|_\infty$ is finite
by Theorem~\ref{thm:par-exp-decay}.
Use the parabolic equations for $u$
to conclude that $\|\Nabla{t}\Nabla{t}\p_tu\|_\infty$
is finite as well.
Hence we are in a position to apply
Theorem~\ref{thm:elliptic-eps} and
Proposition~\ref{prop:4.2p>2}.
We prove the estimate in two steps. 

\medskip
\noindent
{\sc Step~1.} 
{\it There are positive constants $c_1=c_1(p)$
and $\eps_0=\eps_0(p)$ such that
\begin{equation}\label{eq:main-inj-adj}
     \left\|\zeta\right\|_{0,p,\eps}
   \le \left\|\xi\right\|_p
     +\eps^{1/2}\left\|\eta\right\|_p
   \le c_1\left(
     \eps\left\|(\Dd_u^\eps)^*
     \zeta\right\|_{0,p,\eps}
     +\left\|\pi_\eps(\Dd_u^\eps)^*\zeta\right\|_p
     \right)
\end{equation}
for every $\eps\in(0,\eps_0)$
and every compactly supported vector field
$\zeta=(\xi,\eta )\in
\Omega^0(\R\times S^1,u^*TM \oplus u^*TM)$.}

\medskip
\noindent
By Lemmata~\ref{le:nabla-t-xi} and~\ref{le:step3},
there exists a constant $c_2=c_2(p)>0$ such that
\begin{equation}\label{eq:inj-adj}
     \left\|\xi\right\|_p
     +\left\|\Nabla{s}\xi\right\|_p
     +\left\|\Nabla{t}\xi\right\|_p
     +\left\|\Nabla{t}\Nabla{t}\xi\right\|_p
     \le c_2 \left\|(\Dd_u^0)^*\xi\right\|_p
\end{equation}
for every compactly supported
$\xi\in\Omega^0(\R \times S^1,u^*TM)$.
Hence
\begin{equation*}
\begin{split}
     \left\|\xi\right\|_p
    &\le\left\|\xi-\pi_\eps\zeta\right\|_p
     +\left\|\pi_\eps\zeta\right\|_p \\
    &\le\left\|\xi-\pi_\eps\zeta\right\|_p 
     +c_2\left\|(\Dd_u^0)^*\pi_\eps\zeta\right\|_p \\
    &\le\left\|\xi-\pi_\eps\zeta\right\|_p
     +c_2\left\|(\Dd_u^0)^*\pi_\eps\zeta
     -\pi_\eps(\Dd_u^\eps)^*\zeta\right\|_p
     +c_2\left\|\pi_\eps(\Dd_u^\eps)^*
     \zeta\right\|_p \\
    &\le(\kappa_p+c_2c_3) \eps
     \left(\eps^{-1}\left\|\Nabla{t}\xi-\eta\right\|_p
     +\left\|\Nabla{t}\eta\right\|_p\right)
     +c_2\left\|\pi_\eps(\Dd_u^\eps)^*
     \zeta\right\|_p \\
    &\quad+c_2c_3\left(\eps^{1/2}\left\|\xi\right\|_p
     +\eps^2\left\|\eta\right\|_p\right) \\
    &\le(\kappa_p+c_2c_3)c_4\eps
     \left\|(\Dd_u^\eps)^*\zeta\right\|_{0,p,\eps}
     +c_2\left\|\pi_\eps(\Dd_u^\eps)^*
     \zeta\right\|_p \\
    &\quad+(c_2c_3+\kappa_pc_4+c_2c_3c_4)
     \left(\eps^{1/2}\left\|\xi\right\|_p
     +\eps^2\left\|\eta\right\|_p\right)
\end{split}
\end{equation*}
In the fourth step we have used
Lemma~\ref{le:4.3p>2} and Proposition~\ref{prop:4.2p>2} with
a constant $c_3=c_3(p)>0$. The final step follows from
Theorem~\ref{thm:elliptic-eps} for the formal
adjoint operator with a constant $c_4=c_4(p)>0$.
Choose $\eps_0 >0$ so small that
\begin{equation} \label{eq:eps-0-bound}
  (c_2c_3+\kappa_pc_4+c_2c_3c_4){\eps_0}^{1/2}<\frac{1}{2}.
\end{equation}
Then we can incorporate the term $\| \xi \|_p$ 
into the left hand side and obtain
\begin{equation} \label{eq:xi-part}
     \left\|\xi\right\|_p
     \le2(\kappa_p+c_2c_3)c_4\eps
     \left\|(\Dd_u^\eps)^*\zeta\right\|_{0,p,\eps}
     +2c_2\left\|\pi_\eps(\Dd_u^\eps)^*
     \zeta\right\|_p
     +\eps^{3/2}\left\|\eta\right\|_p.
\end{equation}
Similarly,
\begin{equation*}
\begin{split}
     \left\|\eta\right\|_p
    &\le\left\|\eta-\Nabla{t}\pi_\eps\zeta\right\|_p
     +\left\|\Nabla{t}\pi_\eps\zeta\right\|_p \\
    &\le\left\|\eta-\Nabla{t}\pi_\eps\zeta\right\|_p 
     +c_2\left\|(\Dd_u^0)^*\pi_\eps\zeta\right\|_p \\
    &\le(\kappa_p+c_2c_3\eps^{1/2})c_4\eps^{1/2}
     \left\|(\Dd_u^\eps)^*\zeta\right\|_{0,p,\eps}
     +c_2\left\|\pi_\eps(\Dd_u^\eps)^*
     \zeta\right\|_p \\
    &\quad+(c_2c_3+\kappa_pc_4+c_2c_3c_4\eps^{1/2})
     \left(\eps^{1/2}\left\|\xi\right\|_p
     +\eps^2\left\|\eta\right\|_p\right).
\end{split}
\end{equation*}
Use~(\ref{eq:eps-0-bound}) again to obtain
\begin{equation} \label{eq:eta-part}
     \left\|\eta\right\|_p
   \le2(\kappa_p+c_2c_3)c_4\eps^{1/2}
     \left\|(\Dd_u^\eps)^*\zeta\right\|_{0,p,\eps}
     +2c_2\left\|\pi_\eps(\Dd_u^\eps)^*
     \zeta\right\|_p
     +\left\|\xi\right\|_p.
\end{equation}
The assertion of Step~1 now follows from~(\ref{eq:eta-part}) 
and~(\ref{eq:xi-part}).

\vspace{.1cm}
\noindent
{\sc Step 2} 
{\it We prove the theorem.}

\vspace{.1cm}
\noindent
Let $\eps\in(0,\eps_0)$. By~(\ref{eq:elliptic-standard})
for the formal adjoint operator
(with a constant $c_5>0$), we obtain
\begin{equation} \label{eq:inject-D-eps*}
\begin{split}
     \left\|\zeta\right\|_{1,p,\eps}
    &\le c_5\eps^2
     \left\|(\Dd_u^\eps)^*\zeta\right\|_{0,p,\eps}
     +c_5\left\|\zeta\right\|_{0,p,\eps} \\
    &\le c_5(\eps^2+c_1\eps+2\kappa_pc_1)
     \left\|(\Dd_u^\eps)^*\zeta\right\|_{0,p,\eps}
\end{split}
\end{equation}
Here we have also used the estimate~(\ref{eq:main-inj-adj})
of Step~1 and Lemma~\ref{le:4.3p>2}. It follows that 
$(\Dd_u^\eps)^*$ is injective and hence $\Dd_u^\eps$ is onto. 

\smallbreak

Let $\zeta=(\xi,\eta )\in
\Omega^0(\R\times S^1,u^*TM \oplus u^*TM)$
be compactly supported and denote
$$
     \zeta^*:=(\xi^*,\eta^*):=
     (\Dd_u^\eps)^*\zeta.
$$
Recall that $c_6$ is the constant of Lemma~\ref{le:step3}
and $c_3$ is the constant of Proposition~\ref{prop:4.2p>2}.
By Lemma~\ref{le:step3}, with $\xi=\pi_\eps\zeta^*$
and $\eta=\pi_\eps\zeta$, we have
\begin{equation}\label{eq:pi-zeta-*}
\begin{split}
     \left\|\pi_\eps\zeta^*\right\|_p
    &\le c_6
     \left\|\pi_\eps\zeta^*-(\Dd_u^0)^*\pi_\eps\zeta\right\|_p
     +c_6\left\|\Dd_u^0\pi_\eps\zeta^*\right\|_p \\
    &\le c_6\left\|\pi_\eps(\Dd_u^\eps)^*\zeta
     -(\Dd_u^0)^*\pi_\eps\zeta\right\|_p
     +c_6\left\|\Dd_u^0\pi_\eps\zeta^*
     -\pi_\eps\Dd_u^\eps\zeta^*\right\|_p \\
    &\quad
     +c_6\left\|\pi_\eps\Dd_u^\eps\zeta^*\right\|_p \\
    &\le c_3c_6\bigl(\eps^{1/2}\left\|\xi\right\|_p
     +\eps^2\left\|\eta\right\|_p
     +\eps\left\|\Nabla{t}\eta\right\|_p\bigr)
     +c_6\left\|\pi_\eps\Dd_u^\eps\zeta^*\right\|_p \\
    &\quad+c_3c_6\bigl(\eps^{1/2}\left\|\xi^*\right\|_p
     +\eps^2\left\|\eta^*\right\|_p
     +\eps\left\|\Nabla{t}\eta*\right\|_p\bigr) \\
    &\le 2c_3c_6(1+c_4\eps^{1/2})
     \eps^{1/2}\left\|\zeta\right\|_{0,p,\eps}
     +c_6\left\|\pi_\eps\Dd_u^\eps\zeta^*\right\|_p \\
    &\quad
     +3c_3c_6(1+c_4\eps^{1/2})
     \eps^{1/2}\left\|\zeta^*\right\|_{0,p,\eps}
     +c_3c_4c_6\eps
     \left\|\Dd_u^\eps\zeta^*\right\|_{0,p,\eps} \\
    &\le c_7\eps^{1/2}\left\|\zeta^*\right\|_{0,p,\eps}
     +c_3c_4c_6\eps
     \left\|\Dd_u^\eps\zeta^*\right\|_{0,p,\eps}
     +c_6\left\|\pi_\eps\Dd_u^\eps\zeta^*\right\|_p.
\end{split}
\end{equation}
The fourth step follows by applying 
Theorem~\ref{thm:elliptic-eps} twice, with the constant $c_4$,
namely for the operator $(\Dd_u^\eps)^*$
to deal with the term $\Nabla{t}\eta$,
and for the operator $\Dd_u^\eps$
to deal with the term $\Nabla{t}\eta^*$.
The final step follows from~(\ref{eq:inject-D-eps*}).

Now it follows from Lemma~\ref{le:4.3p>2} that
\begin{equation*}
\begin{split}
     \left\|\zeta^*\right\|_{0,p,\eps}
    &\le\left\|\zeta^*-\i\pi_\eps\zeta^*\right\|_{0,p,\eps}
     +\left\|\i\pi_\eps\zeta^*\right\|_{0,p,\eps} \\
    &\le2\kappa_p\eps\left(\eps^{-1}
     \left\|\Nabla{t}\xi^*-\eta^*\right\|_p
     +\left\|\Nabla{t}\eta^*\right\|_p\right)
     +\left\|\pi_\eps\zeta^*\right\|_p
     +\eps\left\|\Nabla{t}\pi_\eps\zeta^*\right\|_p \\
    &\le2\kappa_pc_4\eps
     \left\|\Dd_u^\eps\zeta^*\right\|_{0,p,\eps}
     +\left\|\pi_\eps\zeta^*\right\|_p
     +(2\kappa_p+4\kappa_pc_4) \eps^{1/2}
     \left\|\zeta^*\right\|_{0,p,\eps} \\
    &\le c_4(2\kappa_p+c_3c_6)\eps
     \left\|\Dd_u^\eps\zeta^*\right\|_{0,p,\eps}
     +(c_7+2\kappa_p+4\kappa_pc_4) \eps^{1/2}
     \left\|\zeta^*\right\|_{0,p,\eps} \\
    &\quad
     +c_6\left\|\pi_\eps\Dd_u^\eps\zeta^*\right\|_p.
\end{split}
\end{equation*}
The third step follows from
Theorem~\ref{thm:elliptic-eps}
for the operator $\Dd_u^\eps$
and Lemma~\ref{le:eat-eps}.
The final step uses $(\ref{eq:pi-zeta-*})$.
Choosing $\eps_0>0$ sufficiently
small, we obtain
\begin{equation}\label{eq:zeta-*}
     \left\|\xi^*\right\|_p
   \le\left\|\zeta^*\right\|_{0,p,\eps}
   \le 2c_4(2\kappa_p+c_3c_6)\eps
     \left\|\Dd_u^\eps\zeta^*\right\|_{0,p,\eps}
     +2c_6\left\|\pi_\eps\Dd_u^\eps\zeta^*\right\|_p.
\end{equation}
By~(\ref{eq:elliptic-standard}), we have
\begin{equation*}
     \left\|\zeta^*\right\|_{1,p,\eps} 
   \le c_5\left(\eps^2
     \left\|\Dd_u^\eps\zeta^*\right\|_{0,p,\eps}
     +\left\|\zeta^*\right\|_{0,p,\eps}\right).
\end{equation*}
Combining this with~(\ref{eq:zeta-*}) we 
obtain~(\ref{eq:inverse-zeta}).

We prove~(\ref{eq:inverse-xieta}).
By the triangle inequality and 
Lemmata~\ref{le:4.3p>2} and~\ref{le:eat-eps},
we have
\begin{equation*}
\begin{split}
     \left\|\eta^*\right\|_p
    &\le\left\|\eta^*-\Nabla{t}\pi_\eps\zeta^*\right\|_{0,p,\eps}
     +\left\|\Nabla{t}(\1-\eps\Nabla{t}\Nabla{t})^{-1}
     (\xi^*-\eps^2\Nabla{t}\eta^*)\right\|_p \\
    &\le \kappa_p\eps^{1/2}\left(\eps^{-1}
     \left\|\Nabla{t}\xi^*-\eta^*\right\|_p
     +\left\|\Nabla{t}\eta^*\right\|_p\right)
     +\kappa_p\eps^{-1/2}\left\|\xi^*\right\|_p
     +2\eps\left\|\eta^*\right\|_p \\
    &\le \kappa_pc_4\eps^{1/2}
     \left\|\Dd_u^\eps\zeta^*\right\|_{0,p,\eps}
     +2\kappa_p(1+c_4\eps)\eps^{-1/2}
     \left\|\zeta^*\right\|_{0,p,\eps}
\end{split}
\end{equation*}
The last step follows from
Theorem~\ref{thm:elliptic-eps}
for the operator $\Dd_u^\eps$.
Similarly, 
\begin{equation*}
\begin{split}
     \left\|\Nabla{t}\xi^*\right\|_p
    &\le\left\|\Nabla{t}\xi^*-\eta^*\right\|_p
     +\left\|\eta^*\right\|_p \\
    &\le c_5\eps
     \left\|\Dd_u^\eps\zeta^*\right\|_{0,p,\eps}
     +c_5\eps\left\|\xi^*\right\|_p
     +(1+c_5\eps^3)\left\|\eta^*\right\|_p.
\end{split}
\end{equation*}
Combining the last two estimates with~(\ref{eq:zeta-*})
proves~(\ref{eq:inverse-xieta}).
Since all constants appearing in the proof depend
on $L^\infty$ norms of derivatives of $u$, 
they are invariant under $s$-shifts of $u$.
This proves Theorem~\ref{thm:inverse}.
\end{proof}

\medskip\noindent{\bf Acknowledgement.}
Thanks to Katrin Wehrheim for pointing
out to us the work of Marcinkiewicz and Mihlin
and to Tom Ilmanen for providing the idea for the 
proof of Lemma~\ref{le:apriori-basic-eps}.


\end{document}